\documentclass[11pt]{article}

\usepackage{publication_style}

\hypersetup{pdftitle={Counting Classical Fullerenes I: Exact Enumeration and Generic Achiral Symmetry}}
\title{Counting Classical Fullerenes I:\\Exact Enumeration and Generic Achiral Symmetry}
\author{Kerem Kat\\[2pt]
 \small Independent Researcher\\
 \small\texttt{keremkat@gmail.com}}
\date{}

\begin{document}
\maketitle

\begin{abstract}
We give an exact arithmetic formula for the number of classical fullerenes
on $2m$ vertices. The
oriented count is known; the new contribution counts achiral fullerenes,
those admitting an orientation-reversing symmetry.  We transport these
symmetries to a lattice, classify eleven types, and express their
contribution as seven weighted orbit counts.
The complete formula is printed here, together with an exact
algorithm producing all counts through index $M$ in quasilinear bit time
with fast integer multiplication. We also prove that the
proportion of achiral fullerenes with a unique orientation-reversing
symmetry tends to one as the size increases.  That symmetry is
either a reflection or a fixed-point-free involution, with cumulative
proportions $527/529$ and $2/529$.  Among all fullerenes on $2m$ vertices,
the achiral proportion has order $m^{-5}$; the same order holds cumulatively.
The exceptional achiral classes with nontrivial orientation-preserving
symmetry number $O(m^2)$.
\end{abstract}

\section{Introduction and main results}\label{sec:scope}

The familiar truncated-icosahedral $\mathrm C_{60}$ cage is one of
$1812$ classical fullerene graphs on sixty vertices, up to isomorphism. Of these, $92$
are achiral. The other $1720$ each give two mirror-image classes when
the orientation of the sphere is remembered, so the oriented count is
$2\cdot1720+92=3532$. These numbers, reproduced in
\cref{app:small-counting-examples}, illustrate the distinction at the
heart of this paper: an oriented count determines the ordinary count
once the achiral graphs have been counted.

\begin{figure}[htbp]
\centering
\includegraphics[width=.82\textwidth]{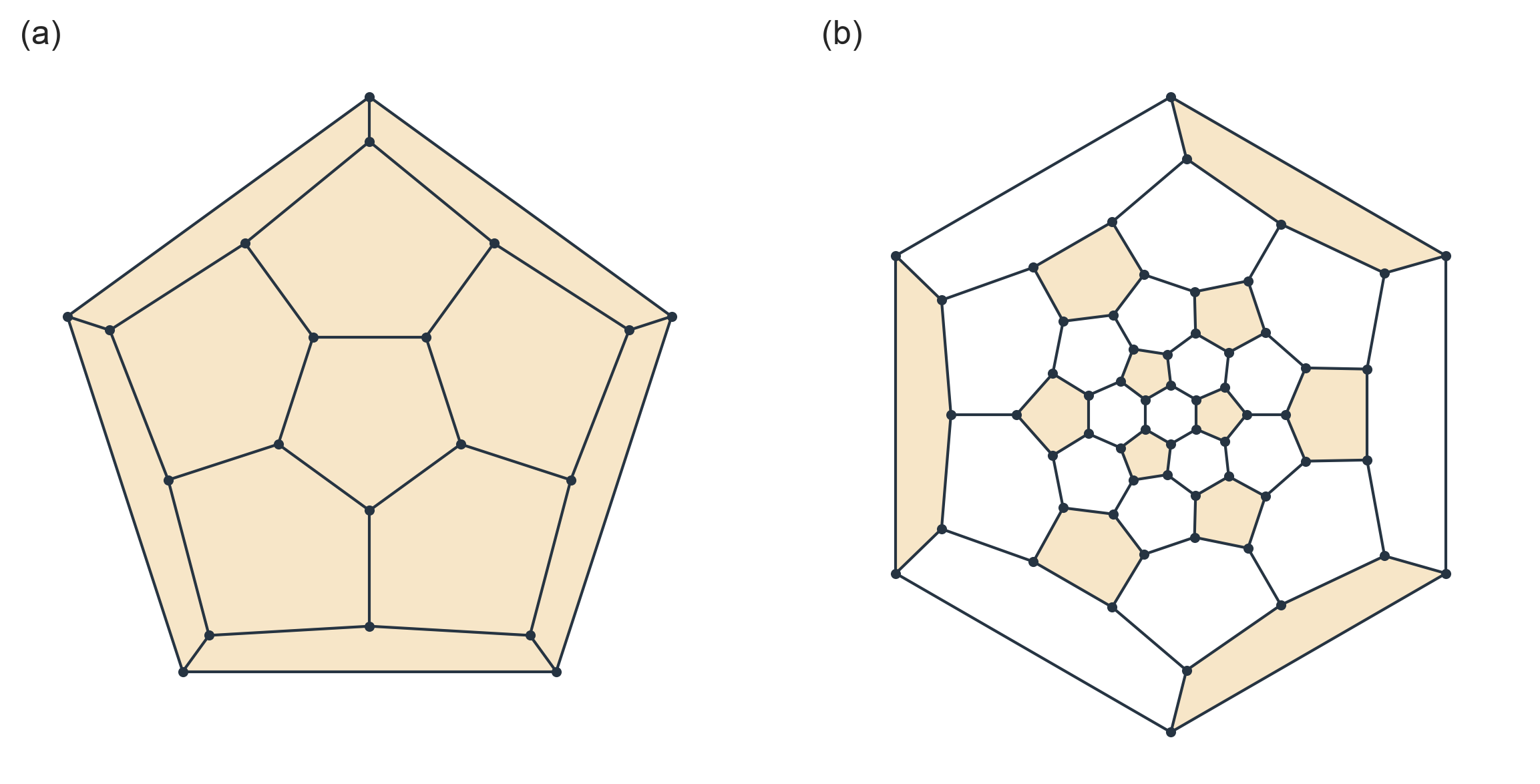}
\caption{Classical fullerene graphs in planar Schlegel projection: (a) the dodecahedral $\mathrm C_{20}$ graph and (b) the truncated-icosahedral $\mathrm C_{60}$ graph. Shading marks the bounded pentagonal faces; the unbounded face is also a face of the cage. The incidence data give respectively $20$ and $60$ cubic vertices, twelve pentagons in each graph, and zero and twenty hexagons. Lengths in the projection carry no molecular metric meaning.}
\label{fig:counting-fullerenes}
\end{figure}

We count carbon frameworks as graphs, independently of their chemical
realization. A
\emph{classical fullerene graph} is a simple, three-connected cubic planar
graph whose faces have sizes five and six.  Euler's relation forces exactly
twelve pentagons.  The degree sum $3|V|=2|E|$ forces $|V|$ to be even,
so we write $|V|=2m$; there are then $m-10$ hexagons.
Figure~\ref{fig:counting-fullerenes} shows the dodecahedral and
truncated-icosahedral examples, with their bounded pentagonal faces shaded. Write
$\iso(2m)$ for the number of ordinary isomorphism classes on $2m$ vertices,
$\BB_m$ for the oriented count, and $a_m$ for the number of achiral classes.
An oriented isomorphism preserves the orientation of the sphere.  An
achiral graph admits an isomorphism reversing that orientation.  Each
chiral graph gives two oriented classes, whereas each achiral graph gives
one.  Thus Burnside's orbit-counting lemma \cite[\S118]{burnside1897} gives
\begin{equation}\label{eq:burnside-count}
  \iso(2m)=\frac{\BB_m+a_m}{2}.
\end{equation}

Engel and Smillie~\cite{engel2025} established an arithmetic formula for
$\BB_m$. We complete the ordinary count by evaluating the classes fixed
by orientation reversal. One graph may have several reversing symmetries,
so counting symmetries requires a weight that counts each graph once.
We therefore mark a reversing symmetry and weight the resulting pair
by the reciprocal of its stabilizer. The weights of all pairs over
one achiral graph add to one.

The dual of a fullerene is a sphere triangulation. Giving its triangles
the unit equilateral metric concentrates all curvature at the twelve
vertices dual to pentagons. Period coordinates encode this metric by
a positive vector in an Eisenstein lattice; a reversing symmetry
becomes an antiunitary operator fixing that vector. This passage
preserves the stabilizers, so the graph count becomes a lattice mass.
Real period geometry leaves eleven operator classes, grouped into
seven arithmetic sectors. Collision mirrors record merged cone
points; weighted inclusion--exclusion removes them.

Marking rotations gives the same kind of calculation for the oriented
count, using five unitary masses.

\begin{theorem}[Exact enumeration]\label{thm:count-overview}
For every positive integer $m$, the ordinary count is
$\iso(2m)=(\BB_m+a_m)/2$, where $\BB_m$ and $a_m$ are the finite
arithmetic expressions in \cref{app:complete-coefficient-formula}.
That appendix gives all $45$ terms of the oriented contribution and all
$181$ terms of the achiral contribution, with their coefficients and
evaluation conventions.  The functions are divisor sums, class numbers,
quadratic-character convolutions and square indicators.
The proof uses the arithmetic evaluations in Papers~II--IV
\cite{fullereneInvolutivePaper,fullereneRotaryPaper,fullereneOrientedPaper}.
\end{theorem}

The sector decomposition is \cref{thm:exact-all-index-grand-formula}.

As an independent check the formula reproduces
the published census through $400$ vertices, and
\cref{app:small-counting-examples} gives selected exact counts, including
the count for two million vertices.
The published point-group census also tests the seven reversing sectors
separately through $150$ vertices: its $462$ exact comparisons are
\rot{prop:rotary-point-group-census}. These checks test the formulas;
their proofs do not use census coefficients.

\begin{theorem}[Rarity and typical symmetry]\label{thm:symmetry-overview}
As $m\to\infty$, the proportion of ordinary $2m$-vertex fullerene classes
that are achiral is $\asymp m^{-5}$.  Among all classes with at most $2x$
vertices, the achiral proportion is $\asymp x^{-5}$.
The number of achiral classes on $2m$ vertices with a nontrivial
orientation-preserving automorphism is
$O(m^2)$.
Every other achiral class has a unique reversing symmetry, which is an
involution.  Among achiral classes with at most $2x$ vertices, the
proportions with a unique reflection and a unique fixed-point-free
involution tend respectively to $527/529$ and $2/529$.
\end{theorem}

Here $f\asymp g$ means that $f/g$ is bounded above and below by positive
constants for sufficiently large arguments.  The precise main terms and
proofs are \cref{thm:achiral-asymptotics,cor:asymptotic-chirality,thm:generic-achiral-symmetry,cor:cumulative-symmetry-shares}.
The proportions $527/529$ and $2/529$ are cumulative: the proportions at
one size depend on the power of two dividing $m$.  The identification of
the two involution classes with reflections and with fixed-point-free
involutions is \cref{lem:dodecagon-model}.
The exceptional proportion tends to zero: $a_m\asymp m^4$, so the
proportion of achiral classes with extra orientation-preserving symmetry
is $O(m^{-2})$.

The rarity of extra symmetry has a geometric explanation. An additional
rotation confines the period to a smaller fixed lattice. Its vector
counts grow more slowly, giving the exceptional bound in
\cref{thm:symmetry-overview}.

The first two weight blocks show the shape of the achiral formula:
\[
\begin{aligned}
a_m={}&\frac{697}{622080}\sigma_4(m)
 -\frac{41}{19440}\sigma_4(m/2)
 +\frac{697}{19440}\sigma_4(m/4)\\
&-\frac{41}{5184}H_4(4m)-\frac{41}{324}H_4(m)
 -\frac1{10368}H_4(12m)-\frac1{648}H_4(3m)\\
&-\frac1{128}H_4(4m/3)-\frac18H_4(m/3)
 +\text{terms of weight at most }4.
\end{aligned}
\]
Here $\sigma_4$ is the fourth divisor sum and $H_4$ is Cohen's
weight-$9/2$ coefficient, normalized in \cref{sec:coefficient-functions};
nonintegral arguments contribute zero. The first line has weight five.
The remaining 172 terms are printed in
\cref{app:complete-coefficient-formula}.

The special geometry of twelve points matters here. Among the four
equal-weight packages considered below, only the twelve-point package
uses the half-integral relaxation of the Deligne--Mostow condition.
Its order-three collision monodromy identifies different real chambers;
this is what lets one geometric angle type give one arithmetic class.

We count classical $(5,6)$-fullerenes. Green, Sreelesh and Arora give
oriented and graph-isomorphism counts, including symmetry subclasses,
for generalized $(3,6)$-fullerenes
\cite{greensreelesharora2026enumeration}; Green, Sreelesh and Gomatam
refine their symmetry classification
\cite{greensreeleshgomatam2026symmetries}. These papers concern triangular
and hexagonal faces, rather than pentagonal and hexagonal faces, and
include the two-connected \emph{godseye} graphs in their convention.

In sector notation, the result for every integer $m\geq1$ is
\[
 \boxed{\displaystyle
 \iso(2m)=\frac12\left(\BB_m+
     \sum_{j\in\{\mathrm{split},\mathrm{antipodal},2,3,4,10,12\}}
       M_j^{\rm open}(m)\right).}
\]
Each $M_j^{\rm open}$ is a strict-centralizer weighted mass.  The
subscripts follow two conventions: for $j\in\{2,3,4\}$ the subscript is
the step $t$ of the root-model representative $Q_t$ of
\cref{eq:phase-ramp}, whose orders are $6$, $4$ and $6$, whereas
$j\in\{10,12\}$ is the order of the symmetry, and these two terms
combine the four order-ten and the two order-twelve classes.  Their
Schoenflies operation names give the following dictionary, with
each operation understood up to inverse:
\[
\begin{array}{c|l}
j&\text{reversing operations}\\\hline
\mathrm{split},\mathrm{antipodal}&\text{reflection }\sigma,\ \text{inversion }i\\
2,3,4&S_6,\ S_4,\ S_3\\
10&S_{10},\ S_{10}^{3},\ S_5,\ S_5^{3}\\
12&S_{12},\ S_{12}^{5}
\end{array}
\]
Thus $M_{C_s}=M_{\rm split}$, $M_{C_i}=M_{\rm antipodal}$ and
$M_{S_6}=M_2$, $M_{S_4}=M_3$, $M_{S_3}=M_4$.
Here $C_s=\langle\sigma\rangle$ and $C_i=\langle i\rangle$.
There are eleven operation classes in eight named families but seven
arithmetic sectors: $M_{10}$ pools the $S_{10}$ and $S_5$ families.
The powers distinguish operations that can generate the same cyclic
subgroup. These names describe the reversing witness, not the
full point group of a graph, which may contain several witness types.
The step labels are retained for compatibility with the integral matrices.
Their
explicit arithmetic definitions and the assembled proof
are \cref{thm:exact-all-index-grand-formula}. This paper proves the
classification and normalization that make those evaluations count
fullerenes, assembles the formula from them, and derives its asymptotic
consequences.

This paper explains why seven masses answer the counting question.
Paper~II \cite{fullereneInvolutivePaper} develops the mass calculation for
involutions; Paper~III \cite{fullereneRotaryPaper} uses
its lattice and orthogonal-mass results to treat the remaining rotary
fixed spaces. Paper~IV \cite{fullereneOrientedPaper} evaluates the complex masses of rotations.
Paper~V \cite{fullereneContactsPaper} proves the pentagon-contact
criterion, single-contact mass and density bound using published
inputs independently of the other papers. Paper~VI \cite{fullereneIprPaper} assembles the
ordinary and oriented isolated-pentagon formulas.
The arithmetic proofs in Papers~II--IV
\cite{fullereneInvolutivePaper,fullereneRotaryPaper,fullereneOrientedPaper}
are independent of the counting theorem established here.
In particular, \cref{lem:fullerene-scalar-fibre-collapse} uses the
explicit order-ten operators and conjugating words constructed in
Paper~III \cite{fullereneRotaryPaper}. Their construction is independent
of this paper's class count: the operators are derived directly in the integral lattice.

The external theorem inputs to this proof chain are collected below;
the cited companion results verify their local hypotheses and
normalizations where they are used.
\begin{center}\small
\begin{tabularx}{\textwidth}{@{}p{.24\textwidth}X@{}}
\toprule Input & Exact source and use in the series\\\midrule
Completed cone periods &
\cite[Theorem~0.1, Theorem~3.4, Propositions~3.5--3.6 and
Theorem~4.1]{thurston1998}; the collision and real-equivariance
interfaces are \cref{lem:stable-collision-specialization}.\\
Configuration uniformization &
\cite[Theorem~10.18.2]{delignemostow1986} and
\cite[Theorem~3.12]{mostow1986halfIntegral}; for the full twelve-point
unitary group, \cite[Theorems~5.1(a), 5.2]{allcockBasak2025}.
See \cref{thm:twisted-inertia,lem:period-bridge}.\\
Completed-mass modularity &
\cite[Theorem~3.8]{engel2025}; \uni{prop:unitary-eisenstein-restriction}
separates this input from the proof that the masses are Eisenstein.\\
Convergent theta identity &
\cite[\S51, Proposition~8; \S52, Theorem~5]{weil1965};
\uni{prop:unitary-eisenstein-restriction} checks the degree-one unitary
range and converts the integral to strict vector mass.\\
Real representation masses &
The Hauptsatz, formula~(2), of \cite{siegel1936ii} and
\cite[S\"atze~1--2]{kneser1961};
\inv{thm:siegel-class-mass} checks the per-class and rank-three
conditions.\\
Finite arithmetic real classes &
\cite[Proposition~3.8(ii), p.~136]{borelserre1964};
\cref{lem:decorated-finite-real-classes} verifies arithmeticity and
the algebraic involution before identifying the nonabelian
cohomology set with odd involution classes.\\\bottomrule
\end{tabularx}
\end{center}

\section{Previous work and the achiral reduction}\label{sec:frontier}
\subsection{Previous counting results}

The complete generator \texttt{fullgen} of Brinkmann and Dress
\cite{brinkmanndress1997} assembles
fullerenes from patches. The later generator \texttt{buckygen}, by Brinkmann,
Goedgebeur and McKay, uses recursive expansions with a canonical
construction path and was more than $3.5$ times faster in the benchmarks
of \cite{brinkmann2012}. The reference versions used in our comparisons
are specified in \cref{sec:data-audit}.
Both generate the graphs themselves; the
arithmetic formulas discussed below evaluate their counts directly.

Thurston's period construction identifies convex equilateral triangulations
of the sphere with positive vectors in an Eisenstein lattice of signature
$(1,9)$, modulo a lattice automorphism group \cite[Thm.~0.1]{thurston1998}.
The acting group is the full unitary group, a result of Laza
\cite[Cor.~2.6]{laza2008triangulations} recorded in \cite[Prop.~3.8]{engelsmillie2018}.  A fullerene is dual to the open curvature
stratum with twelve distinct vertices of curvature one.  Mass-weighted counts
of the completed moduli space are governed by a weight-ten Eisenstein series
\cite{engelsmillie2018}.

Engel and Smillie subsequently handled the open stratum and the unweighted
oriented orbit count; their article includes a computational appendix by
Goedgebeur.  Their result places $\BB(q)$ in an explicit finite-dimensional
sum of modular and quasimodular spaces and then prunes the
answer to a combination of Eisenstein series of weights at most ten and
specified levels \cite{engel2025}. This makes $\BB_m$ arithmetically
computable. Their Remark~3.29 proposes an analytical calculation of the
underlying masses. \Uni{thm:oriented-five-sector-reduction} and
\uni{prop:oriented-collision-scalar-assembly} provide such a derivation
through the identity and four rotation sectors.
Algorithmically, ordinary counts are known
through $400$ vertices \cite{brinkmann2012,oeisA007894}.  Engel and Smillie
state the ordinary case as future work in \cite[\S1.1, p.~576]{engel2025}, and it
is not resolved in the subsequent open-problem survey
\cite[\S2, p.~2913]{bille2025open}.
That survey's displayed asymptotic constant needs a normalization
qualification: it uses $n$ for the number of vertices but prints
$C_{\rm or}=809/(2^{15}3^{13}5^2)$ for the lower envelope of
$\iso(n)/n^9$. In that convention the constant is $C_{\rm or}/2^{10}$,
with an additional $\zeta(9)$ for the upper envelope, along feasible
even $n$. The factor $2^{-1}$ forgets orientation and $2^{-9}$ changes
the index from $m$ to $n=2m$.
With $2m$ vertices, the coefficient is $C_{\rm or}/2$, as in
\cite{oeisA007894} and \cref{cor:asymptotic-chirality} below.

\subsection{The period lattice}
\label{sec:period-lattice-notation}
We fix notation for that lattice once and use it throughout.  Let $\Lambda$
be the Eisenstein fullerene lattice: a free $\mathcal{E}=\ZZ[\zeta]$-module
of rank ten, $\zeta=e^{\pi i/3}$, carrying a Hermitian form of signature
$(1,9)$, whose underlying $\ZZ$-lattice has rank twenty.  Write $J$ for the
$\ZZ$-linear endomorphism given by multiplication by $\zeta$, so that
$J^2=J-1$ and $2J-1=\sqrt{-3}$.  Write $U(\Lambda)$ for the group of
$\mathcal{E}$-linear isometries of $\Lambda$; an \emph{antiunitary} map is an
$\mathcal{E}$-antilinear isometry, that is, a $\ZZ$-linear isometry $A$ with
$AJ=J^{-1}A$.  For such an $A$ we write $\Lambda^A$ for its fixed sublattice.
A vector $v\in\Lambda$ is \emph{positive} when $(v,v)>0$, and its
\emph{quadratic index} is $m=(v,v)/2$.  Here $(v,w)$ denotes the real
even unimodular form $B$ of \cref{sec:integral-representatives}, for which
roots have norm $-2$; the Hermitian form itself gives them norm $-3$.  Let $\mathbb B=\mathbb{CH}^9$ be the
complex hyperbolic ball of positive complex lines in
$\Lambda\otimes_{\mathcal E}\mathbb C$, let the \emph{mirrors} be the
orthogonal complements of the norm-$(-2)$ root lines, that is the collision
divisors, and let $\mathbb B^\circ\subset\mathbb B$ be the complement of all
mirrors.  A positive vector or line is \emph{mirror-free}, equivalently
\emph{globally open}, when it lies in $\mathbb B^\circ$.  Call a line
\emph{rational} when it is spanned by a vector of $\Lambda$.  Since
$\mathcal E$ is a principal ideal domain, a rational line $\ell$ meets
$\Lambda$ in a free rank-one $\mathcal E$-module, whose generators form one
orbit of the units $\mathcal E^\times=\langle\zeta\rangle$; we call them the
\emph{primitive} vectors of $\ell$.  \Cref{tab:notation} summarizes this
notation, and \cref{tab:normalization-dictionary} relates the lattices and
forms used across the series. Unless explicitly restricted, vector
masses include all contents and counting identities apply to every
positive integral index $m$.

The correspondence between fullerenes and lattice vectors, with its index
normalization, is used in this form throughout.

\begin{proposition}[Fullerene period correspondence]
\label{prop:fullerene-period-correspondence}
Let $m\geq1$ and let $X_m^\circ$ be the set of positive vectors
$v\in\Lambda$ of quadratic index $(v,v)/2=m$ whose line is mirror-free.
\begin{enumerate}[label=(\alph*)]
  \item The triangulations of the sphere with at most six triangles at
        every vertex and with $2m$ triangles, up to orientation-preserving
        isomorphism, are in bijection with the $U(\Lambda)$-orbits of the
        positive vectors of $\Lambda$ of quadratic index $m$, imprimitive
        vectors included \cite[Thm.~0.1]{thurston1998}.  That the acting
        group is the full unitary group $U(\Lambda)$ is Laza's theorem
        \cite[Cor.~2.6]{laza2008triangulations}, recorded as
        \cite[Prop.~3.8]{engelsmillie2018}.
  \item A fullerene graph on $2m$ vertices is the dual of such a
        triangulation whose twelve vertices of degree five are distinct,
        and the mirrors are the collision divisors on which two of these
        vertices coincide. Here (a) uses Thurston's triangle edge-pairing
        convention for the completed space. The identification of this
        open curvature stratum with fullerene maps is explicitly stated
        in \cite[\S1.2, p.~577]{engel2025}; the classical cubic planar
        graphs and their duals are three-connected, with their graph/map
        convention recalled in \cite[\S1]{brinkmann2012}.
        Thus the dual in this open stratum is a simple, three-connected
        graph, rather than an arbitrary generalized dual map.
        Hence the orientation-preserving isomorphism
        classes of fullerenes on $2m$ vertices correspond bijectively to
        $U(\Lambda)\backslash X_m^\circ$, and
        $\BB_m=|U(\Lambda)\backslash X_m^\circ|$.
  \item The stabilizer $\operatorname{Stab}_{U(\Lambda)}(v)$ is an
        extension of the orientation-preserving automorphism group of the
        triangulation of $v$ by the local collision groups
        \cite[Prop.~2.6]{engelsmillie2018}; on the mirror-free locus the
        collision groups are trivial, so the orientation-preserving
        automorphism group of a fullerene is the strict stabilizer
        $\operatorname{Stab}_{U(\Lambda)}(v)$ of its period vector.
  \item In the Hermitian form $\star$ of \cite{engel2025}, a
        triangulation with $2n$ triangles has $v\star v=3n$
        \cite[Thm.~1.2]{engel2025}, and the real form here satisfies
        $(v,v)=\tfrac23(v\star v)$; thus $n=(v,v)/2$ is the quadratic
        index, and the fullerene has $2n$ vertices.
\end{enumerate}
\end{proposition}

\subsection{Burnside and weighted witnesses}
\begin{lemma}[Burnside reduction]\label{lem:burnside}
For every $m\geq1$,
\begin{equation}\label{eq:burnside}
  a_m=2\iso(2m)-\BB_m,
  \qquad
  \iso(2m)=\frac{\BB_m+a_m}{2}.
\end{equation}
\end{lemma}

\begin{proof}
Let $X_m=U(\Lambda)\backslash X_m^\circ$ be the finite set of
orientation-preserving isomorphism classes of fullerenes on $2m$ vertices
(\cref{prop:fullerene-period-correspondence}).  Reversal of the orientation defines an
involution $\iota$ on $X_m$. By Whitney's uniqueness theorem for the
spherical embedding of a three-connected planar graph, every abstract
graph isomorphism induces an orientation-preserving or
orientation-reversing map isomorphism; this application is recalled in
\cite[\S1]{brinkmann2012}. Its orbits are therefore the ordinary isomorphism classes,
so $|X_m/\langle\iota\rangle|=\iso(2m)$ and $|X_m|=\BB_m$.  The fixed points
of $\iota$ are exactly the achiral classes, hence there are $a_m$ of them.
Burnside's lemma gives
\[
  \iso(2m)=\frac{|X_m|+|\operatorname{Fix}(\iota)|}{2}
             =\frac{\BB_m+a_m}{2}.
\]
Rearrangement gives the first identity.
\end{proof}

We first count the fixed locus of orientation reversal. The corresponding
weighted identity for rotations will give the oriented contribution in
\cref{prop:analytical-oriented-count}.

The normalization of that fixed locus is essential.  Put
$\Gamma=U(\Lambda)$, let $\Gamma c$ be its orientation-reversing antiunitary
coset, and let $X_m^\circ$ be the mirror-free positive integral vectors
of quadratic index $m$ of \cref{prop:fullerene-period-correspondence}.
Stabilizers of such vectors are finite: the orthogonal complement of a positive vector is
negative definite, so its real unitary stabilizer is compact; intersecting
that compact group with the discrete integral group gives a finite group.
Fixing a line is weaker than fixing its integral period: if $Av_0=v_0$,
then $A(\lambda v_0)=\bar\lambda v_0$, so a scalar lift of this witness fixes $\lambda v_0$
only when $\lambda/\bar\lambda$ is an allowed scalar unit. This is the period
distinction behind handed Goldberg--Coxeter subdivisions. It does not
by itself exclude every other reversing witness; the witness sum must
test strict vectors, including imprimitive ones.

\begin{table}[htbp]
\centering\small
\begin{tabularx}{\textwidth}{@{}lX@{}}
\toprule Symbol & Meaning\\\midrule
$m$ & Half the number of vertices; the quadratic index of a period vector.\\
$\iso(2m),\BB_m,a_m$ & Ordinary, oriented, and achiral fullerene counts.\\
$\Lambda,\Gamma$ & The rank-ten Eisenstein period lattice and $U(\Lambda)$.\\
$\mathbb B^\circ,X_m^\circ$ & Collision-free positive lines and integral
vectors of index $m$, respectively.\\
$A,C_\Gamma(A)$ & An antiunitary reversing symmetry, or \emph{witness},
and the group of integral unitary maps commuting with it.\\
$M_A^{\rm open}(m)$ & Sum of reciprocal strict stabilizer orders on the
collision-free fixed locus of $A$; the sector masses $M_j^{\rm open}(m)$ of
the boxed formula sum the classes of one sector.\\
$E_m,U_m^{\rm sp},U_m^{\rm ap}$ & Achiral classes with extra preserving
symmetry, a unique reflection, or a unique antipodal involution.\\
Strict stabilizer & Maps fixing the vector itself; scalar maps are not
discarded by passing to projective space.\\
Collision flat; normal & An intersection of collision mirrors; its
perpendicular lattice, taken with positive sign.\\
State & A congruence subspace recording which lattice symmetries extend;
for involutions it is a maximal isotropic subspace modulo three.\\
Glue & An identification of discriminant groups specifying a primitive
embedding of a normal and its complement in the full lattice.\\
M\"obius value & The coefficient of a flat in the intersection poset's
inclusion--exclusion formula.\\
Atom & A minimal nonzero collision normal in the intersection arrangement.\\
Ambient lattice & The full period lattice, before restriction to a fixed
locus or a collision flat.\\
Weight component & The terms of one arithmetic weight in the formula;
their weight controls their growth.\\
\bottomrule
\end{tabularx}
\caption{Notation and terminology for the series. Arithmetic functions are defined in
\cref{app:complete-coefficient-formula}.}
\label{tab:notation}
\end{table}
\FloatBarrier

\Inv{par:worked-d4-contribution} follows one $D_4$ collision normal
through the state sum, even gluing and full-group vector mass.
The corresponding identity is \inv{cor:even-all-glue-fubini}.

\begin{table}[htbp]
\centering\small
\begin{tabularx}{\textwidth}{@{}>{\raggedright\arraybackslash}p{0.33\textwidth}X>{\raggedright\arraybackslash}p{0.2\textwidth}@{}}
\toprule
Lattice and form & Norm of a vector of quadratic index $m$ & Where used\\\midrule
$(\Lambda,h)$, the Hermitian form with Gram matrix $H_{\mathcal E}$ &
 roots have $h(\alpha,\alpha)=-3$; positive vectors have $h(v,v)=3m$ &
 \cref{sec:integral-representatives}; \inv{sec:common-lattice}\\
$(\Lambda,B)$, $B=\tfrac23\operatorname{Re}h$, the even unimodular lattice
 $II_{2,18}$ &
 roots have $B(\alpha,\alpha)=-2$; $B(v,v)=2m$, so $m=B(v,v)/2$ &
 \cref{tab:notation}, \cref{sec:integral-representatives}\\
$(F,B)$, $F=\Lambda^\chi$ the fixed lattice of an involution &
 $F=2I_{1,9}$ (split) or $2II_{1,9}$ (antipodal); $B(v,v)=2m$ &
 \inv{prop:involutive-lattice-glue}\\
$(L,B/2)$, $L=F/2$ the halved fixed lattice &
 $L=I_{1,9}$ or $II_{1,9}$; the full Gram norm $Q(x)=x^{\mathsf T}Sx$
 equals $m$ &
 \inv{prop:involutive-lattice-glue}\\
primitive polynomial $q=Q/g$ of a complement $S$,
 $g=\gcd(S_{ii},2S_{ij})$ &
 $q(x)=t$ with $t=m/g$; rows vanish unless $g\mid m$ &
 \inv{lem:general-row-normalization-interface}\\
$(\Lambda_j^{\rm rot},B/2)$, the rotary half-lattices
 $\Lambda_2^{\rm rot}=U(3)\oplus A_2(-1)$,
 $\Lambda_3^{\rm rot}=U(2)\oplus D_4(-1)$,
 $\Lambda_4^{\rm rot}=\langle1,-3,-3,-3\rangle$ &
 Gram norm equal to $m$ &
 \rot{thm:siegel-q234-ambient}, where they are denoted $H_j$\\
order-ten fixed lattices $\Lambda^{A_k}$ &
 with $B$: $10U$ ($k$ odd) or $10I_{1,1}$ ($k$ even), norm $2m$;
 with $B/2$: $5U$ or $5I_{1,1}$, norm $m$ &
 \cref{prop:mirrorfree-spectrum}; \rot{thm:order10-golden-centralizer}\\
\bottomrule
\end{tabularx}
\caption{Normalization dictionary for the series: the same vector of
quadratic index $m$ has these norms in the lattices and forms used across
the series.}
\label{tab:normalization-dictionary}
\end{table}

\begin{lemma}[Centralizer-weighted witness formula]
\label{lem:witness-formula}
Under the correspondence of \cref{prop:fullerene-period-correspondence},
the achiral correction is
\begin{equation}\label{eq:witness-formula}
 a_m=
 \sum_{[A]_\Gamma\subset\Gamma c}
 \ \sum_{[v]\in C_\Gamma(A)\backslash(X_m^\circ)^A}
 \frac{1}{|\operatorname{Stab}_{C_\Gamma(A)}(v)|}.
\end{equation}
Thus every antiunitary conjugacy class enters with coefficient one; its
centralizer quotient supplies the normalization.  Multiple
orientation-reversing symmetries of one fullerene do not overcount it.
\end{lemma}

\begin{proof}
This is the groupoid form of Burnside's lemma.  Consider witness pairs
$(A,v)$ with $Av=v$, with arrows given by the simultaneous action
$h(A,v)=(hAh^{-1},hv)$ of $\Gamma$. Its cardinality is the sum of
reciprocal pair-stabilizer orders. For one achiral $\Gamma$-orbit, let $H$ be the unitary
stabilizer of a representative and choose one antiunitary witness.  Every
point of the orbit has $|H|$ witnesses and the orbit has groupoid mass
$1/|H|$; equivalently, before quotienting, the orbit-size and witness-count
factors cancel.  The total contribution is one.  A chiral orbit has no
witness and contributes zero.  Grouping witness pairs first by the
$\Gamma$-conjugacy class of $A$ gives the two sums and the stabilizer weights
in \cref{eq:witness-formula}.
For instance, if the orientation-preserving automorphism group $H$ of an
achiral graph has order two, its two reversing symmetries form a coset of
$H$, and the conjugation action of $H$ on them has weighted orbit count
$\sum_{[s]}1/|C_H(s)|=2/|H|=1$.
\end{proof}

The identity is stated for strict vectors and strict stabilizers, whereas
the period correspondence of \cref{sec:period-lattice-notation} is
projective.  The passage between the two, including the absorption of the
scalar units $\langle\zeta\rangle$, is carried out in
\cref{rem:projective-vs-strict,lem:period-bridge,lem:fullerene-scalar-fibre-collapse}.

\part{Real periods on equal-weight ball quotients}
\label{part:real-periods}
The counting formula requires an exhaustive list of reversing symmetry
classes, with the stabilizers used in their masses. For twelve cone points
we prove that there are exactly eleven projective classes. The four
equal-weight packages provide a common setting in which to separate the
geometric angle types from the lattice conjugacy classes; the odd collision
monodromy of the twelve-point package is the step that makes one angle
type give one class.

Equal-weight Deligne--Mostow ball quotients come with antiholomorphic real
structures, but the enumeration of orientation-reversing inertia requires
more than the existence of real points.  One must preserve exact isotropy
order, centralizers, connected components of the mirror complement, and the
passage between anti-M\"obius maps and antiunitary lattice actions.  This
part proves that functorial statement uniformly for the four equal-weight
packages, derives the orbit equation for all compatible finite-order rotary
reflections, and records the complete twelve-point class list; the next part
specializes it to the fullerene package and supplies the arithmetic that the
other packages do not share.  The general real-period philosophy is
classical, and recent work treats the real binary-quintic components
explicitly \cite{deGaayFortman2025}; in the twelve-point case the ball
quotient and its mirror braid generators are available from Deligne--Mostow
theory and the refined analysis of Allcock--Basak \cite{allcockBasak2025}.
The classification here is analytical: the orbit equation, connectedness of
configuration spaces, and unordered collision monodromy prove completeness.

The fullerene count uses eight statements from this part.
\Cref{thm:twisted-inertia} turns an orientation-reversing symmetry of a
configuration into a compatible projective antiunitary class fixing its
period point, with exact order and isotropy preserved, and
\cref{lem:stable-collision-specialization} extends the identification to
the collision points of the completed quotient, with explicit local
collision charts.  \Cref{thm:equal-weight-spectra} lists the possible
angle types and orbit partitions of such a symmetry,
\cref{thm:one-partition} shows that each type and partition carries at
most one projective class, and
\cref{lem:type-constant,thm:equal-weight-uniqueness} supply the criterion
by which a type with several partitions still carries exactly one class.
\Cref{lem:type-dim} gives the dimension of each fixed locus; for the
order-ten types it is a geodesic.  The result for twelve points is
\cref{thm:open-class-list}: each of the eleven angle types carries
exactly one compatible projective class, and these are the eleven classes
behind the seven sectors of the count.  A reader who wants only the count
may take these statements as given and resume at \cref{sec:fixed-loci},
the first section of the next part.

\section{Equal-weight packages and real strata}
\label{sec:absorbed-generalizations}

We begin with the four equal-weight Deligne--Mostow packages satisfying the
integrality condition, then pass from their configuration-space real
structures to the compatible antiunitary monodromy cosets.
For fullerene enumeration one may set $N=12$ throughout. The other three
cases are retained in the common argument to isolate exactly where odd
collision monodromy, and hence the twelve-point conjugacy conclusion,
is used; no arithmetic count for those packages is imported.

\subsection{Equal-weight ball packages}

\begin{definition}[Ball package]
\label{def:ball-package}
A \emph{ball package} is a Deligne--Mostow weight vector
\(\mu=(\mu_1,\dots,\mu_N)\) with \(\sum\mu_i=2\) and \(\Sigma\mathrm{INT}\),
together with its Hermitian lattice \(L_\mu\) of signature \((1,N-3)\),
arithmetic group \(\Gamma_\mu\), ball \(\mathbb B=\mathbb{CH}^{N-3}\), root
mirrors \(\mathcal H_\mu\), open stratum \(\mathbb B^\circ\), and the period isomorphism
of orbifolds
\[
  \iota_\mu\colon
  \bigl[M_{0,N}/\operatorname{Aut}(\mu)\bigr]
  \;\longrightarrow\;
  \mathbb B^\circ/P\Gamma_\mu
\]
from \cite{delignemostow1986,mostow1986halfIntegral,thurston1998}, where
\(\operatorname{Aut}(\mu)\) is the group of permutations preserving \(\mu\).
Write \(c\) for complex conjugation on the moduli space.
This definition concerns period geometry, without assigning an integer
index to a general package. In particular, for \(N=5\) the Hermitian
norm takes values in the real subfield of \(\mathbb Q(\zeta_5)\).
Integer-indexed counting below uses the fullerene package and its
specified area normalization.
\end{definition}

The fullerene package is \(\mu=(1/6)^{12}\).

\begin{lemma}[Equal-weight \(\Sigma\mathrm{INT}\)]
\label{lem:equal-weight-int}
Among equal-weight vectors \(\mu=(2/N)^N\) with \(N\geq5\), the condition
\(\Sigma\mathrm{INT}\) holds if and only if \(N-4\) divides \(8\).
The complete list in this range is
\[
  N\in\{5,6,8,12\}.
\]
The corresponding balls are \(\mathbb{CH}^2\), \(\mathbb{CH}^3\),
\(\mathbb{CH}^5\), and \(\mathbb{CH}^9\).  Plain \(\mathrm{INT}\) would
require \(N/(N-4)\in\mathbb Z\), that is \(N-4\mid4\), and would admit only
\(N\in\{5,6,8\}\).  The fullerene package is therefore the one equal-weight
package that exists solely by Mostow's half-integral relaxation, available
because all twelve weights are equal
\cite[\S2, p.~93, condition~(2.2)]{mostow1986halfIntegral}.
\end{lemma}

\begin{proof}
Here \(1-\mu_i-\mu_j=(N-4)/N\), so
\((1-\mu_i-\mu_j)^{-1}=N/(N-4)\).  Half-integrality is
\(2N/(N-4)\in\mathbb Z\).  Since \(2N=2(N-4)+8\), this holds if and only
if \(N-4\) divides \(8\).  The condition \(\mu_i+\mu_j<1\) under which
\(\Sigma\mathrm{INT}\) is imposed reads \(4/N<1\), so \(N\geq5\); the
positive divisors of \(8\) then give \(N-4\in\{1,2,4,8\}\), that is
\(N\in\{5,6,8,12\}\).  The excluded small cases are not divisor failures:
at \(N=4\) one has \(\mu_i+\mu_j=1\) exactly and the condition is vacuous,
and at \(N=3\) the ball has dimension zero.
\end{proof}

\begin{proposition}[Coefficient ring of an equal-weight package]
\label{thm:coeff-ring}
Let \(\mu=(2/N)^N\) with \(N\in\{5,6,8,12\}\), and let \(d\) be the least
common denominator of the weights.  Then \(d=N/\gcd(2,N)\), so
\[
  d=5,3,4,6
\]
respectively.  The Deligne--Mostow lattice \(L_\mu\) is a Hermitian lattice
over the ring of integers of \(\mathbb Q(\zeta_d)\).  For \(N=6,8,12\) this
field is imaginary quadratic, with ramified prime \(p=3,2,3\), the one
rational prime whose ideal is the square of a prime ideal of the field.
For \(N=5\) the field is \(\mathbb Q(\zeta_5)\), ramified at \(5\).
\end{proposition}

\begin{proof}
This is the standard description of the Deligne--Mostow period lattice
\cite{delignemostow1986}.
Each weight is \(2/N\), so \(d=N/\gcd(2,N)\).  The cyclic cover is
\(y^d=\prod_i(x-z_i)^{d\cdot 2/N}\), and the period lattice is the
primitive \(\zeta_d\)-eigenspace of \(H_1\) of that cover, with coefficients
in \(\mathbb Z[\zeta_d]\) \cite{delignemostow1986}.  The fields
\(\mathbb Q(\zeta_3)=\mathbb Q(\zeta_6)\) and \(\mathbb Q(\zeta_4)=\mathbb Q(i)\)
are imaginary quadratic, ramified at \(3\) and \(2\).  The prime \(5\) ramifies
in \(\mathbb Q(\zeta_5)\).
\end{proof}

\subsection{Twisted inertia}

The period correspondence must preserve the symmetries and their
stabilizers, since those stabilizers are the denominators in
\cref{lem:witness-formula}.  The next theorem supplies precisely this
compatibility.  Its holomorphic period-map hypothesis is the
Deligne--Mostow uniformization input
\cite{delignemostow1986,mostow1986halfIntegral}; equivariance under
orientation reversal and the strict-lift conclusion are proved below.

\paragraph{Period-map input (DM).}\label{par:dm-input}
We use the orbifold, rather than coarse-space, uniformization: on the
open stratum the marked period covering has kernel equal to the kernel
of monodromy, and its image group is \(P\Gamma_\mu\). For \(N=5,6,8\)
the ordered form is \cite[Theorem~10.18.2]{delignemostow1986};
passing to unordered equal weights uses the finite permutation
extension in \cite[\S3.11 and Theorem~3.12]{mostow1986halfIntegral}.
For \(N=12\), the latter theorem gives the half-integral extension,
with the full-unitary identification specified in
\cref{lem:period-bridge}(a). The quotient by
\(\operatorname{Aut}(\mu)=S_N\) is already included in this extended
monodromy group; it is not a further quotient of \(P\Gamma_\mu\).
For the other packages we use the actual monodromy image and make no
replacement of it by a larger full unitary group.

\begin{theorem}[Twisted inertia in the compatible monodromy coset]
\label{thm:twisted-inertia}
Let \(\mu\) be a ball package, let \(d\) be the least common denominator of
the weights, assume \(N\geq5\), let \(\zeta=e^{2\pi i/d}\), write
\(\mathcal O_d=\mathbb Z[\zeta]\), and let \(C_Z\) be the cyclic cover
\(y^d=\prod_i(x-z_i)^{d\mu_i}\) with deck generator
\(\tau(x,y)=(x,\zeta y)\).  Under (DM), the period
map \(P:\widetilde M\to\mathbb B^\circ\) of a primitive eigenspace of
\(\tau^*\) is a Galois covering onto \(\mathbb B^\circ\), its monodromy
\(\rho\) is onto \(P\Gamma_\mu\), and its deck group is \(\ker\rho\);
equivalently, it induces the asserted orbifold isomorphism with
\(\mathbb B^\circ/P\Gamma_\mu\).
The conclusions have two parts.

\emph{Conjugation.}
Then conjugation \((x,y)\mapsto(\bar x,\bar y)\) induces an antiholomorphic
isomorphism \(C_Z\to C_{\bar Z}\) intertwining \(\tau\) with \(\tau^{-1}\).
It supplies a projective antiunitary \(\bar\chi_0\) normalizing
\(P\Gamma_\mu\).  Put
\[
 \widetilde{P\Gamma}_\mu
 :=P\Gamma_\mu\sqcup P\Gamma_\mu\bar\chi_0
 \subset P\widetilde U(L_\mu).
\]
\emph{Symmetry pairs and their stabilizers.}
The period isomorphism \(\iota_\mu\) is equivariant for \(c\) and identifies
the \(c\)-twisted inertia groupoid of
\([M_{0,N}/\operatorname{Aut}(\mu)]\) with the groupoid of pairs
\((\bar A,[v])\), where \(\bar A\) lies in the compatible odd coset
\(P\Gamma_\mu\bar\chi_0\) and fixes
\([v]\in\mathbb B^\circ\).  Exact projective order and isotropy are
preserved.  This identifies all projective antiunitaries only in the special
case \(P\Gamma_\mu=PU(L_\mu)\).
\end{theorem}

\begin{proof}
Write \(M=M_{0,N}\) for ordered configurations, \(\mathcal M=[M/\operatorname{Aut}(\mu)]\),
and \(\widetilde M\to M\) the universal cover, with orbifold fundamental
group \(\pi\) of \(\mathcal M\).  Let \(P:\widetilde M\to\mathbb B^\circ\)
be the period map of the chosen \(\zeta^r\)-eigenspace, with monodromy
\(\rho:\pi\to P\Gamma_\mu\), assumed a Galois covering with
\(\rho\) surjective.

(a) Conjugation on \(\mathbb P^1\) commutes with \(\operatorname{Aut}(\mu)\)
and normalizes \(\mathrm{PGL}_2(\mathbb C)\), so it descends to an
involution \(c\) of \(\mathcal M\) which is not the identity.  Choose a
conjugation-fixed base configuration \(Z_0\) with trivial stabilizer and a
lift \(\tilde x_0\in\widetilde M\).  The unique lift \(\tilde c\) of \(c\)
fixing \(\tilde x_0\) satisfies \(\tilde c^2=1\), normalizes \(\pi\), and
lies outside \(\pi\): it is antiholomorphic, whereas the elements of
\(\pi\) act holomorphically. Put
\(\tilde\pi=\pi\rtimes\langle\tilde c\rangle\).

(b) Let \(\kappa_Z:C_Z\to C_{\bar Z}\), \((x,y)\mapsto(\bar x,\bar y)\).
The exponents \(d\mu_i\) are real, so \(\bar Z\) carries the same cover
type.  Then \(\kappa_Z\tau=\tau^{-1}\kappa_Z\).  Define
\[
  K_Z(u)=\overline{(\kappa_Z^{-1})^*u}
  \qquad\text{on }H^1(C_Z,\mathbb C).
\]
This is \(\mathbb C\)-antilinear, preserves integral lattices, and, because
\(\kappa_Z\tau=\tau^{-1}\kappa_Z\), maps each \(\tau^*\)-eigenspace to
itself after coefficient conjugation.  Pullback by the antiholomorphic map
sends \((1,0)\)-classes to \((0,1)\)-classes and conjugation sends them
back, so \(K_Z\) preserves the Hodge line of the chosen eigenspace.
With \(h(u,v)=i\int v\wedge\bar u\) on the chosen eigenspace, up to a
positive constant, the form is linear in its second variable and
\[
 h(K_Zu,K_Zv)
 =i\int_{C_{\bar Z}}\overline{(\kappa_Z^{-1})^*v}\wedge(\kappa_Z^{-1})^*u
 =i\,\overline{\deg(\kappa_Z^{-1})\int_{C_Z}v\wedge\bar u}
 =\overline{h(u,v)},
\]
because \(\deg\kappa_Z^{-1}=-1\).  Thus
\(K_Z\) is an antiunitary isometry of the primitive eigenlattice, and
\(K_{\bar Z}K_Z=1\).  The maps \(\kappa_Z\) are the fibres of an
antiholomorphic involution of the total space of the family of curves over
\(M\) covering \(c\); hence \(K\) is an isomorphism of local systems over
\(c\) and commutes with Gauss--Manin transport: for a path \(\gamma\) from
\(Z_0\) to \(Z\) with conjugate path \(\bar\gamma\) from \(Z_0\) to
\(\bar Z\), \(K_Z\circ T_\gamma=T_{\bar\gamma}\circ K_{Z_0}\).  Let
\(\chi_0\) be the antiunitary isometry of \(L_\mu\) corresponding to
\(K_{Z_0}:\Lambda_{Z_0}\to\Lambda_{Z_0}\) under the marking at
\(\tilde x_0\).  For \(x\in\widetilde M\) over \(Z\) with marking
\(T_x\), the point \(\tilde c(x)\) lies over \(\bar Z\) with marking
\(T_{\tilde c(x)}=K_ZT_x\chi_0^{-1}\) by the transport identity, so, writing
\(H^{1,0}(C_Z)_{\zeta^r}\) for the Hodge line of the chosen eigenspace,
\begin{align*}
 P(\tilde c(x))
 &=T_{\tilde c(x)}^{-1}\bigl(H^{1,0}(C_{\bar Z})_{\zeta^r}\bigr)\\
 &=\chi_0T_x^{-1}K_Z^{-1}\bigl(H^{1,0}(C_{\bar Z})_{\zeta^r}\bigr)\\
 &=\chi_0T_x^{-1}\bigl(H^{1,0}(C_Z)_{\zeta^r}\bigr)
 =\bar\chi_0P(x).
\end{align*}
Uniqueness holds because \(P(\widetilde M)=\mathbb B^\circ\) is open in
\(\mathbb B\) and an isometry of the ball is determined by its restriction
to an open set.  For \(\gamma\in\pi\), equivariance gives
\(P(\tilde c\gamma x)=\bar\chi_0\rho(\gamma)P(x)\) and
\(P(\tilde c\gamma\tilde c^{-1}\cdot\tilde cx)
 =\rho(\tilde c\gamma\tilde c^{-1})\bar\chi_0P(x)\), so the two
projective isometries agree on the open set \(P(\widetilde M)\) and
\[
  \rho(\tilde c\gamma\tilde c^{-1})
  =\bar\chi_0\rho(\gamma)\bar\chi_0^{-1}.
\]

(c) Set \(\tilde\rho(\tilde c)=\bar\chi_0\). Equivariance and
\(\tilde c^2=1\) give \(\bar\chi_0^2P(x)=P(x)\) for every
\(x\in\widetilde M\). Since \(P\) is onto \(\mathbb B^\circ\),
the unitary \(\bar\chi_0^2\) fixes that open set pointwise and is the
projective identity. Hence \(\tilde\rho\) is a homomorphism
\(\tilde\pi\to P\widetilde U(L_\mu)\), with image precisely
\(\widetilde{P\Gamma}_\mu\) because \(\rho\) hits \(P\Gamma_\mu\) and
\(\bar\chi_0\) represents its compatible odd coset, and with kernel
\(\ker\rho\), since an element \(\gamma\tilde c\) is never in the
kernel: an antiunitary is not the identity.  As \(P\) is
\(\tilde\rho\)-equivariant and is a Galois covering of \(\mathbb B^\circ\)
with deck group \(\ker\tilde\rho\), it induces
\(\widetilde M/\ker\tilde\rho\cong\mathbb B^\circ\), equivariant for
\(\tilde\pi/\ker\tilde\rho\cong\widetilde{P\Gamma}_\mu\).
Equivariance of \(P\) thus gives an isomorphism of quotient orbifolds
\[
  \bigl[\widetilde M/\tilde\pi\bigr]
  \;\cong\;
  \bigl[\mathbb B^\circ/\widetilde{P\Gamma}_\mu\bigr].
\]

(d) An isomorphism of quotient orbifolds identifies isotropy groups.
Explicitly, \(\tilde\rho\) is injective on \(\tilde\pi_x\) because
\(\ker\rho\) acts freely on \(\widetilde M\), and it is surjective onto
the isotropy group of \(P(x)=[v]\): if \(h[v]=[v]\), choose
\(g\in\tilde\pi\) with \(\tilde\rho(g)=h\); then \(P(gx)=P(x)\), so
\(gx=kx\) for some \(k\in\ker\rho\), and \(k^{-1}g\in\tilde\pi_x\)
maps to \(h\).  The coset condition and the orders are preserved since
\(\tilde\rho\) is a group isomorphism on \(\tilde\pi_x\).  The
\(c\)-twisted inertia of the configuration orbifold has objects
\((Z,\sigma)\) with \(\sigma\) anti-M\"obius, \(\sigma(Z)=Z\), and
morphisms holomorphic conjugacies; on the ball it has objects
\((\bar A,[v])\) with
\(\bar A\in P\Gamma_\mu\bar\chi_0\) and \(\bar A[v]=[v]\) in
\(\mathbb B^\circ\).  No equation \(\bar A^2=1\) is imposed: twisted
inertia records every odd arrow.  Exact orders match by full faithfulness.
\end{proof}

The preceding identification also specifies the topology of the symmetry
moduli. For a compatible class \(\bar A\), its part of the moduli of pairs
is the quotient orbifold
\[
 \bigl[(\operatorname{Fix}(\bar A)\cap\mathbb B^\circ)/
                 C_{P\Gamma_\mu}(\bar A)\bigr].
\]
Indeed, in a local period chart an odd arrow is an element of a discrete
group, so fixing that arrow cuts out its fixed locus. An arrow between two
such pairs conjugates their odd elements; after choosing the representative
\(\bar A\), the remaining arrows are precisely its centralizer. The
equivariant local period maps and their inverses restrict to these fixed
loci. Thus the identification is an isomorphism of real-analytic
orbifolds, not just a bijection of isomorphism classes, and its underlying
quotient spaces have the same connected components and dimensions.

The projective correspondence does not by itself choose an integral
lift. The following statement gives existence on primitive vectors and
uniqueness on every nonzero vector, and preserves the stabilizer used
in the witness mass.

\begin{lemma}[Strict lifts]\label{lem:primitive-strict-lifts}
In the setting of \cref{thm:twisted-inertia}, suppose that
\(\Gamma_\mu\) is the full inverse image of \(P\Gamma_\mu\) and
contains the scalar isometry kernel. Let \((\bar A,[v])\) be a
compatible projective symmetry pair, and let \(v\) be a nonzero lattice
vector on \([v]\). A strict lift of \(\bar A\) fixing \(v\), if it
exists, is unique and has the same order as \(\bar A\).
If \(v\) is \(\mathcal O_d\)-primitive, meaning
\(\mathbb Q(\zeta_d)v\cap L_\mu=\mathcal O_dv\), then there is a unique
strict lift \(A\) of \(\bar A\) fixing \(v\). The projective isotropy
of the pair is identified with
\(\operatorname{Stab}_{C_{\Gamma_\mu}(A)}(v)\).
No existence assertion is made here for an imprimitive vector.
\end{lemma}

\begin{proof}
Two strict lifts fixing the same nonzero vector differ by a scalar
that fixes that vector, hence are equal. If \(A^e=\varepsilon I\),
then \(v=A^ev=\varepsilon v\), so \(\varepsilon=1\); strict order
therefore equals projective order. Neither argument uses primitivity.

Now let \(v\) be \(\mathcal O_d\)-primitive. If an
integral lift satisfies \(A_0v=\lambda v\), primitivity gives
\(\lambda\in\mathcal O_d\), while applying the same argument to \(A_0^{-1}\)
gives \(\lambda^{-1}\in\mathcal O_d\). Antiunitarity gives
\(\lambda\bar\lambda=1\), so \(\lambda\) is a scalar isometry; multiplying by \(\lambda^{-1}\)
gives a strict lift fixing \(v\). If a projective stabilizer arrow \(g\)
sends \(v\) to \(\varepsilon v\), put \(g'=\varepsilon^{-1}g\), so
\(g'v=v\). Projective commutation gives
\(g'Ag'^{-1}=uA\) for a scalar isometry \(u\).
Evaluating at \(v\) gives \(v=uv\), hence \(u=1\).
Thus the normalized lift strictly commutes with \(A\), and the
construction is reversible. It identifies
\(\operatorname{Stab}_{C_{\Gamma_\mu}(A)}(v)\) with the M\"obius symmetries
of \(Z\) commuting with \(\sigma\).

Finally, the statement does not depend on the choice of generator.
Another \(\mathcal O_d\)-primitive generator of the same line and the same
index is \(uv\) with \(u\in\mathcal O_d^\times\) and \(u\bar u=1\), hence
a root of unity \(u=\pm\zeta^b\).  By antilinearity the lift fixing
\(\zeta^bv\) is \(\zeta^{2b}A=\zeta^bA\zeta^{-b}\), so replacing \(v\) by
\(\zeta^bv\) conjugates the strict lift by the scalar isometry
\(\zeta^b\in\Gamma_\mu\), and the sign changes nothing because
\(A(-v)=-v\).  Scalars are central in \(\Gamma_\mu\), so
\(C_{\Gamma_\mu}(\zeta^{2b}A)=C_{\Gamma_\mu}(A)\), and a linear map fixes
\(uv\) if and only if it fixes \(v\); the stabilizer
\(\operatorname{Stab}_{C_{\Gamma_\mu}(A)}(v)\) is therefore the same
subgroup for every generator.
\end{proof}

The argument does not use \(N=12\) or \(d=6\) except as the special-case
values of the same letters.  A lift fixing a primitive generator need not
fix an imprimitive multiple of that generator.
For example, with \(d=6\), let \(A_0\) be a
strict antiunitary fixing a primitive vector \(v_0\) and put
\(v=\sqrt{-3}\,v_0=(2\zeta-1)v_0\).  Antilinearity gives
\(A_0v=\overline{\sqrt{-3}}\,v_0=-v\), so the imprimitive vector \(v\) is
fixed by the lift \(-A_0=\zeta^3A_0\) of the same projective class and
not by \(A_0\).  A coefficient formula must therefore sum the fixed-vector
masses of all strict lifts above a compatible projective class
(\cref{def:strict-scalar-fibres}), and a separate argument is needed to
collapse that sum to one strict class
(\cref{lem:fullerene-scalar-fibre-collapse}).

\subsection{Normal forms and orbit partitions}

\begin{proposition}[Rotary-reflection normal form]
\label{thm:rotary}
Every finite-order anti-M\"obius transformation of \(\mathbb P^1\) is
holomorphically conjugate to
\[
  A_{n,k}=\operatorname{diag}(R_{\pi k/n},-1),\qquad n\geq1,\quad
  k\in\{0,1,\dots,2n-1\}.
\]
Here \(\mathbb P^1\) is identified with the unit sphere
\(S^2\subset\mathbb R^3\), and the matrix acts through the block subgroup
\(O(2)\times O(1)\subset O(3)\): \(R_{\pi k/n}\) rotates the equatorial
plane by the angle \(\pi k/n\) and \(-1\) reverses the polar axis.
Its exact order is \(2n\) if and only if
\(\operatorname{lcm}(2n/\gcd(k,2n),2)=2n\), equivalently
\(\gcd(k,2n)\in\{1,2\}\) with \(n\) odd in the second case.  Holomorphic
conjugacy identifies \(k\) with \(2n-k\), so that \(0\leq k\leq n\) is a
fundamental domain for this identification.  The case \(k=0\) occurs only
for \(n=1\), and is a plane reflection.  The case \((n,k)=(1,1)\) is the
antipodal map \(-I\).
\end{proposition}

\begin{proof}
This is the classical classification of finite orientation-reversing
isometries of the sphere.  A finite-order element of
\(\mathrm{PGL}_2(\mathbb C)\rtimes\langle c\rangle\)
outside \(\mathrm{PGL}_2(\mathbb C)\) generates a finite group of
isometries of hyperbolic three-space through the Poincar\'e extension;
that group fixes the circumcentre of any of its orbits, so the element is
conjugate into the stabilizer \(O(3)\) of a point, acting on
the Riemann sphere, and the orientation-reversing finite-order elements of
\(O(3)\) are rotary reflections, including angle zero (plane reflections)
and the central inversion.  The square is
\(A_{n,k}^2=\operatorname{diag}(R_{2\pi k/n},1)\), of order
\(n/\gcd(k,n)\) for \(k\neq0\).  Odd powers remain improper, so the order
is even and equals \(\operatorname{lcm}(2n/\gcd(k,2n),2)\).  Setting this
equal to \(2n\) yields the gcd condition.  For \(k=0\) one has
\(A_{n,0}=\operatorname{diag}(I,-1)\) of order \(2\), so exact order
\(2n\) forces \(n=1\).  Conjugation by a holomorphic half-turn about an axis
perpendicular to the rotation axis sends \(k\) to
\(-k\equiv 2n-k\pmod{2n}\).
\end{proof}

\begin{theorem}[Orbit equation]
\label{thm:equal-weight-spectrum}
Let \(A_{n,k}\) have exact order \(2n\).  Write
\(g_k=\gcd(k,2n)\).  Off the axis, a point has \(\langle A_{n,k}\rangle\)-orbit
of size \(2n\) unless it lies on the equator and \(g_k=2\), in which case
the orbit has size \(n\).  Equivalently, equatorial \(n\)-orbits exist if
and only if \(A_{n,k}^n\) is reflection through the equator, which holds
precisely when \(g_k=2\).  An invariant set of \(N\) distinct points
therefore satisfies
\begin{equation}
\label{eq:orbit}
  2b+2ng+ns=N,
\end{equation}
where \(b\in\{0,1\}\) records the axial pair, \(g\geq0\) is the number of
generic \(2n\)-orbits, and \(s\geq0\) is the number of equatorial
\(n\)-orbits, with \(s=0\) unless \(g_k=2\).

Two degenerate involutions are recorded separately.  For the plane
reflection \(A_{1,0}\) the axis is not distinguished among conjugate
pairs, and the partitions are \(2g+s=N\) with \(s\equiv N\pmod{2}\).  For
the antipodal map \(A_{1,1}\) every orbit has size two, and there is a
single partition if and only if \(N\) is even.

Conversely every solution is realized by placement.  The possible orders
of a finite anti-M\"obius map preserving \(N\) distinct points are exactly
the values \(2n\) for which some exact-order \(k\) admits a solution.
\end{theorem}

\begin{proof}
The map \(A_{n,k}\) sends the north pole to the south pole, so the two
axis points form one orbit of size two.  Off the axis and off the equator
the group has order \(2n\) and acts freely.  On the equator,
\(A_{n,k}\) reduces to rotation by \(\pi k/n\), of orbit size
\(2n/g_k\).  This equals \(n\) if and only if \(g_k=2\).  In that case
\(n\) is odd by \cref{thm:rotary}, and
\(A_{n,k}^n=\operatorname{diag}(I,-1)\).  Counting an invariant set of
\(N\) distinct points gives \eqref{eq:orbit}.

For \(A_{1,0}\) the equator is pointwise fixed and every off-equator orbit
has size two.  For \(A_{1,1}=-I\) every orbit is an antipodal pair, which
can cover \(N\) distinct points only if \(N\) is even.  Distinctness of
the \(N\) points is the open condition in \(M_{0,N}\).
\end{proof}

\begin{theorem}[Spectra of the equal-weight packages]
\label{thm:equal-weight-spectra}
For the four equal-weight ball packages of \cref{lem:equal-weight-int},
the exact-order angle types \((n,k)\) with \(0\leq k\leq n\), their
orbit-size partitions, and the number of partitions are as follows.
Special types (\(g_k=2\)) are marked \(s\).  The last column names the
integral representatives of \cref{sec:fixed-loci}: the root-monomial
operators \(Q_t\) of \cref{eq:phase-ramp}, matched with these types in
\cref{lem:dodecagon-model}, and the order-ten operators
\(A_k\) constructed in \cref{sec:integral-representatives}; a dash marks a
package for which no integral representative is constructed in this
paper.
\begin{center}
\small
\begin{tabular}{@{}ccccp{43mm}cc@{}}
\toprule
\(N\) & \(2n\) & \(k\) & & partitions & \(\#\) & Representative \\
\midrule
\(5\) & \(2\) & \(0\) & &
  \([2^2,1]\), \([2,1^3]\), \([1^5]\) & \(3\)  & --- \\
      & \(6\) & \(2\) & \(s\) & \([3,2]\) & \(1\)  & --- \\
      & \(10\) & \(2,4\) & \(s\) & \([5]\) each & \(1\)  & --- \\
\addlinespace
\(6\) & \(2\) & \(0\) & &
  \([2^3]\), \([2^2,1^2]\), \([2,1^4]\), \([1^6]\) & \(4\)  & --- \\
      & \(2\) & \(1\) & & \([2^3]\) (antipodal) & \(1\)  & --- \\
      & \(4\) & \(1\) & & \([4,2]\) & \(1\)  & --- \\
      & \(6\) & \(1\) & & \([6]\) & \(1\)  & --- \\
      & \(6\) & \(2\) & \(s\) & \([6]\), \([3^2]\) & \(2\)  & --- \\
\addlinespace
\(8\) & \(2\) & \(0\) & &
  \([2^4]\), \([2^3,1^2]\), \([2^2,1^4]\), \([2,1^6]\), \([1^8]\) & \(5\)  & --- \\
      & \(2\) & \(1\) & & \([2^4]\) (antipodal) & \(1\)  & --- \\
      & \(4\) & \(1\) & & \([4^2]\) & \(1\)  & --- \\
      & \(6\) & \(1\) & & \([6,2]\) & \(1\)  & --- \\
      & \(6\) & \(2\) & \(s\) & \([6,2]\), \([3^2,2]\) & \(2\)  & --- \\
      & \(8\) & \(1,3\) & & \([8]\) each & \(1\)  & --- \\
\addlinespace
\(12\) & \(2\) & \(0\) & &
  \([2^6]\), \([2^5,1^2]\), \([2^4,1^4]\), \([2^3,1^6]\), \([2^2,1^8]\),
  \([2,1^{10}]\), \([1^{12}]\) & \(7\)  & $Q_0$ \\
      & \(2\) & \(1\) & & \([2^6]\) (antipodal) & \(1\)  & $Q_6$ \\
      & \(4\) & \(1\) & & \([4^3]\) & \(1\)  & $Q_3$ \\
      & \(6\) & \(1\) & & \([6^2]\) & \(1\)  & $Q_2$ \\
      & \(6\) & \(2\) & \(s\) & \([6^2]\), \([6,3^2]\), \([3^4]\) & \(3\)  & $Q_4$ \\
      & \(10\) & \(1,3\) & & \([10,2]\) each & \(1\)  & $A_1,A_3$ \\
      & \(10\) & \(2,4\) & \(s\) & \([10,2]\), \([5^2,2]\) each & \(2\)  & $A_2,A_4$ \\
      & \(12\) & \(1,5\) & & \([12]\) each & \(1\)  & $Q_1,Q_5$ \\
\bottomrule
\end{tabular}
\end{center}
In the reflection rows the number of parts equal to \(1\) is \(s\), the
number \(r\) of real points of the configuration, that is, of points on
the fixed circle of the reflection; for \(N=12\) this
\(r\in\{0,2,\dots,12\}\) is the index used in \cref{thm:split-class}.
The numbers of angle types are \(4,5,7,11\).  For \(N=12\) the eleven
types recover the special-case table: the two order-six types are \(Q_2\)
(\(k=1\)) and \(Q_4\) (\(k=2\)); the four order-ten types are the odd-angle
pair \(k=1,3\) and the even-angle pair \(k=2,4\); the two order-twelve
types are \(Q_1,Q_5\).
\end{theorem}

\begin{proof}
Apply \cref{thm:rotary,thm:equal-weight-spectrum} with \(N\in\{5,6,8,12\}\).
Range over \(n=1,\dots,N\) and holomorphic representatives
\(k\in\{0,\dots,n\}\) of exact order \(2n\).  For each pair, list
nonnegative integer solutions of \eqref{eq:orbit} with \(s=0\)
unless \(g_k=2\), using the two involutive conventions of
\cref{thm:equal-weight-spectrum}.  Translating \((b,g,s)\) as
\([ (2n)^g\, n^s\, 2^b ]\) gives the displayed partitions.  For \(n>1\),
at most one of \(b\in\{0,1\}\) occurs: \eqref{eq:orbit} forces
\(2b\equiv N\pmod n\), which separates the two choices when $n>2$.
For $n=2$ the type is nonspecial, so $s=0$ and the congruence modulo
$2n=4$ separates them instead.
\end{proof}

\begin{theorem}[One class per partition]
\label{thm:one-partition}
Under (DM), fix an angle type \((n,k)\) and one
admissible partition \((b,g,s)\).  The moduli space of symmetry pairs
\((Z,\sigma)\), with \(Z\) an unordered configuration of \(N\) distinct
points and \(\sigma\) an anti-M\"obius map of angle type \((n,k)\)
preserving \(Z\) with orbit partition \((b,g,s)\), taken up to holomorphic
conjugacy, is connected: generic
orbits are parametrised by distinct points of a connected surface,
equatorial orbits by distinct points of a circle, and the axial pair is rigid.
The compatible projective \(P\Gamma_\mu\)-class attached to \((Z,\sigma)\)
by \cref{thm:twisted-inertia} is locally
constant on that space, hence constant.  Thus each pair (angle type,
partition) carries at most one compatible projective antiunitary class.
\end{theorem}

\begin{proof}
For \(n>1\), quotient first by the rotation \(\sigma^2\) on the
complement of the poles.  In the coordinate \(w=z^n\), the remaining
involution is \(w\mapsto(-1)^k/\bar w\).  For odd \(k\) its quotient
is an open M\"obius band.  For even \(k\), the free-orbit quotient is
a punctured disc and the equatorial quotient is its boundary circle.
For plane reflection the free-orbit quotient is a disc; for the
antipodal map it is \(\mathbb{RP}^2\).  Configuration spaces of finitely
many distinct unordered points on any of these connected surfaces, or
on a circle, are connected.  The free, equatorial and axial loci are
disjoint, so their prescribed configuration spaces form a connected
product.  Taking the residual holomorphic quotient preserves
connectedness.  The class is locally constant: by the orbifold
isomorphism of \cref{thm:twisted-inertia}, the pairs \((Z,\sigma)\) near
a given pair lie in one orbifold chart of the \(c\)-twisted inertia
stack, in which the odd arrow is represented by a single element
\(\bar A\) of the compatible coset, so every pair of that chart receives
the class of \(\bar A\). A conjugacy class labels an open-and-closed
part of the twisted inertia stack, which can have several connected
components. Conjugacy classes in the
compatible odd coset \(P\Gamma_\mu\bar\chi_0\) are discrete, so a
locally constant class function on a connected space is constant.
\end{proof}

\subsection{Dimensions and collision specialization}
\begin{lemma}[Dimension of a type]
\label{lem:type-dim}
For a rotary type with \(n>1\), the real dimension of the moduli of pairs
\((Z,\sigma)\) of partition \((b,g,s)\) is \(2g+s-1\).  All partitions of
one type have the same dimension, equal to
\(\dim\bigl(\operatorname{Fix}(\bar A)\cap\mathbb B\bigr)\) for a
compatible class \(\bar A\) of that type; this equality follows from the
local orbifold isomorphism of \cref{thm:twisted-inertia}.  For the
plane reflection and for the antipodal map the dimension is \(N-3\).
\end{lemma}

\begin{proof}
After the normal form of \(\sigma\), residual holomorphic conjugacy is
rotation about the axis, a circle.  Generic orbits contribute two real
parameters each, equatorial orbits one.  For a fixed type, \(n>1\)
determines a unique \(b\) by the last paragraph of
\cref{thm:equal-weight-spectra}, and \eqref{eq:orbit} then makes
\(2g+s=(N-2b)/n\) constant.  The identification of
\cref{thm:twisted-inertia} is an isomorphism of orbifolds, so it matches
the moduli of pairs \((Z,\sigma)\) of the type with the disjoint union,
over the compatible classes \(\bar A\) of that type, of
\((\operatorname{Fix}(\bar A)\cap\mathbb B^\circ)/C_{P\Gamma_\mu}(\bar A)\);
the discrete group \(C_{P\Gamma_\mu}(\bar A)\) acts properly on the real
hyperbolic space \(\operatorname{Fix}(\bar A)\cap\mathbb B\), so the
quotient has the dimension of that space, which is therefore \(2g+s-1\)
(and \(N-3\) in the two involutive cases below).  The reflection case is
the real locus of
\(M_{0,N}/S_N\), of dimension \(N-3\).  For the antipodal map every orbit
is a pair, the residual holomorphic conjugacy is \(SO(3)\), and
\(2\cdot(N/2)-3=N-3\) again.
\end{proof}

\begin{lemma}[Specialization through simultaneous stable collisions]
\label{lem:stable-collision-specialization}
For each equal-weight package \(N\in\{5,6,8,12\}\), use the completed
Deligne--Mostow--Thurston period identification: points of the ball quotient
correspond to weighted configurations whose collision clusters have total
weight strictly less than one.  This is the metric completion theorem and
its local orbifold construction in
\cite[Thm.~3.4, Props.~3.5--3.6 and Thm.~4.1]{thurston1998}; for
\(N=12\), the identification with the full stable configuration quotient
is also stated in \cite[\S2]{heckman2015}.
Let \(x_0\in\mathbb B\) correspond to such a configuration \(Z_0\).
The compatible extended stabilizer at \(x_0\) has a natural map to the
group of holomorphic and antiholomorphic automorphisms of \(Z_0\).
Local collision monodromy acts trivially under this map.  If a fixed
compatible odd element \(\bar A\) fixes mirror-free points tending to
\(x_0\), the corresponding anti-M\"obius maps specialize, up to a common
change of marking, to the image \(\sigma_0\) of \(\bar A\).

This statement allows several clusters at once; the map on boundary
stabilizers need not be injective.  At disjoint double
collisions, the transverse local quotient has coordinates
\[
 t_j=\xi_j^q,\qquad q=10,6,4,3
 \quad\text{for }N=5,6,8,12,
\]
after absorbing nonvanishing analytic units.  Here \(t_j\) is the
unordered smoothing coordinate and \(\xi_j\) is a period coordinate.
These equations hold simultaneously for all colliding pairs.
This normalization is holomorphic; the real-equivariant normalization
$t=\epsilon\xi^q$, retaining $\epsilon\in\{\pm1\}$, is made in
\cref{lem:equatorial-merge}.
\end{lemma}

\begin{proof}
\emph{Boundary stabilizers.}
The completed theorem concerns all stable clusters, not only the smooth
part of the collision divisor.  Its local construction cuts out disjoint
small discs containing the clusters and replaces each by a cone whose
curvature is the sum of the curvatures in that cluster
\cite[Thm.~3.4 and its proof]{thurston1998}.
Choose the discs on the limiting sphere first, and then a sufficiently
small common neighbourhood in which their boundaries remain disjoint
and no marked point crosses a boundary.  Label the points temporarily.
For each cluster the interior marking can vary by its local braid
group; variations supported in different discs commute. Their period
images fix the lifted collision point $x_0$. The period group is discrete,
whereas the stabilizer of $x_0$ in the complex hyperbolic isometry group
is compact, so these images are finite. This is also the finite local
group in Thurston's orbifold construction \cite[Thm.~4.1]{thurston1998}.
An arrow fixing the
limiting period permutes clusters of equal total weight and acts on the
common exterior.  Collapsing each disc to its cone point therefore gives
one automorphism of the collapsed sphere.  Changing an interior marking
composes the arrow with a local collision element and leaves this
automorphism unchanged.  Thus the map is well defined even at an
intersection of several collision strata; one does not specialize the
clusters successively or choose an order for them.  This is the local stabilizer construction
of \cite[Props.~3.5--3.6 and Thm.~4.1]{thurston1998}; in the twelve-point case its holomorphic
exact sequence is recorded in \cite[Prop.~2.14]{engel2025}.

\emph{Odd arrows and transverse charts.}
Conjugation first sends a marked
cluster chart at \(Z_0\) to the chart at \(\bar Z_0\): it conjugates its
exterior marking, sends each chosen disc to its conjugate disc, and
reverses each positively oriented braid meridian. Cutting and replacing
a disc by a cone depends only on its boundary metric and total curvature,
both preserved by this operation. Conjugating before collapse therefore
gives exactly the conjugate of the collapsed sphere, with the conjugate
exterior marking. Composing this operation with a holomorphic chart arrow
gives every compatible odd arrow, by the compatible-coset description in
\cref{thm:twisted-inertia}. Apply the holomorphic stabilizer construction
to that second arrow. Its composite with conjugation descends to an
anti-M\"obius automorphism of the collapsed sphere.

Changing an interior marking changes the composite by a local collision
element. Conjugation normalizes the local collision groups (reversing
meridians and permuting the discs), so their kernel still acts trivially
on the collapsed exterior. Composition of two odd arrows agrees with
the holomorphic construction, since conjugation squares to the identity
on the marked charts. This proves that the extended stabilizer map is
well defined and respects composition. Equivariance on the completed
period space also follows by continuity: the conjugation isometries
agree on the dense open stratum by \cref{thm:twisted-inertia}, hence on
its metric completion. The finite collision kernel is retained before
specialization; boundary inertia is not identified with ordinary
configuration inertia.

For a double collision, if \(z_j\) is labelled separation, its
vanishing period is a unit times \(z_j^{1-4/N}\).  Since
\(t_j=z_j^2\) and \(q=2/(1-4/N)\) for the four packages, this gives
\(t_j=u_j\xi_j^q\), where \(u_j\) is a nonvanishing analytic unit.
Shrink the common period chart to choose a holomorphic \(q\)-th root
of each \(u_j\), and absorb these roots into the respective \(\xi_j\).
The choices can be made simultaneously: they concern finitely many
nonvanishing functions on the same small chart.  Disjoint discs give
independent separations and commuting meridians, so the equations hold
as a product chart, with the remaining coordinates tangent to the
collision stratum.  If cluster permutations occur, this chart is taken
before quotienting by those finite permutations.  Larger clusters use
the full local collision group, not this product of cyclic groups.

\emph{Specialization of geometric arrows.}
Choose configuration
representatives tending to \(Z_0\) in a local stable slice.  Normalize
their weights to total mass one.  Every atom of the limit then has mass
less than \(1/2\).  Anti-M\"obius maps preserving these configurations
cannot escape to infinity: otherwise their matrix representatives,
after composition with conjugation and projective normalization, have
a rank-one subsequential limit.  Away from its single exceptional
point, that limit maps everything to one point.  More than half the
configuration mass would then accumulate at that image, contradicting
stability.  Every subsequential limit is therefore a nondegenerate
anti-M\"obius automorphism of \(Z_0\).  To identify it, note that
stability forces at least three distinct limiting clusters: two atoms
of mass strictly less than \(1/2\) cannot carry the total mass one.
Choose one labelled point tending to each of three such cluster
centres.  The marked arrow \(\bar A\) sends these points into the
three target clusters prescribed by the collapsed automorphism
\(\sigma_0\).  Their source and target limits are distinct triples.
The limiting anti-M\"obius map must take the first triple to the
second, so it equals \(\sigma_0\).  This argument applies to every
subsequence, and compactness therefore gives convergence of the whole
family.  Internal braid changes may permute labels within a cluster
but cannot change its centre or the prescribed target cluster.  A
change of exterior marking conjugates both maps by the same map.
This proves specialization from every incident mirror-free component,
including arbitrary simultaneous stable collisions.
\end{proof}

\begin{lemma}[The angle type is constant on a fixed locus]
\label{lem:type-constant}
For an equal-weight package, under (DM) use
the completed period identification and specialization of
\cref{lem:stable-collision-specialization}.  Let \(\bar A\) be a compatible
projective class whose fixed locus meets \(\mathbb B^\circ\).  Then the
pairs \((Z,\sigma)\) attached to the mirror-free points of
\(\operatorname{Fix}(\bar A)\cap\mathbb B\) all have one angle type.  In
particular every chamber of \(\operatorname{Fix}(\bar A)\cap\mathbb B^\circ\)
has that type, and a sum over the compatible classes of one type counts
every chamber of their fixed loci exactly once.
\end{lemma}

\begin{proof}
Put \(F=\operatorname{Fix}(\bar A)\cap\mathbb B\), a connected real
hyperbolic space, and let \(e\) be the order of \(\bar A\).
By twisted inertia every mirror-free point of \(F\) carries an
anti-M\"obius map of exact order \(e\).  We first record why its angle
type cannot change in a nondegenerate limit; the argument uses only the
convergence of the maps.

Suppose \(\sigma_i\to\sigma_0\) are anti-M\"obius maps with
\(\sigma_i^e=1\).  Then \(\sigma_0^e=1\).  Regard them as isometries
of hyperbolic three-space through the Poincar\'e extension, which
identifies \(\mathbb P^1\) with the ideal boundary
\(\partial_\infty\mathbb H^3\), and choose a fixed point \(p_0\) of
\(\sigma_0\).  The unique minimizer \(p_i\) of
\[
 p\longmapsto\sum_{j=0}^{e-1}
       \operatorname{dist}(p,\sigma_i^j p_0)^2
\]
is fixed by \(\sigma_i\), by cyclic invariance and strict convexity.
Since all \(\sigma_i^j p_0\to p_0\), also \(p_i\to p_0\).
Orientation-preserving isometries tending to the identity carry \(p_i\)
to \(p_0\).  After conjugation the maps therefore lie in one compact
\(O(3)\) stabilizer.  Its orientation-reversing elements satisfying
\(g^e=1\) have finitely many eigenvalue multisets.  These determine their
\(SO(3)\)-conjugacy classes, namely the angle types of
\cref{thm:rotary}.  Continuity makes the eigenvalue multiset constant
near the limit.  Thus both the angle type and the exact order persist:
the eigenvalue multiset of \(\sigma_0\) equals that of \(\sigma_i\), whose
exact order is \(e\), so the exact order of \(\sigma_0\) is \(e\) and does
not drop to a proper divisor of \(e\).

At any \(x_0\in F\), \cref{lem:stable-collision-specialization} now
gives one limiting conjugacy class \(\sigma_0\) for every incident
mirror-free component.  The preceding argument makes all their angle
types equal.  The mirrors are locally finite and no mirror contains
\(F\), since \(F\) meets \(\mathbb B^\circ\).  Consequently their
complement is dense in \(F\), and the angle type extends locally
constantly across every collision point.  Connectedness of \(F\)
proves the assertion.  In particular, this argument covers simultaneous
collisions and intersections of real walls as well as single walls.
\end{proof}

\subsection{Conjugacy and uniqueness}
\begin{theorem}[Uniqueness criterion]
\label{thm:equal-weight-uniqueness}
For an equal-weight package, under (DM) use \cref{thm:one-partition}
and the completed period identification in
\cref{lem:stable-collision-specialization}, so that
\cref{lem:type-constant} applies.  For a fixed angle type,
\[
\begin{split}
  \#\{\text{admissible partitions}\}
  &=\sum_{\bar A}
    \#\bigl\{
      C(\bar A)\text{-orbits of components of }
      \operatorname{Fix}(\bar A)\cap\mathbb B^\circ
    \bigr\},
\end{split}
\]
the sum running over compatible projective \(P\Gamma_\mu\)-classes of that
type, with \(C(\bar A)=C_{P\Gamma_\mu}(\bar A)\).  Consequently if a single
constructed class realises as many component orbits as there are partitions,
that type carries exactly one compatible projective class.
\end{theorem}

\begin{proof}
By \cref{thm:one-partition} each partition is one connected configuration
space, by \cref{lem:type-constant} every chamber of the fixed locus of a
class of the given type has that type, and \(\iota_\mu\) identifies the moduli of pairs \((Z,\sigma)\) of
the given type with the disjoint union, over classes \(\bar A\), of
\((\operatorname{Fix}(\bar A)\cap\mathbb B^\circ)/C(\bar A)\).  Components
of the mirror-free fixed locus are the partitions.  Comparing cardinalities
gives the identity, and a matching lower bound from one class saturates it.
\end{proof}

\begin{lemma}[Unordered collision monodromy and odd-order wall conjugacy]
\label{lem:equatorial-merge}
Under (DM), use the completed period identification
of \cref{lem:stable-collision-specialization}, with equal weights
\(\mu=(2/N)^N\) in one of the four packages.  Let
\(\sigma=A_{n,k}\) have \(g_k=2\), and let \(A\) be a strict antiunitary
covering \(\sigma\).  Let a generic \(2n\)-orbit approach the equator,
producing \(n\) pairwise collisions
\(\pi_j=\bigl\{\sigma^j(p),\sigma^{j+n}(p)\bigr\}\) for
\(j=0,\dots,n-1\).  Write \(R_j\) for monodromy of an \emph{unordered}
collision meridian and put \(P=R_0R_1\cdots R_{n-1}\).  Its local
reflection order is
\[
 q=\frac{2N}{\gcd(N-4,2N)}=10,6,4,3
 \quad\text{for }N=5,6,8,12,
\]
respectively; since \((N-4)\mid8\) implies \((N-4)\mid2N\), this equals
\(2N/(N-4)=2/(1-4/N)\), the form used in the proof of
\cref{lem:stable-collision-specialization}.  The \(R_j\) commute,
\(APA^{-1}=P^{-1}\), and compatible
markings of the two real smoothings give antiunitaries \(A\) and
\(P^{-1}A\).  If \(q\) is odd these are conjugate: choose
\(j\) with \(2j\equiv-1\pmod q\), then
\(P^jAP^{-j}=P^{-1}A\).  The two partitions are
\((b,g,s)\) and \((b,g-1,s+2)\).  Among the four packages this odd-order
argument applies only to \(N=12\).

\end{lemma}

The single-collision specialization $P=R_0$ is the plane-reflection
case $A_{1,0}$, where the wall is an ordinary real double collision.
The restriction to $N=12$ concerns the odd-order conjugacy conclusion,
not the local monodromy calculation for the other packages.

\begin{figure}[htbp]
\centering
\includegraphics[width=0.5\textwidth]{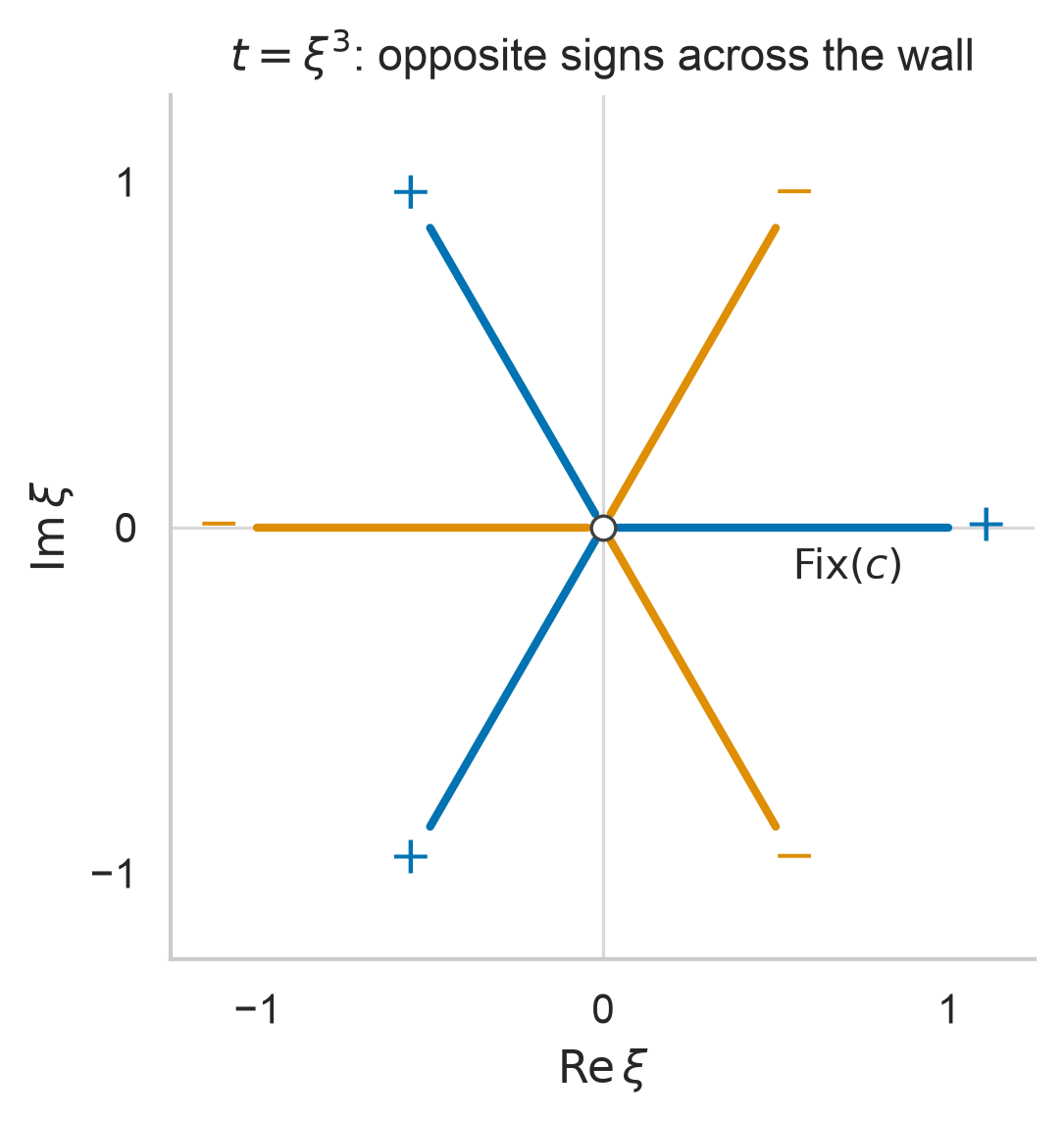}
\caption{The transverse period chart of \cref{lem:equatorial-merge}.
The rays of $\xi^q\in\mathbb R$ form the fixed lines of
$P^kc$, where $c(\xi)=\bar\xi$ and $P(\xi)=e^{2\pi i/q}\xi$.
Blue $+$ and amber $-$ denote the sign of the real smoothing parameter
$t=\xi^q$; the open origin is the collision. The figure uses
$\epsilon=+1$ in $t=\epsilon\xi^q$; at $q=3$, replacing $\xi$ by
$-\xi$ gives this convention when $\epsilon=-1$.
At $q=3$, the two sides of
one fixed line have opposite signs, and the conjugation
$P^jcP^{-j}=P^{2j}c$ joins the lifted real structures.}
\label{fig:counting-wall-crossing}
\end{figure}

\begin{proof}
The \(n\) pairs are disjoint and exhaust the orbit.  Half-twists supported
on disjoint collision discs commute, so their images \(R_j\) commute.
If \(z\) is the separation of the labelled points, the vanishing period
is a nonzero analytic unit times \(z^{1-4/N}\).  The unordered coordinate
is \(t=z^2\), so its meridian multiplies the period by
\(\exp(\pi i(1-4/N))\), of the displayed order.  A pure braid makes two
such loops and cannot replace this unordered meridian.

The power \(\sigma^n\) is equatorial reflection, and \(\sigma^j(p)\)
collides with \(\sigma^{j+n}(p)\) as the generating point meets the
equator.  Because \(\sigma\) is antiholomorphic it reverses orientation, so
it sends a positive meridian about \(\pi_j\) to a negative meridian about
\(\pi_{j+1}\).  Twisted inertia therefore gives
\(A R_j A^{-1}=R_{j+1}^{-1}\), indices modulo \(n\).  Hence
\[
  A\,P\,A^{-1}
  =\prod_{j=0}^{n-1}R_{j+1}^{-1}
  =P^{-1},
\]
with indices modulo \(n\), the last equality because the \(R_j\) commute.
Thus \(P^q=1\) and
\(P^jAP^{-j}=P^{2j}A\).

We first identify the two smoothings in the real transverse slice.
Choose a local coordinate \(w\) centred at one collision with
\(\sigma^n(w)=-\bar w\), and propagate it to the other discs by
\(\sigma\). The unordered pair \(\{w,-w\}\) has coordinate \(t=w^2\).
For positive real \(t\) its two points lie off the equator and are
interchanged by \(\sigma^n\); for negative real \(t\) they lie on the
equator and are fixed separately. Propagation around the \(n\) discs
therefore gives respectively one generic \(2n\)-orbit and two special
\(n\)-orbits.

The real normalization of the period coordinate can be made explicit.
In the product chart of \cref{lem:stable-collision-specialization},
take the transverse component of \(\operatorname{Fix}(\bar A^2)\).
The integer \(n\) is odd, so \(\bar A^2\) cyclically permutes all \(n\)
normal lines. Its \(n\)th power is the identity, and hence the product
of the multipliers around that cycle is one. Its fixed normal space
therefore has complex dimension one, with every coordinate determined
by the first. Linearizing the finite stabilizer at the collision point
gives a holomorphic transverse slice with this tangent line, while
keeping the tangential parameters fixed. The element \(\bar A^2\)
fixes this slice pointwise, and \(\bar A\) restricts to an
antiholomorphic involution. Choose a
holomorphic coordinate \(z\) for which this restriction is
\(z\mapsto\bar z\): if a starting coordinate has this linear part, replace
it by the average of itself and the conjugate of its pullback by \(A\).
Its derivative at the origin remains nonzero. On the configuration
slice the coordinate \(t=w^2\) has involution \(t\mapsto\bar t\).
Equivariance of the period map now gives
\[
 t(z)=z^q u(z),\qquad
 u(\bar z)=\overline{u(z)},\qquad u(0)\in\mathbb R^\times.
\]
Set \(\epsilon=\operatorname{sign}u(0)\). On a smaller invariant disc
the function \(\epsilon u\) has a holomorphic \(q\)-th root \(b\) with
\(b(0)>0\). Uniqueness of this root with that value gives
\(b(\bar z)=\overline{b(z)}\). Replacing \(z\) by
\(\xi=z b(z)\) thus preserves the real fixed line and gives the exact
relation \(t=\epsilon\xi^q\). This supplies the real coefficient rather
than inferring it from symmetry alone.

The two sides \(\pm\xi\) of this fixed line consequently give opposite
signs of \(t\) precisely when \(q\) is odd. In that case both partitions
\((b,g,s)\) and \((b,g-1,s+2)\) occur in adjacent chambers of the same
fixed locus, without changing its marking. When \(q\) is even, both
sides give the same smoothing.
Figure~\ref{fig:counting-wall-crossing} shows the fixed lines at $q=3$
and the opposite smoothing signs on the two sides of each line.

For completeness, the marked transition also fixes the inverse in
\(P^{-1}A\). Keep one exterior marking and continue the smoothing
parameter clockwise through a half-circle about zero in every collision
disc. A branch of its \(q\)-th root turns by \(-\pi/q\).
Reflection in the resulting line acts on a normal coordinate as
\[
 \xi\longmapsto e^{-2\pi i/q}\bar\xi,
\]
whereas before continuation it is \(\xi\mapsto\bar\xi\).
Since the positive unordered meridian acts by \(e^{2\pi i/q}\),
the new lift is \(R_j^{-1}\) times the old lift on the \(j\)th
normal direction. The exterior marking is unchanged; multiplying the
commuting local factors gives \(P^{-1}A\). Equivalently the
continuation path followed by the inverse of its conjugate is the
negative unordered meridian in each disc. Finally
\(P^jAP^{-j}=P^{2j}A=P^{-1}A\) when \(2j\equiv-1\pmod q\).
Such a \(j\) exists for odd \(q\), agreeing with the fixed-line
description above; for even \(q\), the other smoothing belongs to the
twisted fixed locus and this conjugacy argument does not apply.

For a plane reflection the colliding pair is a single conjugate pair (or
two real points), \(n=1\), and \(P\) is one complex reflection of order
\(q\).  The same braid and orientation computation applies.
\end{proof}

\begin{theorem}[Antipodal projective class]
\label{thm:antipodal-class}
Under (DM), the antipodal type exists if and only if
\(N\) is even, and then carries exactly one compatible projective
\(P\Gamma_\mu\)-class, with partition \([2^{N/2}]\).  It is not the split
class: the anti-M\"obius map has no fixed points on \(\mathbb P^1\).
\end{theorem}

\begin{proof}
\Cref{thm:equal-weight-spectrum} gives existence precisely for \(N\) even, by
placing \(N/2\) antipodal pairs.  There is one partition, so
\cref{thm:one-partition} gives uniqueness.  A plane reflection has a
circle of fixed points; \(-I\) has none.
\end{proof}

\subsection{The twelve-point classification}
\begin{theorem}[The twelve-point split projective class]
\label{thm:split-class}
Under (DM), use the completed period identification
of \cref{lem:stable-collision-specialization}, and let \(N=12\).  The
plane-reflection type carries exactly one compatible projective
\(P\Gamma_\mu\)-class.  It realises all seven partitions
\([2^{(12-r)/2},1^r]\), indexed by the number \(r\) of real points,
\(r\in\{0,2,4,\dots,12\}\).
\end{theorem}

\begin{proof}
Adjacent partitions differ by one conjugate pair meeting the real line, or
two real points colliding, which is the \(n=1\) case of
\cref{lem:equatorial-merge} with the odd local order \(q=3\).  The
resulting monodromy conjugates the
antiunitary of one chamber to that of the next.  The partition graph is a
path of length six, so all seven reflection partitions lie in
one class.  M\"obius maps commuting with conjugation preserve the real
line, hence preserve \(r\): different partitions are different
\(C(\bar A)\)-orbits of chambers.  \Cref{thm:equal-weight-uniqueness} saturates.
Existence is placement of \(r\) distinct real points and
\((12-r)/2\) conjugate pairs, which lie in \(M_{0,12}\) and hence in
\(\mathbb B^\circ\).
\end{proof}

\begin{theorem}[Complete twelve-point open class list]
\label{thm:open-class-list}
Under (DM), use the completed period identification
of \cref{lem:stable-collision-specialization}, and let \(N=12\).  Each
angle type of \cref{thm:equal-weight-spectra} that meets \(\mathbb B^\circ\)
carries exactly one compatible projective \(P\Gamma_\mu\)-class.  The list
\(\mathcal C(\mu)\) of these projective classes is complete, of cardinality
\(11\).  For \(N=5,6,8\), the numbers \(4,5,7\) count angle types
only; see the remark below.
\end{theorem}

\begin{proof}
Single-partition types are unique by \cref{thm:one-partition}, and are
realised by placement.  The reflection type is \cref{thm:split-class}.  The
antipodal type is \cref{thm:antipodal-class}.

The remaining multi-partition types are the special types \(g_k=2\) with
more than one solution of \eqref{eq:orbit}.  For \(n>1\) the value
of \(b\) is unique, so \(2g+s\) is unique, and the equatorial move
\((g,s)\mapsto(g-1,s+2)\) connects all partitions of the type in a single
path.  \Cref{lem:equatorial-merge} conjugates adjacent vertices.  M\"obius
maps commuting with \(\sigma\) preserve the equator, hence preserve whether
a point is equatorial, so distinct partitions remain distinct
\(C(\bar A)\)-orbits.  \Cref{thm:equal-weight-uniqueness} saturates.

The types and partition counts are those of
\cref{thm:equal-weight-spectra}.  For \(N=12\) the list is the special-case
eleven-class table.
\end{proof}

\begin{remark}[The smaller packages]
The one-class-per-type conclusion is specific to the odd collision order
\(q=3\) of \(N=12\).  For \(N=5\) the plane-reflection type carries three
projective antiunitary classes \cite[Lemma~7.2]{deGaayFortman2025}, and
the even local orders \(q=10,6,4\) of \(N=5,6,8\) give no wall-merging
argument; only the orbit equation and the per-partition connectedness of
\cref{thm:one-partition} carry over to those packages.
\end{remark}

\begin{definition}[Strict scalar fibres and witness masses]
\label{def:strict-scalar-fibres}
For any compatible projective class \(\bar A\), choose the
strict integral lift \(\chi_0\) furnished by the period local system and put
\[
 \mathcal L_{\Gamma_\mu}(\bar A)
 =\left\{[A]_{\Gamma_\mu}: A\in\Gamma_\mu\chi_0,
          \ [P(A)]_{P\Gamma_\mu}=\bar A\right\}.
\]
This is the finite \emph{strict scalar fibre} over \(\bar A\).  For a
strict antiunitary \(A\in\Gamma_\mu\chi_0\) and an index \(n>0\), with the
strict vector stabilizers \(\operatorname{Stab}_{C_{\Gamma_\mu}(A)}(v)\) of
\cref{lem:witness-formula,lem:primitive-strict-lifts}, define the strict
witness masses, the fixed-vector sets and the aggregates over the scalar
fibre by
\begin{align*}
 M_A^{\mathrm{open}}(n)
 &=\sum_{[v]\in C_{\Gamma_\mu}(A)\backslash X_{A,n}^{\mathrm{open}}}
       \frac1{|\operatorname{Stab}_{C_{\Gamma_\mu}(A)}(v)|},\\
 M_A^{\mathrm{amb}}(n)
 &=\sum_{[v]\in C_{\Gamma_\mu}(A)\backslash X_{A,n}^{\mathrm{amb}}}
       \frac1{|\operatorname{Stab}_{C_{\Gamma_\mu}(A)}(v)|},\\
 X_{A,n}^{\mathrm{open}}
 &=\{v\in L_\mu:Av=v,\ (v,v)/2=n,\ [v]\in\mathbb B^\circ\},\\
 X_{A,n}^{\mathrm{amb}}
 &=\{v\in L_\mu:Av=v,\ (v,v)/2=n\},\\
 \mathfrak M_{\bar A}^{\mathrm{open}}(n)
 &=\sum_{[A]\in\mathcal L_{\Gamma_\mu}(\bar A)}M_A^{\mathrm{open}}(n),\\
 \mathfrak M_{\bar A}^{\mathrm{amb}}(n)
 &=\sum_{[A]\in\mathcal L_{\Gamma_\mu}(\bar A)}M_A^{\mathrm{amb}}(n).
\end{align*}
\end{definition}
The aggregate is intrinsic.  Choosing another strict \(\chi_0\) only
relabels the same odd arrows.  A separate scalar-normalizer argument is
needed before replacing \(\mathfrak M_{\bar A}^{\mathrm{open}}\) by one
strict-lift mass.

\begin{remark}[Projective versus strict lifts]\label{rem:projective-vs-strict}
The twisted-inertia equivalence of \cref{thm:twisted-inertia} is projective.
Passing from projective pairs to strict antiunitary fixed-vector pairs is a
separate arithmetic scalar-fibre calculation, and the stabilizer
isomorphism of \cref{thm:twisted-inertia} is exact; for the fullerene
package that calculation is \cref{lem:period-bridge}(e) and
\cref{lem:fullerene-scalar-fibre-collapse}.
\end{remark}

The real-period correspondence has identified the projective symmetry
types and preserved their stabilizers. To use the witness formula on
integral vectors, the next step is to choose strict representatives and
check their scalar lifts.

\part{Real periods and strict symmetry classes}
\section{Antiunitary representatives and involutive gluing}
\label{sec:fixed-loci}
\subsection{Integral representatives and involutive glue}
\label{sec:integral-representatives}
We use the explicit integral normalization of \inv{sec:common-lattice},
repeated here to identify the period objects without a notation lookup.
The first ten roots $\alpha_0,\ldots,\alpha_9$ form an Eisenstein basis.
Its Hermitian Gram $H_{\mathcal E}$ is tridiagonal, with
\[
 (H_{\mathcal E})_{ii}=-3,\qquad
 (H_{\mathcal E})_{i,i+1}=1+\zeta,\qquad
 (H_{\mathcal E})_{i+1,i}=1+\bar\zeta.
\]
The real bilinear form is
$B(x,y)=\frac23\operatorname{Re}(x^*H_{\mathcal E}y)$.
In the ordered integral basis
$(\alpha_0,J\alpha_0,\ldots,\alpha_9,J\alpha_9)$ it is even
unimodular of signature $(2,18)$; the tridiagonal determinant recurrence
verifies these assertions.  The two remaining roots are determined by
\[
 \sum_{i=0}^{11}\alpha_i=0,\qquad
 \sum_{i=0}^{11}\bar\zeta^{\,i}\alpha_i=0.
\]
The coefficient determinant $\bar\zeta^{10}-\bar\zeta^{11}$ is a unit,
so these relations are integral.  This fixes every coordinate and
operator normalization used below.

For a collision root $\alpha_j$, with $h(x,y)=x^*H_{\mathcal E}y$
conjugate-linear in the first variable and $h(\alpha_j,\alpha_j)=-3$, the
order-three \emph{root triflection} is the unitary map
\[
 T_j(x)=x+(1-\zeta^2)\,\frac{h(\alpha_j,x)}{3}\,\alpha_j
       =x-(1-\zeta^2)\,\frac{h(\alpha_j,x)}{h(\alpha_j,\alpha_j)}\,\alpha_j,
\]
which fixes $\alpha_j^\perp$ pointwise and multiplies $\alpha_j$ by
$\zeta^2$; it is integral because every entry of $H_{\mathcal E}$ lies in
the prime ideal $(1+\zeta)\mathbb Z[\zeta]$ above $3$.  The same symbol $T_j$
is used for these operators throughout this section: $T_{2i}$ in
\cref{prop:real-top}, $T_j$ with $j$ even in \cref{lem:q4-merge}, and
$T_0,\dots,T_9$ in the order-ten construction at the end of this
subsection.

\paragraph{Chain marking and monodromy convention.}
Put $\theta=1+\zeta$.  For developed cone points $p_0,\ldots,p_{11}$,
choose the side of the cut chain for which the holonomies
$g_j(w)=\zeta(w-p_j)+p_j$ satisfy $g_0\cdots g_{11}=1$, with the
rightmost map applied first.  Equivalently, $\sum_j\zeta^jp_j=0$.
Normalize the period vector by
\[
 h(\alpha_j,v)=\theta(p_{j+1}-p_j),\qquad j\in\mathbb Z/12\mathbb Z.
\]
The pairing order matters: $v$ is in the linear slot.  In the right
change-of-marking convention, the counterclockwise Hurwitz exchange is
\[
 (g_j,g_{j+1})\longmapsto(g_jg_{j+1}g_j^{-1},g_j),\qquad
 (p_j,p_{j+1})\longmapsto
 ((1-\zeta)p_j+\zeta p_{j+1},p_j).
\]
Writing $d_j=p_{j+1}-p_j$, the adjacent differences become
$d_{j-1}+\zeta d_j$, $-\zeta d_j$, and $d_{j+1}+d_j$.
The Gram matrix therefore gives $v'=v+\zeta d_j\alpha_j=T_j^{-1}v$;
the inverse re-marking gives $T_jv$.
We use the inverse, left action for continuation of periods: period
coordinates pull back by the inverse marking, as in
\cref{thm:twisted-inertia}.  Thus the positive unordered meridian of
\cref{lem:equatorial-merge}, whose holomorphic normal coordinate is
multiplied by $e^{2\pi i/3}=\zeta^2$, acts by $R_j=T_j$.

The unitary orbit is independent of the cut tree. Two embedded cut
trees on the same punctured sphere are related by branch slides and
changes of basepoint: put them in transverse position, slide along
outermost segments to remove crossings, and then match their leaf
branches successively. A branch slide changes consecutive based
peripheral loops by the Hurwitz exchange above or its inverse, hence
changes the period by $T_j^{-1}$ or $T_j$. Isotopies do not change
the developed periods, and changing the base development rotates them
by a unit of $\mathcal E$; a common translation cancels in their
differences. All these operations therefore preserve the
$U(\Lambda)$-orbit used in \cref{prop:fullerene-period-correspondence}.

On the twelve collision-root coordinates $\alpha_0,\ldots,\alpha_{11}$,
define the integral antiunitary actions
\begin{equation}\label{eq:phase-ramp}
 Q_{t,k}\colon\alpha_j\longmapsto\zeta^{k-j}\alpha_{j+t},
 \qquad t\in\mathbb Z/12\mathbb Z,\quad k\in\mathbb Z/6\mathbb Z.
\end{equation}
Their squares send $\alpha_j$ to $\zeta^{-t}\alpha_{j+2t}$, so their
exact orders are $24/\gcd(2t,12)$.  Steps $t$ and $-t$ are conjugate,
and a unitary rotation removes the phase $k$ inside the root-monomial
normalizer.  Write $Q_t=Q_{t,0}$ and $F_t=\Lambda^{Q_t}$ for
$t=0,\ldots,6$.  Their orders are $2,12,6,4,6,12,2$.
\Cref{lem:dodecagon-model} identifies $Q_0$ with the plane-reflection
(split) type and $Q_6$ with the antipodal type.
The order-ten representatives are constructed at the end of this
subsection.  The
classification below concerns classes meeting the mirror-free locus;
torsion contained in mirrors is irrelevant to the count.

The next proposition is the concrete integral root-triflection
realization, on $\Lambda$, of the braid-monodromy merge of
\cref{thm:split-class}: the unordered collision monodromies $R_j$ of
\cref{lem:equatorial-merge}, of local order $q=3$, are the triflections
$T_j$ in this left-continuation convention.

\begin{proposition}[Split wall conjugacy and arithmetic centralizers]
\label{prop:real-top}
Put $\chi_0=Q_0$ and $\chi_{i+1}=T_{2i}^{-1}\chi_i$ for $0\leq i<6$, with
$T_{2i}$ the root triflections.  Then
\[
 T_{2i}\chi_iT_{2i}^{-1}=\chi_{i+1}
 \quad(0\leq i<6),
\]
so $\chi_0,\dots,\chi_6$ are conjugate integral antiunitary involutions;
they are the integral form on $\Lambda$ of the six wall crossings of
\cref{thm:split-class}, with $T_{2i}$ the collision monodromy of
\cref{lem:equatorial-merge}.  Their fixed lattice is $2I_{1,9}$.  The
antipodal representative $Q_6$ has fixed lattice $2II_{1,9}$, so it is not
conjugate to the split class.

For either class, put $F=\Lambda^\chi$, $F^-=\Lambda^{-\chi}$,
$D=2J-1$, and $K=\ker(D\bmod3)\subset F/3F$.  Then $K$ is maximal totally
isotropic in the ten-dimensional quadratic space $F/3F$ over
$\mathbb F_3$, which has plus type; write $O^+_{10}(\mathbb F_3)$ for its
orthogonal group, the superscript referring to the Witt type and not to a
positive cone.  Restriction identifies the strict arithmetic centralizer
with
\[
 C_{U(\Lambda)}(\chi)
 \cong
 \operatorname{red}_3^{-1}
 \bigl(\operatorname{Stab}_{O^+_{10}(\mathbb F_3)}K\bigr).
\]
The integral action is transitive on these states, and the centralizer has index
\[
 2\prod_{j=1}^{4}(3^j+1)=183680
\]
in $O(I_{1,9})$ for the split class and in $O(II_{1,9})$ for the antipodal
class. This is the number of maximal isotropic five-dimensional subspaces
of $F/3F$.
\end{proposition}

\begin{proof}
The even roots $\alpha_0,\alpha_2,\dots,\alpha_{10}$ are pairwise
orthogonal, so their triflections commute, and $Q_0T_jQ_0^{-1}=T_j^{-1}$
because $Q_0$ preserves each root line antilinearly.  Hence
$\chi_iT_{2i}\chi_i^{-1}=T_{2i}^{-1}$ and
$T_{2i}\chi_iT_{2i}^{-1}=T_{2i}^{-1}\chi_i=\chi_{i+1}$, which proves the
wall conjugacies.  The positive meridian fixes the collision hyperplane
and has normal eigenvalue $\zeta^2$, so it is $T_{2i}$ rather than
$T_{2i}^{-1}$ in the continuation convention fixed above.
The fixed lattices, the strict extension criterion, the transitivity on
ternary states and the index are \inv{prop:involutive-lattice-glue}.
That proposition states the index inside the positive-cone orthogonal
group; since $-I$ commutes with $\chi$ and lies in the centralizer, and
$O(F)=O^+(F)\times\{\pm I\}$, the index in the full group $O(F)$ is the
same number.
Their different parities also distinguish the split and antipodal classes.
\end{proof}

The four order-ten representatives $A_1,\dots,A_4$ of
\rot{thm:order10-golden-centralizer}, one for each angle parameter
$k=1,2,3,4$, are constructed from
\[
 \varepsilon=T_0(T_0T_1\cdots T_9),\qquad
 g=\zeta\varepsilon^2,
\]
where $|\varepsilon|=60$, $\varepsilon^{10}=\zeta$, and $|g|=5$.
The $g$-fixed lattice has rank four, signature $(2,2)$, determinant $625$,
and level five; its rank-sixteen negative complement has determinant $625$
and minimum four.  For angle parameter $k=1,2,3,4$, the fixed lattice of
$A_k$ is $10U$ of level ten when $k$ is odd and $10I_{1,1}$ of level twenty
when $k$ is even. These are lattices for the unhalved bilinear form $B$:
the index $m=B(v,v)/2$ is $10xy$ on $10U$ and $5(x^2-y^2)$ on
$10I_{1,1}$. With the half-form $B/2$ used in Paper~III \cite{fullereneRotaryPaper}, the same
lattices are $5U$ and $5I_{1,1}$. The four $A_k$ are pairwise nonconjugate in $U(\Lambda)$,
and every class contains mirror-free positive lines.  Their full strict
centralizers are all the cyclic group $C_{10}$ of order ten.  The $10U$
plane has no positive global root
wall; the $10I_{1,1}$ plane has exactly the two primitive wall rays of
quadratic values $5$ and $15$.  These integral invariants, the strict
centralizers and the pairwise nonconjugacy are proved in
\rot{thm:order10-golden-centralizer}, and the wall list in
\rot{thm:order10-analytical-walls}.

\subsection{The fullerene period correspondence}

The two completeness statements below are transported from the
Deligne--Mostow picture by a single bridge.  We isolate it because its parts
have different inputs: the holomorphic orbifold isomorphism is quoted,
whereas its extension to orientation-reversing arrows follows from the
natural conjugation on the cyclic-cover Gauss--Manin local system.  The
relevant groupoid is twisted inertia, with no involutive cocycle condition;
this distinction is what retains anti-M\"obius elements of orders
$4,6,10,12$.  Parts \textup{(b)}--\textup{(e)} are the specialization to
$N=12$, $d=6$ of \cref{thm:twisted-inertia,lem:primitive-strict-lifts},
restated in the notation of the fullerene lattice $\Lambda$ with
$\Gamma_\mu=U(\Lambda)$.

\begin{lemma}[Real-compatible Deligne--Mostow bridge]\label{lem:period-bridge}
Let $\mu_1=\dots=\mu_{12}=1/6$, so $\sum_i\mu_i=2$ and
$(1-\mu_i-\mu_j)^{-1}=3/2$ satisfies the half-integral condition
$\Sigma\mathrm{INT}$.  For twelve distinct points
$Z=\{z_1,\dots,z_{12}\}\subset\mathbb P^1$ let $C_Z$ be the cyclic cover
$y^6=\prod_i(x-z_i)$.
\begin{enumerate}[label=(\alph*)]
  \item \emph{(External input.)}  The period map of $C_Z$ induces an
        isomorphism of orbifolds
        \[
          M_{0,12}/S_{12}\;\cong\;\mathbb B^\circ/PU(\Lambda),
        \]
        and the projectivizations of the ordinary and affine braid
        triflection representations have image $PU(\Lambda)$, which
        follows from \cite[Theorems~5.1(a) and~5.2, with the full isometry
        group defined in \S2.7]{allcockBasak2025}; see also the displayed period isomorphism
        and braid description in \cite[\S\S2,4]{heckman2015}.
        The loop-to-braid-generator comparison is
        made explicit in \cite[Thm.~5.2]{allcockBasak2025}, and the underlying
        cone-metric description is \cite[Thm.~0.1]{thurston1998}.
        Its extension over the entire stable collision locus is the completed
        identification used in
        \cref{lem:stable-collision-specialization}.  Only this projective
        statement is used here.
  \item \emph{(Cover-level calculation.)}  Conjugation
        $(x,y)\mapsto(\bar x,\bar y)$ is an antiholomorphic isomorphism
        $C_Z\to C_{\bar Z}$ intertwining the $\zeta$-eigenspace
        decompositions of the two covers through $\zeta\mapsto\bar\zeta$;
        it is the antilinear isometry $K_Z$ of \cref{thm:twisted-inertia}.
  \item The isomorphism in \textup{(a)} is equivariant for complex
        conjugation and induces an equivalence of the conjugation-twisted
        inertia groupoids.  Thus a pair $(A,[v])$, with $A$ an antiunitary
        isometry of $\Lambda$ and $[v]\in\mathbb B^\circ$ fixed by $A$,
        corresponds to a pair $(Z,\sigma)$ of twelve distinct points and an
        anti-M\"obius transformation preserving them, and the
        correspondence is a bijection between $PU(\Lambda)$-classes of
        pairs $(\bar A,[v])$ and $\mathrm{PGL}_2(\mathbb C)$-classes of
        pairs $(Z,\sigma)$.
  \item Exact orders are preserved: if $A$ has exact order $e$ and fixes
        some $v\neq0$, then $\bar A$ has exact order $e$.
  \item \emph{(Primitive pointed lifts.)}  On rational lines, sending
        $(\bar A,[v])$ to the unique lift of $\bar A$ fixing a chosen
        primitive vector $v$ of $[v]$ is a bijection onto the
        $U(\Lambda)$-classes of pairs $(A,v)$ with $A$ antiunitary and
        $v\in\Lambda$ primitive, positive and mirror-free, $Av=v$; it
        identifies the projective automorphism group of $(\bar A,[v])$
        with $\operatorname{Stab}_{C_{U(\Lambda)}(A)}(v)$.
\end{enumerate}
\end{lemma}

\begin{proof}
Part~\textup{(a)} is the quoted Deligne--Mostow theorem.  Parts
\textup{(b)}--\textup{(d)} are \cref{thm:twisted-inertia} for the
fullerene package $\mu=(1/6)^{12}$, $d=6$,
$\mathcal O_6=\mathbb Z[\zeta]$, with part~\textup{(a)} as its
hypothesis; the extension $\tilde\rho$ of the monodromy to the odd coset
reads
\begin{equation}\label{eq:real-monodromy}
 \bar\rho(c_*\beta)=
 \bar\chi_0\,\bar\rho(\beta)\,\bar\chi_0^{-1}
 \qquad\text{in }PU(\Lambda)
\end{equation}
for every orbifold braid $\beta$, and on the spherical braid generators it
says that conjugation reverses every positive half-twist while
$\bar\chi_0$ sends its triflection to the inverse triflection.  The same
framed real-period mechanism is used for adjacent ball quotients in
\cite[Thm.~3.3 and Lemmas~3.4--3.5]{allcock2010} and
\cite[\S\S3.3--3.6]{chu2011}.  For the order statement, if
$A^{e'}=\varepsilon I$ for a proper divisor $e'$ of $e$ and an Eisenstein
unit $\varepsilon$, then $Av=v$ gives $v=A^{e'}v=\varepsilon v$, so
$\varepsilon=1$, contradicting exactness.

Part~\textup{(e)} is \cref{lem:primitive-strict-lifts} with $d=6$ and
$\Gamma_\mu=U(\Lambda)$, which contains the scalar isometries $\mu_6$ and
is the full inverse image of $PU(\Lambda)$.  The one step specific to
$d=6$ is that the map is well defined on $U(\Lambda)$-classes.  Suppose $g\in U(\Lambda)$ carries
$(\bar A,[v])$ to $(\bar A',[v'])$, with $Av=v$ and $A'v'=v'$.  Then
$gv=\zeta^bv'$ for some $b\in\mathbb Z/6\mathbb Z$, and $gAg^{-1}$ is a
lift of $\bar A'$ fixing $\zeta^bv'$; since
$\zeta^jA'(\zeta^bv')=\zeta^{j-b}v'$, that forces $gAg^{-1}=\zeta^{2b}A'$.
Replacing $g$ by $\zeta^{-b}g$, which is again in $U(\Lambda)$, and using
antilinearity to move the scalar across, gives
$(\zeta^{-b}g)A(\zeta^{-b}g)^{-1}=\zeta^{-2b}\zeta^{2b}A'=A'$.  Thus
projectively conjugate pointed pairs have $U(\Lambda)$-conjugate strict
lifts; the inverse map is projection, so the map on classes is a
bijection.
\end{proof}

\subsection{Projective class uniqueness and the order-six merge}
The regular dodecagon realizes every root-monomial representative at once.

\begin{lemma}[The dodecagon model]\label{lem:dodecagon-model}
Let $Z_{12}=\{\omega_j=e^{\pi ij/6}:j\in\mathbb Z/12\mathbb Z\}$ be the
regular dodecagon on the unit circle and, for $t\in\mathbb Z/12\mathbb Z$,
let
\[
 \sigma_t(z)=e^{\pi it/6}/\bar z .
\]
Then $\sigma_t$ is an anti-M\"obius transformation preserving $Z_{12}$,
acting on its points by $\omega_j\mapsto\omega_{j+t}$, and the pair
$(Z_{12},\sigma_t)$ corresponds under \cref{lem:period-bridge}(c) to the
projective class of $Q_t$.  In particular:
\begin{enumerate}[label=(\roman*)]
 \item $\sigma_0(z)=1/\bar z$ is the reflection in the unit circle, which
       it fixes pointwise, so $Q_0$ represents the plane-reflection type of
       \cref{thm:split-class}, in the chamber with $r=12$ real points;
 \item $\sigma_6(z)=-1/\bar z$ is the antipodal map, without fixed points,
       so $Q_6$ represents the antipodal type of \cref{thm:antipodal-class};
 \item $\sigma_t$ is the rotary reflection $A_{n,k}$ with $k/n=t/6$, so
       $Q_3$, $Q_2$, $Q_4$ and $Q_1,Q_5$ represent the types
       $(2,1)$, $(3,1)$, $(3,2)$ and $(6,1),(6,5)$ of
       \cref{thm:equal-weight-spectra}, with orbit partitions
       $[4^3]$, $[6^2]$, $[3^4]$ and $[12]$.
\end{enumerate}
Since each type carries one compatible projective class
(\cref{thm:open-class-list}), $\bar Q_t$ is that class.
\end{lemma}

\begin{proof}
Since $\bar\omega_j=\omega_j^{-1}$, one has
$\sigma_t(\omega_j)=e^{\pi it/6}\omega_j=\omega_{j+t}$, and $\sigma_t$ is
the composite of the reflection $z\mapsto1/\bar z$ in the unit circle with
the rotation by $\pi t/6$ about the polar axis, which is the normal form
$A_{n,k}$ with $\pi k/n=\pi t/6$.  Its square is the rotation by $\pi t/3$,
so its order is $2$, $12$, $6$, $4$, $6$, $12$, $2$ for $t=0,\dots,6$,
matching the orders of $Q_t$, and the orbit sizes on $Z_{12}$ are
$12/\gcd(t,12)$, which gives the displayed partitions.  Parts (i) and (ii)
are the cases $t=0$ and $t=6$.

For the identification with $Q_t$, mark $Z_{12}$ as in
\cref{lem:period-bridge}(a). Here a marking is a choice over the
orbifold point $Z_{12}$, whose dihedral isotropy group has order $24$;
it is not a choice of a point with trivial stabilizer.
The based orbifold identification of
\cite[Theorem~5.2 and Lemma~5.4]{allcockBasak2025} retains that isotropy.
Changing the marking conjugates the resulting operators.
The twelve counterclockwise half-twists $\beta_j$ along the
arcs from $\omega_j$ to $\omega_{j+1}$ are the affine braid generators of
\cite[Thm.~5.2]{allcockBasak2025}.  Their left continuation action on periods
is the triflection $T_j$ in the collision root $\alpha_j$; their right
Hurwitz re-marking action is $T_j^{-1}$, as fixed in
\cref{sec:integral-representatives}.  The
rotation $z\mapsto e^{\pi i/6}z$ of $Z_{12}$, which conjugates $\beta_j$
to $\beta_{j+1}$, acts as the unitary $R\colon\alpha_j\mapsto\alpha_{j+1}$,
with $RT_jR^{-1}=T_{j+1}$.  The map $\sigma_t$ carries the arc from
$\omega_j$ to $\omega_{j+1}$ onto the arc from $\omega_{j+t}$ to
$\omega_{j+t+1}$ and reverses the orientation of the sphere, so
$\sigma_t\beta_j\sigma_t^{-1}=\beta_{j+t}^{-1}$.  By the equivariance
\eqref{eq:real-monodromy}, the antiunitary $A$ attached to
$(Z_{12},\sigma_t)$ by \cref{lem:period-bridge}(c) therefore satisfies
$AT_jA^{-1}=T_{j+t}^{-1}$ in $PU(\Lambda)$; a scalar multiple of
$T_{j+t}^{-1}$ that is again a complex reflection of order three is
$T_{j+t}^{-1}$ itself, so the identity holds in $U(\Lambda)$.  Comparing
the root lines $\operatorname{im}(T_j-1)=\mathcal E\alpha_j$ shows that
$A$ maps $\mathcal E\alpha_j$ onto $\mathcal E\alpha_{j+t}$, hence
$A\alpha_j=u_j\alpha_{j+t}$ with units $u_j$, because $A$ is an integral
isometry and the roots are primitive.  Antiunitarity and the cyclic Gram
$h(\alpha_j,\alpha_{j+1})=1+\zeta$ of \cref{sec:integral-representatives}
give $\bar u_ju_{j+1}(1+\zeta)=\overline{1+\zeta}=\bar\zeta(1+\zeta)$, so
$u_{j+1}=\bar\zeta u_j$ and $A=Q_{t,k}$ with $\zeta^k=u_0$, in the
notation of \cref{eq:phase-ramp}.  Finally $RQ_{t,k}R^{-1}=Q_{t,k+1}$, so
every $Q_{t,k}$ is conjugate to $Q_t=Q_{t,0}$, and the pair corresponds to
the projective class of $Q_t$.  The configuration $Z_{12}$ consists of
twelve distinct points, so its period point is mirror-free and
$\bar Q_t$ is a compatible class of the angle type of $\sigma_t$; by
\cref{thm:open-class-list} it is the only one.
\end{proof}

\begin{remark}[A parity check of the identification]
\label{rem:antipodal-parity}
A fixed-point-free involution of a fullerene graph on $2m$ vertices acts
freely on its $3m$ edges, since an edge mapped to itself would contain a
fixed midpoint, so $m$ is even.  This agrees with the identification: the
row of $Q_6$ is supported on even $m$, because its fixed lattice
$2II_{1,9}$ is even and the row is evaluated at $m/2$
(\cref{cor:fixed-dyadic-symmetry-shares}), whereas the row of $Q_0$ is
positive at odd indices, for instance $a_{15}=3$ and
$M^{\rm open}_{\rm split}(15)=13/5$ by the formula of
\cref{thm:exact-all-index-grand-formula}.
\end{remark}

The scalar-fibre lemma and the counting theorem use the class list of
\cref{thm:open-class-list} in the following form.

\begin{lemma}[Projective per-type uniqueness]
\label{lem:type-uniqueness}
For the fullerene package, \(\Gamma_\mu=U(\Lambda)\), each angle type of
\cref{thm:equal-weight-spectra} carries exactly one compatible projective
\(PU(\Lambda)\)-class (\cref{thm:open-class-list}), and
\cref{thm:equal-weight-uniqueness} reads
\begin{equation}\label{eq:merge-count}
\begin{split}
 \#\{&\text{admissible partitions of the type}\}\\
 &=\sum_{\bar A}\#\bigl\{
 C_{PU(\Lambda)}(\bar A)\text{-orbits of components of }
 \operatorname{Fix}(\bar A)\cap\mathbb B^\circ\bigr\}.
\end{split}
\end{equation}
Thus the projective centralizers have seven component orbits for the
split type, three for \(Q_4\), two for each even order-ten type, and one
for each of the other seven types.  These assertions concern projective
classes; the passage to strict classes for all integral witnesses is
\cref{lem:fullerene-scalar-fibre-collapse}.
\end{lemma}

\begin{proof}
\Cref{lem:stable-collision-specialization} supplies the completed period
identification for the fullerene package, so
\cref{thm:open-class-list,thm:equal-weight-uniqueness} apply; the
partition counts are those of \cref{thm:equal-weight-spectra}.  The next
lemma gives explicit marked conjugators for the three \(Q_4\) partitions.
\end{proof}

The next lemma is the concrete integral root-triflection realization of
the equatorial merge of \cref{lem:equatorial-merge} for the $Q_4$ type,
with $n=3$ and $q=3$: the products $P_0,P_2$ below are the commuting
products $P=R_0R_1\cdots R_{n-1}$ of that lemma, written in the
triflections $T_j$ of \cref{sec:integral-representatives}.

\begin{lemma}[Simultaneous-wall merge for $Q_4$]\label{lem:q4-merge}
Let $T_j$ be the order-three triflection in the collision root $\alpha_j$,
put $A_0=Q_4$, and define
\[
 P_0=T_0T_4T_8,\qquad P_2=T_2T_6T_{10},\qquad
 A_1=P_0^{-1}A_0,\qquad A_2=P_2^{-1}A_1.
\]
Then $A_0,A_1,A_2$ are conjugate exact order-six antiunitaries.  They occur
in mirror-free chambers whose anti-M\"obius orbit partitions are,
respectively,
\[
 [3,3,3,3],\qquad[6,3,3],\qquad[6,6].
\]
Consequently the $Q_4$ angle type carries exactly one compatible projective class and
its projective centralizer has exactly three chamber orbits.
\end{lemma}

\begin{proof}
The phase-ramped action gives
$A_0T_jA_0^{-1}=T_{j+4}^{-1}$.  The six roots with indices
$0,4,8,2,6,10$ are pairwise orthogonal, so their triflections commute.
Hence $P_i^3=1$, $A_0P_iA_0^{-1}=P_i^{-1}$, and $P_0P_2=P_2P_0$.  Therefore
\[
 P_0A_0P_0^{-1}=P_0^{-1}A_0=A_1,
 \qquad
 P_2A_1P_2^{-1}=P_2^{-1}A_1=A_2.
\]
Conjugacy proves antiunitarity and exact order six.

The geometric wall labels can be read off the dodecagon model of
\cref{lem:dodecagon-model}.  Put $\omega=e^{2\pi i/3}$ and
$\sigma(z)=\omega/\bar z=\sigma_4(z)$.  A special orbit is
$S(a)=\{a,\omega a,\omega^2a\}$ with $|a|=1$, while
\[
 G(r,a)=\{r\omega^ja,r^{-1}\omega^ja:0\leq j<3\}
\]
is a generic six-orbit.  As $r\to1$, $G(r,a)$ becomes three double points
on $S(a)$; the other real smoothing is
$S(ae^{i\delta})\cup S(ae^{-i\delta})$.  In a local logarithmic coordinate
$w$, $\sigma^3$ acts by $w\mapsto-\bar w$, so the unordered smoothing
coordinate $w^2$ has opposite signs on these two smoothings.  At the regular
dodecagon the four special orbits are the residue classes modulo four.  The
first smoothing collides the pairs indexed by $\{0,4,8\}$ and the second
those indexed by $\{2,6,10\}$.

Equivalently, on point labels let $q(j)=j+4$ and
\[
 s_0=(0,1)(4,5)(8,9),\qquad
 s_2=(2,3)(6,7)(10,11).
\]
Then
\[
\begin{aligned}
q&=(0,4,8)(1,5,9)(2,6,10)(3,7,11),\\
s_0q&=(0,5,8,1,4,9)(2,6,10)(3,7,11),\\
s_2s_0q&=(0,5,8,1,4,9)(2,7,10,3,6,11),
\end{aligned}
\]
which gives the three displayed partitions. Choose the phases of the
remaining special orbits distinct from the colliding phases. The limiting
configuration then has exactly three disjoint double collisions and no
others. By \cref{lem:stable-collision-specialization}, a sufficiently
small period chart has just the three corresponding collision factors
$t_j=\xi_j^3$. Its two indicated real smoothings are mirror-free.
Thus the geometric construction itself supplies the required wall
crossings and excludes an additional collision at the crossing.

The displayed strict representatives are conjugate, so their three
partition strata realise one projective class.  By
\cref{eq:merge-count}, its projective-centralizer summand is at
least three, while the total over all classes is the three admissible
partitions.  Equality follows and no additional $Q_4$ class exists.
\end{proof}

\subsection{Scalar lifts on primitive and imprimitive shells}
\begin{lemma}[Fullerene scalar-fibre collapse]
\label{lem:fullerene-scalar-fibre-collapse}
Let $\mathcal E=\mathbb Z[\zeta]$, $\Gamma=U(\Lambda)$, and let $A$ be
one of $Q_0,\ldots,Q_6$ or the four representatives $A_1,\ldots,A_4$
of \rot{thm:order10-golden-centralizer}.  Every strict antiunitary projectively $P\Gamma$-conjugate to $A$ is strictly
$\Gamma$-conjugate to $A$.  Thus each of the eleven compatible projective
types has one strict scalar fibre.  This remains true for imprimitive fixed
vectors, and their masses are still formed from strict-centralizer orbits
with strict vector stabilizers.
\end{lemma}

\begin{proof}
It suffices to find, for each representative $A$, an integral unitary
$g_A$ with $g_AAg_A^{-1}=\zeta A$. Since $g_A$ commutes with scalars,
its powers give $g_A^jAg_A^{-j}=\zeta^jA$ for every integer $j$,
covering all six scalar twists.

For the seven root-monomial representatives $Q_t$, $0\leq t\leq6$, let
$R$ be the unitary rotation $R\alpha_j=\alpha_{j+1}$ of the twelve
collision-root coordinates, as in \cref{lem:dodecagon-model}.  From
$Q_{t,0}\alpha_j=\zeta^{-j}\alpha_{j+t}$ one obtains
\[
 RQ_{t,0}R^{-1}\alpha_j
 =\zeta^{1-j}\alpha_{j+t}
 =\zeta Q_{t,0}\alpha_j.
\]
This covers the split and antipodal types, $Q_2,Q_3,Q_4$, and the two
order-twelve types $Q_1,Q_5$.

For the odd-angle order-ten representatives, the unitary
$\varepsilon=T_0(T_0T_1\cdots T_9)$ of \cref{sec:integral-representatives},
of order $60$ with $\varepsilon^{10}=\zeta$, satisfies
\begin{equation}\label{eq:order10-odd-twist}
 \varepsilon^{-1}A_1\varepsilon=\zeta A_1,
 \qquad
 \varepsilon^{-1}A_3\varepsilon=\zeta A_3 .
\end{equation}
The first identity is a multiplication of the integral matrix $C_1$ of
$A_1$, printed in \cref{rot-sec:order10-analytical-reduction}
of~\cite{fullereneRotaryPaper}, with the triflections; the second follows
from $A_3=A_1^7$ and $(\zeta A_1)^7=\zeta A_1^7$, because
$(\zeta A_1)^2=\zeta\bar\zeta A_1^2=A_1^2$ by antilinearity.
For the even-angle representatives, put $V=A_1A_2$.  This is an
integral unitary operator: in the same Eisenstein coordinates its matrix
is $C_1\overline{C_2}$.  Define
\begin{equation}\label{eq:order10-even-twist}
 W=\varepsilon V\varepsilon V^{-1}\varepsilon,
 \qquad h=\zeta^{-1}W.
\end{equation}
Multiplication using the displayed matrices $C_1,C_2$ gives
$WA_2=-A_2W$.  Thus, remembering that $A_2$ is antilinear,
\[
 hA_2h^{-1}=\zeta^{-2}(-A_2)=\zeta A_2,
 \qquad hA_4h^{-1}=(\zeta A_2)^3=\zeta A_4.
\]
Here $A_4=A_2^3$ and $\zeta^3=-1$.  All factors in
\cref{eq:order10-even-twist} are integral unitary, so this is a strict
conjugator on the whole lattice.  Equivalently, it is the explicit
root-triflection word
\[
 h=\varepsilon^{-9}
       T_{\beta_0}(T_{\beta_0}T_{\beta_1}\cdots T_{\beta_9})
       \varepsilon,
 \qquad \beta_j=V\alpha_j.
\]
Indeed $VT_{\alpha_j}V^{-1}=T_{\beta_j}$ and
$\varepsilon^{10}=\zeta$.  This closes both remaining fibres by direct
integral identities.

Finally, suppose $A'x=x$ and $P(A')$ is $P\Gamma$-conjugate to $P(A)$.
Lift the projective conjugator to $g\in\Gamma$, so
$gA'g^{-1}=\zeta^jA$. Taking $h=g_A^{-j}$ gives
$h(\zeta^jA)h^{-1}=A$.  Then $hg$ carries $x$ to a vector fixed by $A$.
Because $hg$ is an integral unitary isometry, it preserves the norm and the
mirror-free open locus; conjugation by $hg$ identifies
$C_\Gamma(A')$ with $C_\Gamma(A)$ and identifies their stabilizers of the
corresponding vectors.  This is an isomorphism of strict fixed-vector
groupoids, valid for every $x$, primitive or not, and it preserves the
mass denominators.
\end{proof}

\subsection{Order-ten and order-twelve representatives}
\begin{corollary}[Exhaustion of the order-ten classes]\label{cor:order10-exhaustion}
The four explicit representatives $A_1,\ldots,A_4$ of
\rot{thm:order10-golden-centralizer} exhaust the strict order-ten classes
having a mirror-free positive fixed point.
\end{corollary}
\begin{proof}
\Cref{thm:open-class-list} gives exactly four order-ten angle types, each carrying one
projective class, and \rot{thm:order10-golden-centralizer} supplies four
pairwise nonconjugate strict representatives of those types.  Since
\cref{lem:fullerene-scalar-fibre-collapse} identifies strict with projective
classes, these four representatives fill the four projective classes.
\end{proof}

\begin{lemma}[Geometric identification of the order-twelve representatives]
\label{lem:order12-geometric-classes}
The operators $Q_1,Q_5$ are not unitary conjugate and exhaust the strict
order-twelve classes with a mirror-free positive fixed vector.  Their
primitive anti-M\"obius angle types are respectively $\{\pm1\}$ and
$\{\pm5\}$ in units of $\pi/6$.
\end{lemma}
\begin{proof}
By \cref{thm:equal-weight-spectra}, the order-twelve angle types for
\(N=12\) are \((6,1)\) and \((6,5)\), the primitive angles \(\pm1\) and
\(\pm5\) in units of \(\pi/6\), each with the single partition \([12]\),
and by \cref{lem:type-uniqueness} each carries one projective class.
\Cref{lem:dodecagon-model} realizes the type of angle \(t\) by
\((Z_{12},\sigma_t)\) with class \(\bar Q_t\) for \(t=1,5\).  The two
types are distinct, so \(Q_1\) and \(Q_5\) are not conjugate, even
projectively, and \cref{lem:fullerene-scalar-fibre-collapse} identifies
strict with projective classes; hence \(Q_1,Q_5\) exhaust the strict
order-twelve classes with a mirror-free positive fixed vector.  The same
conclusion is visible along the fixed loci: by \rot{lem:order12-no-walls}
every positive line of either fixed lattice is mirror-free, so each fixed
locus is a connected hyperbolic line on which the angle type is constant
by \cref{lem:type-constant}.
\end{proof}

\begin{proposition}[Integral representatives of the eleven classes]
\label{prop:mirrorfree-spectrum}
The eleven compatible projective classes of \cref{thm:open-class-list},
grouped by the exact orders \(2,4,6,10,12\) of
\cref{thm:equal-weight-spectra}, have the strict integral representatives
\[
\begin{array}{c|c}
\text{order}&\text{representatives}\\
\hline
2&Q_0\ \text{(split)},\ Q_6\ \text{(antipodal)}\\
4&Q_3\\
6&Q_2,Q_4\\
10&A_1,A_2,A_3,A_4\\
12&Q_1,Q_5,
\end{array}
\]
each class having exactly one strict scalar fibre.  The four order-ten
operators lie outside the root-monomial normalizer, and the mirror-free
mass \(M_{10}^{\mathrm{open}}(m)\) of their four classes is the arithmetic
formula displayed in \cref{thm:exact-all-index-grand-formula}.
\end{proposition}

\begin{proof}
The projective list is \cref{thm:open-class-list},
\cref{lem:period-bridge}(d) matches exact orders, and
\cref{lem:fullerene-scalar-fibre-collapse} identifies the strict scalar
fibres.  \Cref{lem:dodecagon-model} places each root-monomial \(Q_t\) in
the class of its angle type, the split type with \cref{prop:real-top},
the \(Q_4\) type with the three partitions of \cref{lem:q4-merge}, and
the order-twelve types with \cref{lem:order12-geometric-classes}.  The
fixed lattices \(\Lambda_3^{\rm rot},\Lambda_2^{\rm rot},\Lambda_4^{\rm rot}\)
and the strict centralizers of
\(Q_3,Q_2,Q_4\) are established in \rot{thm:siegel-q234-ambient},
\rot{thm:lower-rank-centralizer-bridge} and
\cref{rot-sec:quaternionic-centralizer,rot-sec:analytical-rotary-chamber}
of~\cite{fullereneRotaryPaper}, and those of \(Q_1,Q_5\) in
\rot{prop:order12}.  The order-ten representatives, their integral
invariants, strict centralizers and pairwise nonconjugacy are
\rot{thm:order10-golden-centralizer}; they fill the four order-ten
classes by \cref{cor:order10-exhaustion}, and
\rot{thm:order10-analytical-walls} proves their wall list and the
order-ten mass of \cref{thm:exact-all-index-grand-formula}.
\end{proof}

The strict representatives now exhaust the collision-free witnesses,
with the same stabilizer normalization as the geometric objects. Their
seven grouped masses can therefore be added directly; the next part
substitutes their arithmetic evaluations into Burnside's identity.

\part{The exact count}
\section{Exact counting formula}\label{sec:exact-grand-formula}

The preceding analytical reductions evaluate each sector.
Their sum gives the ordinary fullerene count.

\subsection{Arithmetic notation}
All arithmetic functions and zero conventions are defined once in
\cref{sec:coefficient-functions}. For the sector presentation put
\begin{align*}
 R_3^{(4)}(m)
 &=\frac{3T(m/2)-S(m/2)}{24}-\frac{S(m/4)}{18},\\
 I_3^{(2)}(m)
 &=\frac{11}{144}d_{4,0}(m)+\frac1{16}d_{4,2}(m)
   +\frac1{36}d_{12,0}(m).
\end{align*}

\subsection{Sector formulas and assembly}
In the formulas below, $R_0,\ldots,R_7$ are the finite all-glue sums
of \inv{cor:analytical-lower-split-return}, $R_8$ is the binary formula
of \inv{thm:rank-eight-binary-closed-forms}, and $R_9$ is the unary
formula of \inv{thm:analytical-unary-endpoint}.
The set $\mathcal N_{\rm ap}$ and its ten weights are the antipodal
normal list and state table of \inv{thm:analytical-antipodal-weights}.
All $\mathfrak m_S$ use full future-cone orthogonal stabilizers and
include imprimitive vectors; $M_j^{\rm amb}$ are the ambient masses
of \rot{thm:siegel-q234-ambient}. The functions $T_9,T_{27}$ are the
optimal-embedding class-number mass of \rot{thm:q4-t9-class-number}
and its transfer in \rot{thm:q4-analytical-ternary-transfer}.

\begin{theorem}[Exact counting formula at every index]
\label{thm:exact-all-index-grand-formula}
For every integer $m\geq1$,
\[
 \iso(2m)=\frac12\left(
 \sum_{(c,f)\text{ in the oriented table}}c f(m)
 +\sum_{(c,f)\text{ in the achiral tables}}c f(m)\right),
\]
where the 45 oriented and 181 achiral entries are printed in
\cref{app:complete-coefficient-formula}, using exactly the arithmetic
functions and zero conventions of \cref{sec:coefficient-functions}.
The seven strict antiunitary-sector masses are
\begin{align*}
 M_{\rm split}^{\rm open}(m)
 &=\sum_{r=0}^{9}R_r(m),\\
 M_{\rm antipodal}^{\rm open}(m)
 &=\sum_{N\in\mathcal N_{\rm ap}}
       \frac{\mathcal M_N}{|W(N)|}\,\mathfrak m_{S_N}(m),\\
 M_2^{\rm open}(m)
 &=M_2^{\rm amb}(m)-\frac19\tau(m/6)-\frac1{18}d_{6,4}(m),\\
 M_3^{\rm open}(m)
 &=M_3^{\rm amb}(m)-R_3^{(4)}(m)+I_3^{(2)}(m),\\
 M_4^{\rm open}(m)
 &=M_4^{\rm amb}(m)
   +\frac{N_3^+(m)+N_6^+(m)+\mathbf1_{3\mid m}N_3^+(m/3)}6\\
 &\quad-\frac1{18}d_{6,1}(m)-\frac19d_{6,3}(m)
   -T_9(m)-T_{27}(m)\\
 &\quad
   -\frac19I_1(m)-\frac1{18}I_3(m)+\frac1{18}I_9(m),\\
 M_{10}^{\rm open}(m)
 &=\frac25\left(
      \tau(m/10)+d_{\rm same}(m/5)-I_5(m)-I_{15}(m)
    \right),\\
 M_{12}^{\rm open}(m)&=\frac13\tau(m/12),
\end{align*}
Their sum is
\begin{equation}\label{eq:exact-grand-coefficient}
 \boxed{\begin{aligned}
 a_m={}&M_{\rm split}^{\rm open}(m)+M_{\rm antipodal}^{\rm open}(m)
          +M_2^{\rm open}(m)+M_3^{\rm open}(m)\\
       &+M_4^{\rm open}(m)+M_{10}^{\rm open}(m)+M_{12}^{\rm open}(m)
 \end{aligned}}
\end{equation}
and, with $\iso$ as in \eqref{eq:burnside-count},
\begin{equation}\label{eq:exact-all-index-isomer-count}
 \iso(2m)=\frac{\BB_m+a_m}{2}.
\end{equation}
Equivalently, with $q=e^{2\pi iz}$ for $z\in\mathbb H$ and the generating
series
\[
 \mathcal A(q)=\sum_{m\geq1}a_mq^m,\qquad
 \mathcal B(q)=\sum_{m\geq1}\BB_mq^m,\qquad
 \mathcal I(q)=\sum_{m\geq1}\iso(2m)q^m,
\]
and $\Theta_j^{\rm open}(q)=\sum_{m\geq1}M_j^{\rm open}(m)q^m$, all of
which converge locally uniformly on $\mathbb H$,
\[
 \mathcal A=\Theta_{\rm split}^{\rm open}
 +\Theta_{\rm antipodal}^{\rm open}
 +\Theta_{2}^{\rm open}+\Theta_{3}^{\rm open}
 +\Theta_{4}^{\rm open}+\Theta_{10}^{\rm open}
 +\Theta_{12}^{\rm open},
 \qquad
 \mathcal I=\frac12(\mathcal B+\mathcal A).
\]
Every coefficient on the right is an exact finite arithmetic quantity.
\end{theorem}

The indicator $\mathbf1_{3\mid m}$ in front of $N_3^+(m/3)$ is
redundant under the zero-extension convention and is retained to match
\rot{thm:rotary-rows}. For direct evaluation,
\cref{app:complete-coefficient-formula} expands $a_m$ and $\BB_m$
entirely into functions defined in this paper.

\begin{proof}
Project any strict mirror-free antiunitary witness to the period ball.  By
\cref{lem:period-bridge}(c) and \cref{thm:open-class-list}, its projective
class is one of the eleven classes represented in
\cref{prop:mirrorfree-spectrum}.
\Cref{lem:fullerene-scalar-fibre-collapse} then identifies every strict lift
over each compatible projective type with its displayed representative,
including on imprimitive shells.  These strict classes are disjoint and
exhaust every mirror-free antiunitary witness.  The
centralizer-weighted witness formula therefore partitions $a_m$ into exactly
the seven rows in
\eqref{eq:exact-grand-coefficient}; the four order-ten and two order-twelve
classes have already been summed within their displayed rows.

\Inv{thm:involutive-rows}, with
\inv{thm:involutive-all-index-collisions} identifying the root restrictions,
gives the first two formulas.  \Rot{thm:rotary-rows} gives the remaining
five: in particular \rot{cor:q3-analytical-open-return} supplies the complete
$Q_3$ assembly and \rot{cor:q4-analytical-open-return} the $Q_4$ assembly.
Their ambient masses are \rot{thm:siegel-q234-ambient}; the order-ten
and order-twelve returns are \rot{thm:order10-analytical-walls} and
\rot{thm:analytical-order-twelve}.  Every return uses the same index and
the full strict vector stabilizer, including its ineffective kernel.
The lower-split terms here are orbit masses. At rank seven,
\inv{thm:rank-seven-ternary-class-forms} evaluates them in class-number
form; that evaluation is needed for the printed arithmetic expansion,
but not for the orbit-mass identity itself.
Each chamber shell is
finite at fixed positive norm; every residual term is a finite divisor,
class-number, finite-polar-space, or arithmetic-orbit calculation.  Hence
the displayed sum is defined for every $m$.

The oriented count has the parallel five-sector identity in
\uni{thm:oriented-five-sector-reduction}. Its homogeneous complex masses
and zero-period collision sums yield the forty-five coefficients by
\uni{prop:oriented-collision-scalar-assembly};
\cref{prop:analytical-oriented-count} records this assembly.

The groupoid proof of the witness formula makes the sum an integer: strict
stabilizer denominators cancel orbitwise.  Burnside's lemma then gives
\eqref{eq:exact-all-index-isomer-count} and its parity.  Finally
$0\leq a_m\leq\BB_m$, while the oriented series is a finite
Eisenstein/quasimodular combination of polynomially growing coefficients;
this gives local uniform convergence on $\mathbb H$.
\end{proof}

\begin{remark}[Scope of the generating series]
The symbols $\Theta_A^{\rm open}$ in the theorem denote the exact convergent
sector generating series.  Their arithmetic evaluation is a sum of
contributions of different ranks.
\end{remark}

\paragraph{Census identifiability and independent checks.}
On the 181 achiral columns of the formula, evaluation at $m\leq200$
has rank 166 and a 15-dimensional kernel, supported on 19 weight-$3/2$
split columns. The census determines the coefficients modulo this
kernel; the proof above fixes them through all-index stratum identities.
The exact kernel and later pivot indices are given in
\cref{prop:counting-space-coordinates} and the following discussion.
The first new direction occurs at $m=208$, where
$a_{208}=983001$ and $\iso(416)=189565760391$.

For direct use of the seven rows, the closed forms of $M_2^{\rm amb}$,
$M_3^{\rm amb}$, $M_4^{\rm amb}$, $T_9$ and $T_{27}$ are printed at the end
of \cref{sec:coefficient-functions}, the first subsection of
\cref{app:complete-coefficient-formula}, and the thirteen-term closed form
of $M_{\rm antipodal}^{\rm open}$ is printed in
\cref{prop:counting-space-coordinates}.

\subsection{The oriented contribution}
\begin{proposition}[Analytical oriented count]
\label{prop:analytical-oriented-count}
For every positive integer $m$,
\[
 \BB_m=V_m+O_2(m)+O_3(m)+O_5(m)+O_6(m).
\]
Here $V_m$ is the open period mass weighted by reciprocal unitary
stabilizers, and $O_q(m)$ is the sum of the open strict-centralizer
masses for unitary rotations of order $q$. There is one contributing
class for each of orders $2,3,6$ and two for order five.
The five masses evaluate to exactly the forty-five terms in
\cref{app:oriented-coefficients}.
\end{proposition}
\begin{proof}
For a finite rotational stabilizer $G$, weighting its conjugacy classes
of marked rotations gives
\[
 \sum_{[g]\subset G}\frac1{|C_G(g)|}=1.
\]
The period correspondence therefore counts every oriented fullerene once.
The strict class list and quotient groups are
\uni{thm:oriented-five-sector-reduction}. The identity mass is
\uni{cor:unitary-weighted-open-count}; the rotation masses and their
all-index scalar assembly are
\uni{prop:oriented-collision-scalar-assembly}. The latter reduces the
individual-component formulas by divisor recurrences and character
projections.
\end{proof}

\begin{remark}[Coefficients of the oriented formula]
\label{rem:oriented-misprint}
The oriented formula of Engel and Smillie
\cite[Theorem~3.26]{engel2025} is used here with coefficients $1/18$ and
$1/48$ for $\sigma_1(m/9)$ and $\sigma_1(m/12)$, replacing the printed
$1/64$ and $2/27$. Their archived Sage implementation already uses
the correct values. The derivation and coefficient comparison are in
\uni{prop:oriented-collision-scalar-assembly};
\cref{sec:data-audit} records the source check.
\end{remark}

\part{Asymptotic counting and generic symmetry}
\section{Achiral growth and typical symmetry}
\label{sec:achiral-asymptotics}\label{sec:generic-introduction}

An achiral fullerene may have several orientation-reversing witnesses, so a
weighted fixed-locus count does not automatically describe the distribution
of isomorphism classes.  We separate the exceptional classes with a
nontrivial orientation-preserving automorphism and prove that they have lower
order.  This turns the two leading involutive witness masses into literal
asymptotic class proportions.

\subsection{The achiral main term}

The two rank-ten fixed lattices give the leading term. Collision
corrections have smaller rank and hence smaller growth; averaging the
resulting arithmetic components gives the cumulative count.

\begin{theorem}[Main term and error for the achiral count]
\label{thm:achiral-asymptotics}
Write $m=2^am'$ with $m'$ odd. The combined split and antipodal
weight-five contribution is
\[
 W_5(m)=c(a)\,\sigma_4(m'),\qquad
 c(0)=\frac{697}{622080},\qquad
 c(a)=\alpha\cdot16^{\,a-1}-\beta\quad(a\geq1),
\]
with $\alpha=5617/291600$ and $\beta=21689/9331200$.  Then
\[
 a_m=W_5(m)+O\bigl(m^{7/2}\bigr),
\]
with an implied constant independent of $m$.  Consequently
\[
 \sum_{m\leq x}a_m=\zeta(5)\,\frac{21689}{99532800}\,x^5+O\bigl(x^{9/2}\bigr),
\]
\[
 \liminf\frac{a_m}{m^4}=\frac{10537}{9953280},
 \qquad
 \limsup\frac{a_m}{m^4}=\frac{5617}{447897600}\pi^4,
\]
and $\operatorname{iso}(2m)=\tfrac12\BB_m+\tfrac12W_5(m)+O(m^{7/2})$.
\end{theorem}

\begin{proof}
The weight-five block of the printed table in
\cref{app:complete-coefficient-formula} is
\[
 \frac{697}{622080}\sigma_4(m)-\frac{41}{19440}\sigma_4(m/2)
 +\frac{697}{19440}\sigma_4(m/4);
\]
these are the two rank-ten ambient terms, the split and antipodal masses
of \inv{prop:analytical-rank10-masses}, and every other term of the
table comes from a lower-rank fixed or collision stratum.  Writing $m=2^am'$
and using $\sigma_4(2^am')=\sigma_4(2^a)\sigma_4(m')$ with
$\sigma_4(2^a)=(16^{a+1}-1)/15$, the block equals
$c(a)\sigma_4(m')=W_5(m)$ with the displayed $c(0)$, $\alpha$ and
$\beta$.  \Cref{prop:effective-achiral-remainder} bounds every other
term of the table by an explicit constant times $m^{7/2}$ or $m^3$, so
$|a_m-W_5(m)|\leq C_{7/2}m^{7/2}+C_3m^3$ with the constants of
\cref{app:effective-error}; this is the stated remainder.

For the cumulative constant, sum the three terms of $W_5$ with
\[
 \sum_{m\leq x}\sigma_4(m/j)=\frac{\zeta(5)}{5j^5}x^5+O(x^4),
\]
which gives
$\zeta(5)x^5\bigl(\tfrac{697}{622080}-\tfrac{41}{19440\cdot32}
+\tfrac{697}{19440\cdot1024}\bigr)/5=\zeta(5)\tfrac{21689}{99532800}x^5$,
while the remainder contributes $O(x^{9/2})$.

For the extremes, $W_5(m)/m^4=\bigl(c(a)/16^a\bigr)\,\sigma_4(m')/m'^4$
and $a_m/m^4=W_5(m)/m^4+O(m^{-1/2})$.  For odd $m'$ the second factor
satisfies $1\leq\sigma_4(m')/m'^4<\prod_{p\ \mathrm{odd}}(1-p^{-4})^{-1}
=\tfrac{15}{16}\zeta(4)$, with lower limit $1$ along odd primes and upper
limit $\tfrac{15}{16}\zeta(4)$ along products of the first odd primes
with growing exponents.  The first factor depends only on the dyadic
valuation: $c(0)=\tfrac{697}{622080}$, and for $a\geq1$,
$c(a)/16^a=\alpha/16-\beta/16^a$ increases strictly with $a$ from
$(\alpha-\beta)/16=\tfrac{10537}{9953280}$ at $a=1$ to the limit
$\alpha/16$; since $(\alpha-\beta)/16<c(0)<\alpha/16$, the minimum over
all $a$ is attained at $a=1$ and the supremum is $\alpha/16$.  Hence
$\liminf a_m/m^4=\tfrac{10537}{9953280}$, approached along $m=2p$ with
$p$ an odd prime, and
$\limsup a_m/m^4=\tfrac{\alpha}{16}\cdot\tfrac{15}{16}\zeta(4)
=\tfrac{5617}{447897600}\pi^4$.
Finally $\iso(2m)=\tfrac12(\BB_m+a_m)$ by \cref{lem:burnside}, which
gives the last identity.
\end{proof}

\begin{corollary}[Sign and size of the leading error]
\label{cor:achiral-leading-error-sign}
Let \(S_{9/2}(m)\) denote the six weight-\(9/2\) terms in
Appendix~\ref{app:complete-coefficient-formula}, with the stated Cohen
normalization.  For every positive integer \(m\),
\[
 S_{9/2}(m)\leq
 -\left(\frac{41}{5670}+\frac{\sqrt3}{420}\right)m^{7/2}.
\]
Moreover \(a_m-W_5(m)=S_{9/2}(m)+O(m^3)\), with a uniform implied
constant.  Consequently \(a_m<W_5(m)\) for every sufficiently large
\(m\), and \(W_5(m)-a_m\) is bounded above and below by positive
constant multiples of \(m^{7/2}\) in that range.
\end{corollary}
\begin{proof}
Cohen's definitions \cite[Definitions~2.1--2.2, pp.~272--273]{cohen1975}
give, for \(N>0\) with \(N\equiv0,1\pmod4\),
\[
 H_4(N)=\frac{3}{4\pi^4}
 \sum_{\substack{d^2\mid N\\N/d^2\equiv0,1\pmod4}}
 (N/d^2)^{7/2}L(4,\chi_{N/d^2}),
\]
where \(\chi_M=(M/\,\cdot\,)\) is the possibly imprimitive quadratic
character.  Every Euler factor is positive, and
\[
 L(4,\chi_M)\geq\prod_p(1+p^{-4})^{-1}
 =\frac{\zeta(8)}{\zeta(4)}=\frac{\pi^4}{105}.
\]
Keeping the \(d=1\) summand proves \(H_4(N)\geq N^{7/2}/140\).
Off this support the coefficient is zero.  All six coefficients in
\(S_{9/2}\) are negative, and \(4m,12m\) always lie in the support.
The two corresponding terms alone therefore give
\[
 S_{9/2}(m)\leq
 -\frac{41(4m)^{7/2}}{5184\cdot140}
 -\frac{(12m)^{7/2}}{10368\cdot140}
 =-\left(\frac{41}{5670}+\frac{\sqrt3}{420}\right)m^{7/2}.
\]
Every remaining component has coefficient growth \(O(m^3)\): the
weight-four divisor sums have this bound, and the lower Cohen,
divisor, class-number and indicator terms have smaller exponents,
including any logarithmic factors.  The upper bound follows from
Theorem~\ref{thm:achiral-asymptotics}; the displayed strict negative
leading bound gives the lower bound and the eventual inequality.
\end{proof}

\begin{figure}[htbp]
\centering
\includegraphics[width=\textwidth]{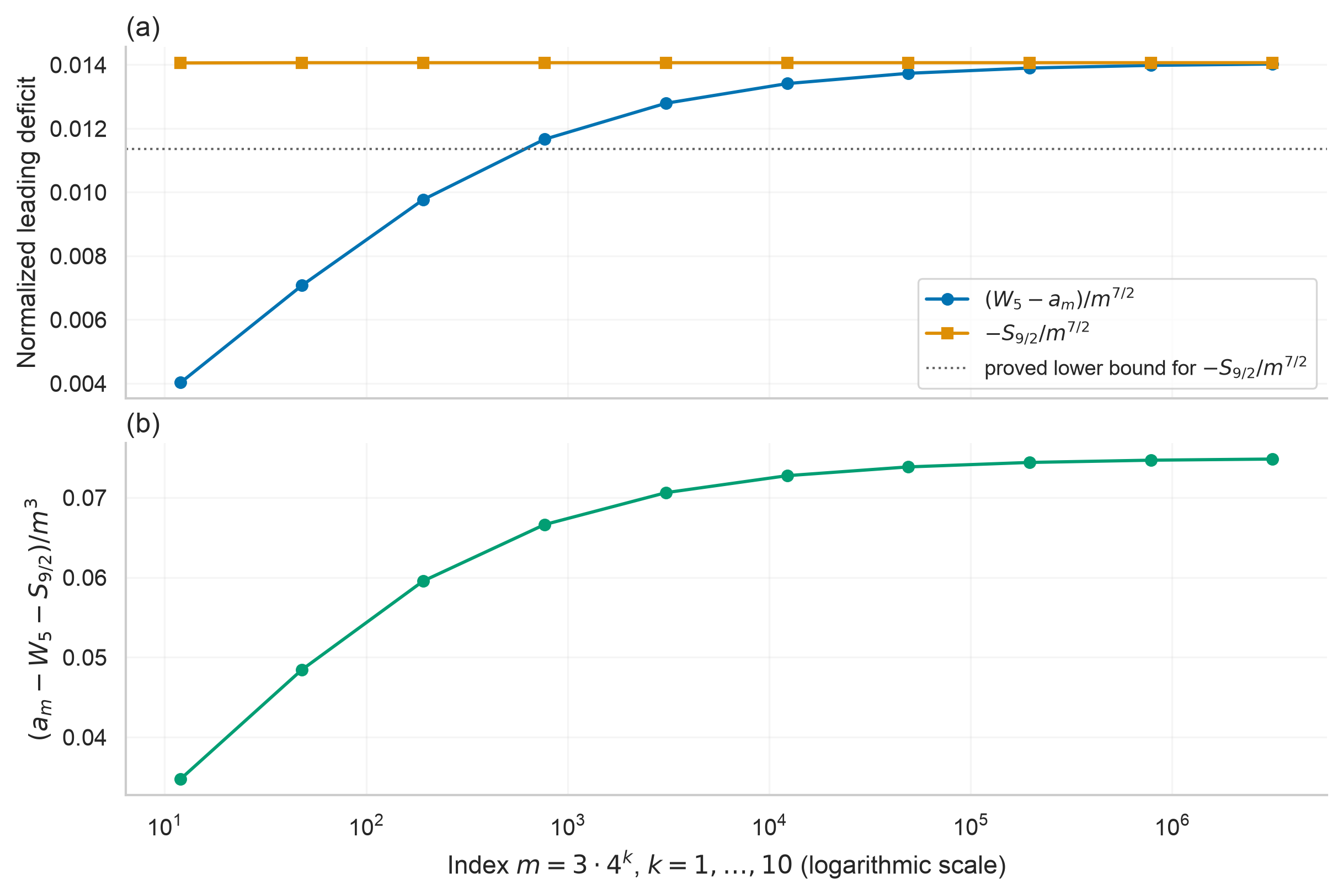}
\caption{The leading collision correction, on the fixed-square-class
family $m=3\cdot4^k$, $1\leq k\leq10$.
(a) The normalized deficit and the negative weight-$9/2$ Cohen block.
The dotted line is the proved lower bound for the Cohen block alone in
\cref{cor:achiral-leading-error-sign}.
(b) The remaining signed contribution, divided by $m^3$.
Every point is evaluated from the printed arithmetic formula with exact
rational subtraction before plotting. Connecting segments guide the eye;
no limiting value or explicit band for the $O(m^3)$ remainder is inferred.}
\label{fig:counting-cancellation}
\end{figure}

Figure~\ref{fig:counting-cancellation} separates the negative Cohen
block from the remaining signed correction along $m=3\cdot4^k$.

\subsection{Rarity among all fullerenes}
\begin{corollary}[Sharp asymptotic rarity of achirality]
\label{cor:asymptotic-chirality}
Set
\[
 C_{\mathrm{or}}=\frac{809}{2^{15}3^{13}5^2}.
\]
Then the oriented and ordinary counts satisfy
\[
 \BB_m=C_{\mathrm{or}}\sigma_9(m)+O(m^8),
 \qquad
 \iso(2m)=\frac{C_{\mathrm{or}}}{2}\sigma_9(m)+O(m^8).
\]
For every $m$ with $\iso(2m)>0$, let
\[
 p_m=\frac{a_m}{\iso(2m)}
\]
be the proportion of ordinary $2m$-vertex fullerene classes that are
achiral.  With $m=2^am'$ and $m'$ odd, one has the sharper pointwise formula
\[
 p_m=
 \frac{2c(a)\sigma_4(m')}{C_{\mathrm{or}}\sigma_9(m)}
 +O(m^{-11/2})
 =\Theta(m^{-5}).
\]
In fact,
\[
 \liminf_{m\to\infty}m^5p_m
 =\frac{64(\alpha-\beta)}{513C_{\mathrm{or}}},
 \qquad
 \limsup_{m\to\infty}m^5p_m
 =\frac{15\alpha}{128C_{\mathrm{or}}}
   \frac{\zeta(4)}{\zeta(9)}.
\]
The corresponding cumulative achiral proportion has a genuine asymptotic:
\[
 \frac{\sum_{m\leq x}a_m}{\sum_{m\leq x}\iso(2m)}
 =\frac{40\cdot21689\cdot3^8}{809}
   \frac{\zeta(5)}{\zeta(10)}x^{-5}
  +O(x^{-11/2}).
\]
Thus the sizewise and cumulative proportions of chiral ordinary fullerenes
tend to one, and the exponent five is sharp.
\end{corollary}

The growth exponent nine of the total count is already present in
Sah's counting estimates \cite[Theorem~3.4]{sah1994generalizedleapfrog} and in the oriented
formula \cite{engel2025}.  The assertion here is the order-five gap
between ordinary and achiral counts, together with its sizewise
envelopes and cumulative constant.  Earlier chirality results for the
icosahedral family \cite[Theorem~2.3]{sah1994generalizedleapfrog} concern that restricted
family; the comparison with other symmetry enumerations is given in
\cref{app:symmetry-comparisons}.

\begin{proof}
The analytically derived oriented table of
\uni{prop:oriented-collision-scalar-assembly}, reproduced in
\cref{app:oriented-coefficients}, has the unique weight-ten summand
$C_{\mathrm{or}}E_{10}(q)$; every other summand has weight at most nine.
Its Fourier normalization is
\[
 [q^n]E_{k,\chi,\psi}(q)
 =\sum_{d\mid n}\chi(n/d)\psi(d)d^{k-1},
\]
so the absolute value is at most $\sigma_{k-1}(n)$.  A fixed dilation merely
replaces $n$ by $n/t$ when $t\mid n$.  For $k\geq3$,
$\sigma_{k-1}(n)\leq\zeta(k-1)n^{k-1}$; the weight-two and weight-one terms
are bounded by $n^2$ and $n$, respectively.  The finite sum of all terms
below weight ten is therefore $O(m^8)$. Since
$[q^m]E_{10}=\sigma_9(m)$, this proves the formula for
$\BB_m$.  The formula for $\iso(2m)$ follows from Burnside and
$a_m=O(m^4)$, which is immediate from
\cref{thm:achiral-asymptotics}.

As $\sigma_9(m)\geq m^9$, the denominator
$\BB_m+a_m=C_{\mathrm{or}}\sigma_9(m)+O(m^8)$ can be inverted uniformly for
large $m$.  Substituting
$a_m=c(a)\sigma_4(m')+O(m^{7/2})$ gives the displayed pointwise formula.
For its sharp envelopes, put
\[
 h_a=\frac{2^{5a}c(a)}{\sigma_9(2^a)},
 \qquad
 g(n)=\frac{n^5\sigma_4(n)}{\sigma_9(n)}\quad(n\ \text{odd}).
\]
Then $m^5p_m=(2/C_{\mathrm{or}})h_ag(m')+O(m^{-1/2})$.  Direct rational
comparison gives $h_1<h_0<\lim_{a\to\infty}h_a$.  For $a\geq1$, substitution
gives
\[
 h_a=\frac{511(\alpha/16-\beta2^{-4a})}{512-2^{-9a}},
\]
and cross-multiplication shows that this sequence is increasing.  Its two
extreme values are
\[
 h_1=\frac{32(\alpha-\beta)}{513},
 \qquad
 \lim_{a\to\infty}h_a=\frac{511\alpha}{8192}.
\]
For odd prime powers,
\[
 g(p^e)=
 \frac{1-p^{-4(e+1)}}{1-p^{-4}}
 \frac{1-p^{-9}}{1-p^{-9(e+1)}}.
\]
Hence $\liminf g(n)=1$, approached along odd primes, and
\[
 \limsup g(n)
 =\prod_{p\ \mathrm{odd}}\frac{1-p^{-9}}{1-p^{-4}}
 =\frac{\zeta(4)}{\zeta(9)}
  \frac{1-2^{-4}}{1-2^{-9}}.
\]
This proves both pointwise limits; in the upper limit the factors $511$
cancel because $(1-2^{-4})/(1-2^{-9})=480/511$.

Finally, the standard divisor-sum estimate gives
\[
 \sum_{m\leq x}\sigma_9(m)
 =\frac{\zeta(10)}{10}x^{10}+O(x^9).
\]
Combine it with the summatory formula in
\cref{thm:achiral-asymptotics} and Burnside, then divide the two asymptotic
series.  Since $99532800=2^{14}3^5 5^2$, the leading constant simplifies to
the one displayed above.
\end{proof}

\begin{figure}[htbp]
\centering
\includegraphics[width=\textwidth]{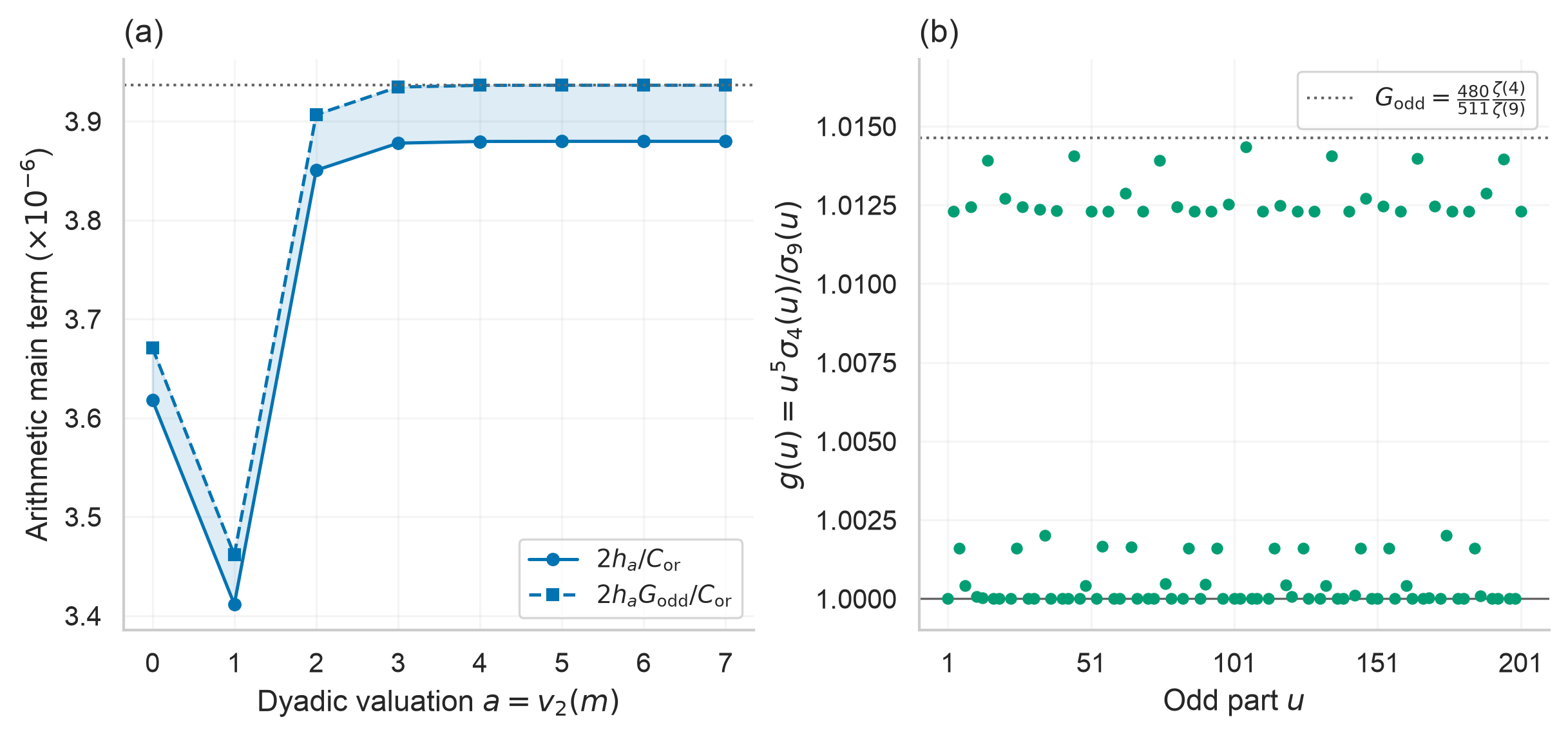}
\caption{The two arithmetic factors in the main term
$F(2^au)=2h_ag(u)/C_{\rm or}$ of \cref{cor:asymptotic-chirality}.
(a) Its exact lower and upper envelopes at each integer dyadic valuation,
scaled by $10^{-6}$; shading between joining segments guides the eye at
the discrete valuations. The dotted line is the limiting upper envelope
as $a\to\infty$.
(b) The exact odd factor at $1\leq u\leq201$, with its upper envelope
$G_{\rm odd}=480\zeta(4)/(511\zeta(9))$.
The bounds enclose the arithmetic main term; the finite counts satisfy
$m^5p_m=F(m)+O(m^{-1/2})$. This separation explains the arithmetic
fluctuations and the distinct lower and upper limits.}
\label{fig:counting-arithmetic-envelopes}
\end{figure}

Figure~\ref{fig:counting-arithmetic-envelopes} shows how the dyadic
envelopes and the odd factor $g(m')$ account for the fluctuations in
the arithmetic main term of $m^5p_m$.

\subsection{Bounds for exceptional symmetry}
\begin{lemma}[Conjugation on the decorated period module]
\label{lem:decorated-period-conjugation}
Let $(H,\mu)$ be a triangulable bistratum of \cite{engel2025} with
$H\ne1$, let
$(\Lambda_\mu,\Gamma_{H,\mu})$ be its Eisenstein period module and
monodromy group, and let $Y^{\mathrm{ex},\circ}_{H,\mu,n}$ be the relative
open locus of orbitriangulations of index $n$ whose $H$-cover has full
orientation-preserving automorphism group exactly $H$.  After choosing a
conjugation-compatible marking, complex conjugation on the decorated period
local system is represented by an integral antiunitary $\chi_\mu$ satisfying
\[
 \chi_\mu^2=1,
 \qquad
 \chi_\mu\Gamma_{H,\mu}\chi_\mu^{-1}=\Gamma_{H,\mu}.
\]
The odd arrows coming from decorated quotient geometry lie in the single
compatible coset $\Gamma_{H,\mu}\chi_\mu$.
\end{lemma}

\begin{proof}
Proposition~2.19 of \cite{engel2025} constructs
\[
 \Lambda_{\mu^u}\subset\Lambda_\mu
 \subset\frac{1}{2(1+\zeta_6)}\Lambda_{\mu^u}
\]
and the finite-index group $\Gamma_{H,\mu}$.  The period conditions in its
proof are stable under conjugation: besides $\mathcal E=\mathbb Z[\zeta_6]$
they use only
\[
 \tfrac12\mathcal E,
 \qquad
 (1+\zeta_6)^{-1}\mathcal E,
\]
and conjugation preserves the first ideal and sends the second to its unit
multiple $\zeta_6(1+\zeta_6)^{-1}\mathcal E$.  It also preserves the
ramification orders and the edge/face decorations, so it acts on the
decorated intermediate moduli cover used in the proof of that proposition.
Choose its labelled quotient points on $\mathbb{RP}^1$ and the cyclic
period cover with real defining coefficients. Coordinatewise conjugation
fixes this decorated base object and squares to one on its integral
local system, supplying a conjugation-compatible base marking.
Naturality of Gauss--Manin transport on this cover supplies an integral
antiunitary $\chi_\mu$, with $\chi_\mu^2=1$ after choosing a fixed base
marking, and proves
$\chi_\mu\Gamma_{H,\mu}\chi_\mu^{-1}=\Gamma_{H,\mu}$.  The ideal calculation
shows that this action stays inside $\Lambda_\mu$.
\end{proof}

\begin{lemma}[A finite class set before imposing the area]
\label{lem:decorated-finite-real-classes}
For a fixed triangulable bistratum $(H,\mu)$ with $H\ne1$, the strict
antiunitary involutions supplied by real decorated quotient geometry on
the exact-$H$ locus meet finitely many $\Gamma_{H,\mu}$-conjugacy classes.
Denote their set by $\mathcal A^{\mathrm{ex},\circ}_{H,\mu}$. It is
independent of the index $n$.
\end{lemma}
\begin{proof}
Put $V=\Lambda_\mu\otimes_{\mathcal E}\mathbb Q(\zeta_6)$ and
$\Gamma=\Gamma_{H,\mu}$. First, $\Gamma$ is arithmetic in the
$\mathbb Q$-algebraic unitary group $U(V)$. In rank one this follows
from finiteness of the integral unitary group. In higher rank,
Thurston's theorem gives finite volume to
$P\Gamma_{\mu^u}\backslash\mathbb B$ \cite[Thm.~0.2]{thurston1998}.
The containing integral group $PU(\Lambda_{\mu^u})$ is discrete, so
the index of $P\Gamma_{\mu^u}$ in it is finite: the covolume is
multiplied by that index, and the containing group's covolume is
positive. Proposition~2.19 of \cite{engel2025}, the finite scalar
kernel, and the commensurable lattices
$\Lambda_{\mu^u}\subset\Lambda_\mu$ therefore make $\Gamma$
commensurable with $U(\Lambda_\mu)$.

By \cref{lem:decorated-period-conjugation}, the integral antiunitary
$\chi_\mu$ normalizes $\Gamma$ and satisfies $\chi_\mu^2=1$.
Conjugation by $\chi_\mu$ defines a $\mathbb Q$-algebraic
automorphism $s$ of $U(V)$ of order dividing two: although
$\chi_\mu$ is antilinear over $\mathbb Q(\zeta_6)$, it is linear
over $\mathbb Q$. Thus $C_2$ acts on the arithmetic group $\Gamma$
through this algebraic automorphism. Borel--Serre's finiteness theorem
\cite[Proposition~3.8(ii), p.~136]{borelserre1964} applies to
precisely this situation and gives a finite set $H^1(C_2,\Gamma)$.
Their definition uses commensurability with integral points; no
maximality or congruence-subgroup assumption is required.

An odd element $A=\gamma\chi_\mu$ is an involution exactly when
\[
 \gamma s(\gamma)=1.
\]
Conjugating $A$ by $h\in\Gamma$ replaces its cocycle by
$h\gamma s(h)^{-1}$. Hence $H^1(C_2,\Gamma)$ is exactly the set of
$\Gamma$-conjugacy classes of involutions in $\Gamma\chi_\mu$.
The classes realized by decorated geometry on the exact-$H$ locus
form a subset. Taking realizations before specifying an index gives
the claimed finite, index-independent set.
\end{proof}

\begin{lemma}[Decorated real-period injection]
\label{lem:decorated-real-period-injection}
Let $(H,\mu)$ be a triangulable bistratum with $H\ne1$, and use the
finite set $\mathcal A^{\mathrm{ex},\circ}_{H,\mu}$ of
\cref{lem:decorated-finite-real-classes}. Every achiral triangulation $T$
with $\operatorname{Aut}^{+}(T)=H$, quotient profile $\mu$, and $2m$
triangles determines a pair $(A,v)$ with
\[
 [A]\in\mathcal A^{\mathrm{ex},\circ}_{H,\mu},\qquad
 Av=v,\qquad \frac{(v,v)}2=\frac{m}{|H|}.
\]
Exact-cover reconstruction gives an injection
\[
 \{[T]\}\hookrightarrow
 \coprod_{[A]\in\mathcal A^{\mathrm{ex},\circ}_{H,\mu}}
 C_{\Gamma_{H,\mu}}(A)\backslash
 \left\{v\in\Lambda_\mu^{A,+}:\frac{(v,v)}2=\frac{m}{|H|}\right\}.
\]
Every image orbit has trivial strict stabilizer and groupoid weight one.
The vector $v$ may be imprimitive.
\end{lemma}
\begin{proof}
Use the compatible conjugation of
\cref{lem:decorated-period-conjugation}.

The group $H=\operatorname{Aut}^{+}(T)$ is the kernel of the orientation
character and is normal.  A reverser descends to $T/H$; all reversers induce
the same quotient map, and its square is the identity.  More concretely, an
orientation-reversing simplicial isometry has developed local form
\[
 z\longmapsto u\bar z+b,\qquad u\in\mu_6.
\]
Translation disappears on period differences, so $a\mapsto u\bar a$ is an
integral antiunitary preserving the three period ideals above.  Transporting
the marking by the geometric symmetry produces a strict $A$ fixing the
period vector $v$.  The argument applies to every period vector,
primitive or not.
Equivalently, under the vector-orbit bijection of
\cite[Prop.~2.19]{engel2025}, conjugation sends the class of $v$ to the
class of $\chi_\mu v$.  Since the exact quotient orbit is fixed, choose
$\gamma\in\Gamma_{H,\mu}$ with $\gamma\chi_\mu v=v$ and put
$A=\gamma\chi_\mu$.  Normalization of $\Gamma_{H,\mu}$ by $\chi_\mu$
gives $A^2\in\Gamma_{H,\mu}$, and $A^2v=v$.  Proposition~2.24 of
\cite{engel2025} makes the exact-locus stabilizer $\Gamma_{H,\mu,v}$
trivial, so $A^2=1$.

The resulting involution belongs to the index-independent finite set
of \cref{lem:decorated-finite-real-classes}.

Finally, the lifting assertion behind Proposition~2.24 of
\cite{engel2025} can be seen directly on the decorated quotient.
The sphere carrying the $H$-action is the universal orbifold cover of
the sphere with the prescribed ramification orders. An automorphism
preserving those orders lifts to that cover; its lifts differ by deck
transformations in $H$. The quotient metric and its vertex, edge and
face decorations pull back to the metric and triangulation upstairs.
Thus a nonidentity orientation-preserving automorphism of the decorated
orbitriangulation would lift to an orientation-preserving automorphism
of its cover outside $H$. The exact-$H$ condition excludes it. This
proves triviality of the decorated quotient stabilizer, rather than
assuming that arbitrary quotient symmetries lift to arbitrary covers.
In particular, for $H=C_3$, a quotient symmetry interchanging the two
branch poles inverts the monodromy generator. Every lift then normalizes
$H$ and lies outside $H$, so it is excluded on the exact-$H$ locus.
The same universal-cover argument lifts an antiholomorphic quotient
symmetry, and its lifts form one $H$-coset. Its square on quotient
periods lies in the trivial strict stabilizer already identified above.
Consequently the associated odd arrow is unique, has square one, and
has trivial strict stabilizer. Proposition~2.24 identifies the resulting
quotient orbitriangulations with their exact-$H$ covers, both of weight
one; forgetting the exact-cover and open conditions gives the displayed
injection. By
\cref{prop:fullerene-period-correspondence}(d), a quotient with $2n$
triangles has $n=(v,v)/2$.  The $H$-cover multiplies area by $|H|$, so
$n=(v,v)/2=m/|H|$, as claimed.
\end{proof}

\begin{lemma}[Growth of a fixed low-rank Lorentzian mass]
\label{lem:low-rank-mass-growth}
For a fixed Lorentzian lattice $F$ of rank $r\leq6$, choose a fixed positive
integer $c_F$ making its index form an integral Gram form, and write
$N=c_Fn\in\mathbb Z$.
For every $\eta>0$, its full positive-vector mass at norm $n>0$ is
\[
 \ll_{F,\eta}
 \begin{cases}
  1,&r=1,\\
  n^\eta,&r=2,\\
  n^{1/2+\eta},&r=3,\\
  n^{1+\eta},&r=4,\\
  n^{r/2-1},&5\leq r\leq6.
 \end{cases}
\]
The same bounds hold for any fixed arithmetic subgroup commensurable with
$O(F)$, with finite ineffective kernel.
\end{lemma}
\begin{proof}
For $r\geq4$, the per-class Siegel--Kneser formula
\inv{thm:siegel-class-mass} expresses the mass as a fixed constant
times $N^{r/2-1}$ and the product of the counting densities of the
integral quadratic polynomial $Q$ obtained by clearing the value lattice.
For rank four, at a prime not dividing $2\det B\,N$ the unimodular
Gauss sum gives a density $1+O(p^{-2})$. At a prime dividing $N$ but
not $2\det B$, the conductor estimate below, with $r=4$, bounds the
density by $(1-p^{-1})^{-1}$. The same estimate bounds each of the
finitely many primes dividing $2\det B$ by a constant depending only
on $F$. Thus the product is bounded by a convergent product times
\(\prod_{p\mid N}(1-p^{-1})^{-1}\ll_\eta N^\eta\).
The first assertion can also be read directly from
\cite[Thm.~3.3 and Cor.~4.2]{yang1998}: when $p\nmid N$, the
nontrivial conductor-one Gauss term is $O(p^{-2})$, and the higher
conductor terms vanish. Here $B$ is the integral polar matrix of $Q$,
as in the next paragraph.

For $r\geq5$ the density product is uniformly bounded in $N$.
Let $B$ be the nonsingular integral polar matrix of $Q$. For $p\nmid a$
and $h\geq1$, write
\[
 G_h(a)=p^{-hr}\sum_{x\bmod p^h}
              e^{2\pi i aQ(x)/p^h}.
\]
Differencing the square of this sum and summing over $x$ gives
\[
 |G_h(a)|^2\leq p^{-hr}\#\ker(B\bmod p^h)
       \leq |\det B|_p^{-1}p^{-hr}.
\]
The last inequality follows from the Smith normal form of $B$.
Fourier expansion of the congruence $Q(x)\equiv N\pmod{p^k}$,
grouped by the conductor $p^h$ of its additive character, therefore gives
\[
 p^{-k(r-1)}\#\{x\bmod p^k:Q(x)\equiv N\pmod{p^k}\}
 \leq1+|\det B|_p^{-1/2}\sum_{h=1}^k p^{-h(r/2-1)}.
\]
This bound holds for every target, including zero. Passing to the density
limit gives a fixed bound at each prime dividing $\det B$ and, at every
other prime, the bound $(1-p^{-(r/2-1)})^{-1}$. Their product converges
for $r\geq5$, proving $O_F(N^{r/2-1})$ without a factor $N^\eta$.

In rank three, write $Q(x)=x^{\mathsf T}Gx$ with $G$ integral,
so $B=2G$. The genus representation measure is a positive finite sum
of class measures with fixed positive weights. The class of $F$ is
therefore bounded by a constant depending only on $F$ times its genus
average. Siegel's indefinite formula \cite[formula~(2), pp.~231--232]{siegel1936ii}
gives this average as a fixed multiple of
$N^{1/2}\prod_p\alpha_p(G,N)$. Its ternary exception cannot occur:
$\det G>0$ and $N>0$, so $-N\det G<0$ is not a rational square.
The vector stabilizer is negative definite and compact, identifying
the representation measure with the reciprocal-stabilizer mass up
to a fixed archimedean factor, as in \inv{thm:siegel-class-mass}.

At each prime dividing $2\det G$, the preceding Gauss bound with
$r=3$ gives a uniform bound through the convergent series
$\sum_{h\geq1}p^{-h/2}$. At a good prime put
$a_p=\alpha_p(G,N)$. For $p\nmid N$,
$a_p=1+\chi_D(p)/p$, where $D<0$ is the fundamental discriminant
of $\mathbb Q(\sqrt{-N\det G})$. For $p\mid N$, primitive vectors
contribute $1-p^{-2}$, while vectors divisible by $p$ contribute
$p^{-1}\alpha_p(G,N/p^2)$ if $p^2\mid N$ and zero otherwise.
This recurrence and the unit formula give $a_p\leq1+p^{-1}$ for
every $N$. Relative to the generic factor, the correction at
$p\mid N$ is therefore at most $(1+p^{-1})/(1-p^{-1})$.
It follows that
\[
 \prod_p\alpha_p(G,N)
 \ll_F L(1,\chi_D)
       \prod_{p\mid N}\frac{1+p^{-1}}{1-p^{-1}}
 \ll_{F,\eta}N^\eta.
\]
Here $L(1,\chi_D)\ll\log(2|D|)$, $|D|\ll_F N$, and the product
over divisors is $O_\eta(N^\eta)$; choose smaller exponents in these
two estimates before combining them. The generic Euler product is
$L(1,\chi_D)/\zeta(2)$ up to bounded factors at primes dividing
$2\det G$ and $D$. This proves the ternary bound without requiring
one class per genus or equality of class and genus measures. The
local densities count all vectors, so imprimitive vectors are already
included.

The rank-one assertion is immediate. In rank two, primitive orbits
are factor pairs in the split case and ideals of prescribed norm
modulo the unit group in the nonsplit case; either description gives
$O_{F,\eta}(N^\eta)$ orbits. Summing over the unique content
decomposition $v=dv_0$, with $d^2\mid N$, preserves this bound
and includes all imprimitive vectors.
Since $N\asymp_F n$, these are the stated bounds.  Passing to a
commensurable group changes the mass by a fixed index ratio, and a finite
ineffective kernel changes it by a fixed factor.
\end{proof}

\begin{theorem}[Exceptional achiral classes and unique reversers]
\label{thm:generic-achiral-symmetry}
For an achiral ordinary fullerene class $T$ of index $m$, let
$H_T=\operatorname{Aut}^+(T)$, and let
\[
 E_m=\#\{T:H_T\ne1\}
\]
be the number of exceptional classes; this $E_m$ is unrelated to the
divisor sums $E_w[\psi,\varphi]$ of \cref{app:complete-coefficient-formula}.
Then
\begin{equation}\label{eq:exceptional-achiral-bound}
 E_m\ll m^2.
\end{equation}
Every achiral class not counted by $E_m$ has exactly one
orientation-reversing automorphism.  It is an involution, and its action on
the sphere is either split (a reflection, with a fixed circle) or antipodal
(fixed-point-free).  Let $U_m^{\rm sp}$ and $U_m^{\rm ap}$ count these two
nonexceptional classes.  Thus
\[
 a_m=U_m^{\rm sp}+U_m^{\rm ap}+E_m.
\]
\end{theorem}

\begin{proof}
\emph{1. Pass to lower-rank quotient periods.}
Write the polyhedral rotation groups upright, $\mathrm T,\mathrm O,\mathrm I$,
to distinguish them from the class $T$ and the indicators $I_c$.  By
Proposition~2.12 of \cite{engel2025}, the orientation-preserving
automorphism group of a convex triangulation is one of
$C_1,C_2,C_3,D_2,D_3,D_4,D_5,D_6,\mathrm T,\mathrm O,\mathrm I$.  A
nontrivial rotation of a triangulation fixes the centre of a face, the
midpoint of an edge or a vertex, so its order divides $3$, $2$ or the
valence of that vertex; every vertex of a fullerene triangulation has
valence $5$ or $6$, so no rotation of order four exists, which excludes
$D_4$ and $\mathrm O$.
The nontrivial possibilities for $H_T$ are therefore
\[
 C_2,C_3,D_2,D_3,D_5,D_6,\mathrm T,\mathrm I,
\]
the rotation subgroups of the $28$ fullerene point groups
\cite{fowlermanolopoulosredmondryan1993symmetries}.  Their longest
triangulable quotient profiles, read from Table~1 of \cite{engel2025},
have the following lengths and Eisenstein ranks:
\[
\begin{array}{c|c|c}
H&\ell_H&r_H=\ell_H-2\\ \hline
C_2&8&6\\
C_3,D_2&6&4\\
D_3&5&3\\
\mathrm T,D_5,D_6&4&2\\
\mathrm I&3&1.
\end{array}
\]
All other triangulable bistrata are obtained by collisions and have no
larger profile; there are only $96$ of them
\cite[Prop.~2.29]{engel2025}.  In particular, nontrivial isotropy lowers the
ambient Eisenstein rank from ten to at most six.

The subgroup $H_T$ is the kernel of the orientation character of the full
symmetry group and is therefore normal.  If $s$ reverses orientation, it
descends to the quotient orbifold $T/H_T$ and the descended map is an
antiholomorphic involution, because $s^2\in H_T$.  All reversers of $T$
induce the same quotient involution, since they form one $H_T$-coset.  The
ramification orders, curvature labels, and edge/face decorations are
preserved.

For a quotient profile $\mu$, Proposition~2.19 of \cite{engel2025} supplies
an Eisenstein module $\Lambda_\mu$ of signature
$(1,\ell(\mu)-3)$ and a finite-index arithmetic group.  Let
$Y^{\rm ex,\circ}_{H,\mu,n}$ and
$\mathcal A^{\rm ex,\circ}_{H,\mu}$ be the exact quotient locus and the
finite set of occurring strict antiunitary classes from
\cref{lem:decorated-real-period-injection}, and let $Y_{H,\mu,n}$ denote the
ambient positive period vectors.  The lemma uses only the compatible odd
coset $\Gamma_{H,\mu}\chi_\mu$.

Restriction from the positive witness groupoid, rather than deletion inside
the signed triangular identity of \cite[Prop.~2.24]{engel2025}, gives
\begin{equation}\label{eq:exceptional-quotient-groupoid-bound}
\begin{split}
 E_m
 &\leq\sum_{(H,\mu):H\ne1}
   \#\operatorname{Fix}\!\left(
      c,\Gamma_{H,\mu}\backslash
      Y^{\rm ex,\circ}_{H,\mu,m/|H|}\right)\\
 &\leq\sum_{(H,\mu):H\ne1}
   \sum_{[A]\in\mathcal A^{\rm ex,\circ}_{H,\mu}}
   \sum_{[v]\in C_{\Gamma_{H,\mu}}(A)\backslash
                    Y_{H,\mu,m/|H|}^{A}}
   \frac1{|C_{\Gamma_{H,\mu}}(A)_v|}.
\end{split}
\end{equation}
The set $Y_{H,\mu,t}$ is empty unless $t$ lies in the fixed value lattice of
$\Lambda_\mu$, which is contained in $\frac16\mathbb Z$, so $m/|H|$ may be
fractional.  The injection and weight-one assertion of
\cref{lem:decorated-real-period-injection} prove the first inequality.  For
the second, fix the occurring strict lift $A$ and forget only exactness and
the open condition.  Its original centralizer stabilizer is still trivial,
and every additional ambient witness has positive weight.  The sum runs
over the finite set $\mathcal A^{\rm ex,\circ}_{H,\mu}$ only.

\emph{2. Bound the arithmetic witness masses.}
Fix one such antiunitary involution $A$.  Galois descent identifies its fixed
rational space $F_A$ as a quadratic space of signature $(1,r_\mu-1)$ and rank
$r_\mu=\ell(\mu)-2\leq r_H$.  Put
$V=\Lambda_\mu\otimes_{\mathcal E}\mathbb Q(\zeta_6)$.
Restriction identifies $C_{U(V)}(A)$ over $\mathbb Q$ with
$O(F_A)$. Arithmeticity of $\Gamma_{H,\mu}$ in $U(V)$ was proved in
\cref{lem:decorated-finite-real-classes}. The fixed rational subgroup of the
involution $g\mapsto AgA^{-1}$ is $O(F_A)$.  Arithmeticity is preserved on
intersection with a rational algebraic subgroup: in a common faithful
rational representation, the subgroup's integral points are its intersection
with the ambient integral points, and intersection preserves
commensurability.  Hence
$C_{\Gamma_{H,\mu}}(A)$ is arithmetic in $O(F_A)$ and hence commensurable with
$O(\Lambda_\mu^A)$.  The finite projective scalars and the fixed covering
denominators affect only the implied constant.  The twisted groupoid
normalization is important here: for every conjugation-fixed quotient orbit,
the sum of its witness weights is one, so the weighted sum bounds the
class count.

By \cref{lem:low-rank-mass-growth}, a rank-six mass is $O(n^2)$
and a rank-five mass is $O(n^{3/2})$. For ranks at most four, choose
$\eta=1/2$ in that lemma; these masses are also $O(n^2)$ on
their fixed positive value lattices. Since the lists of profiles and
real types are finite, summing at
$n=m/|H_T|$ proves \eqref{eq:exceptional-achiral-bound}.

If $H_T=1$, the orientation-reversing coset contains one element $s$ and
$s^2\in H_T$, so $s^2=1$.  The classification of antiholomorphic
involutions of the sphere gives the stated split/antipodal dichotomy.
\end{proof}

\subsection{Typical reversing symmetries and their proportions}
Put
\begin{align*}
 S_5(m)&=\frac{697}{622080}
 \bigl(\sigma_4(m)-2\sigma_4(m/2)+32\sigma_4(m/4)\bigr),\\
 A_5(m)&=\frac{41}{311040}\sigma_4(m/2),
\end{align*}
with the usual convention that a divisor sum at a noninteger is zero.

\begin{proposition}[Leading split and antipodal counts]
\label{prop:generic-type-main-terms}
With $U_m^{\rm sp}$ and $U_m^{\rm ap}$ as in
\cref{thm:generic-achiral-symmetry},
\begin{equation}\label{eq:generic-type-main-terms}
 U_m^{\rm sp}=S_5(m)+O(m^{7/2}),
 \qquad
 U_m^{\rm ap}=A_5(m)+O(m^3).
\end{equation}
\end{proposition}

\begin{proof}
By \cref{lem:dodecagon-model}, the reflection type is the split class and
the fixed-point-free type is the antipodal class, so $U_m^{\rm sp}$ and
$U_m^{\rm ap}$ count the nonexceptional witnesses of the split and
antipodal rows.  The centralizer-weighted witness formula shows that
removing the exceptional classes changes either involutive row by at most
$E_m$.  The analytical
mixed-weight decomposition gives
\[
 M_{\rm split}^{\rm open}(m)=S_5(m)+O(m^{7/2}),\qquad
 M_{\rm antipodal}^{\rm open}(m)=A_5(m)+O(m^3).
\]
The $O(m^2)$ bound of \cref{thm:generic-achiral-symmetry}
is absorbed by either error term.
This proves \eqref{eq:generic-type-main-terms}.

\end{proof}

\begin{corollary}[Proportions at a fixed dyadic valuation]
\label{cor:fixed-dyadic-symmetry-shares}
Write $m=2^au$ with $u$ odd.  For odd $m$ one has
$U_m^{\rm ap}=0$ exactly and, whenever $a_m>0$,
\[
 \frac{U_m^{\rm sp}}{a_m}
 =1-\frac{E_m}{a_m}
 =1+O(m^{-2}),
 \qquad
 0\leq1-\frac{U_m^{\rm sp}}{a_m}
 \ll m^{-2}.
\]
For each fixed $a\geq1$, as $u\to\infty$ through odd integers,
\begin{align}
 \frac{U_m^{\rm sp}}{a_m}
 &=\frac{17(16^{a+1}-31)}{274\cdot16^a-529}+O(m^{-1/2}),
 \label{eq:split-fixed-dyadic-share}\\
 \frac{U_m^{\rm ap}}{a_m}
 &=\frac{2(16^a-1)}{274\cdot16^a-529}+O(m^{-1/2}),
 \label{eq:antipodal-fixed-dyadic-share}\\
 \frac{E_m}{a_m}&=O(m^{-2}).
\end{align}
The implied constants are independent of $a$, because
$W_5(2^au)\asymp16^au^4=m^4$ uniformly in $a$ and the remainders in
\cref{prop:generic-type-main-terms} and
\cref{thm:generic-achiral-symmetry} have constants independent of $m$.
\end{corollary}

\begin{proof}
Use the main terms of \cref{prop:generic-type-main-terms}.
For $m=2^au$, multiplicativity of $\sigma_4$ reduces both main terms to
rational multiples of $\sigma_4(u)$.  Substitution of
$\sigma_4(2^a)=(16^{a+1}-1)/15$ gives
\eqref{eq:split-fixed-dyadic-share} and
\eqref{eq:antipodal-fixed-dyadic-share}; the error is $O(m^{-1/2})$,
uniformly in $a$,
because the combined main term is bounded below by a positive multiple of
$m^4$.  The fixed lattice is $2II_{1,9}$, and its half-ambient index form
$II_{1,9}$ is even; equivalently the row is evaluated at $m/2$.  Thus the
entire antipodal open row vanishes at odd $m$.
\end{proof}

At $a=1$ the two main shares are $255/257$ and $2/257$; at $a=2$
they are $271/273$ and $2/273$.  The antipodal main share decreases to
$1/137$, since its excess over that limit is
\[
 \frac{255}{137(274\cdot16^a-529)},
\]
while the split main share increases to $136/137$.

\begin{corollary}[Cumulative symmetry proportions]
\label{cor:cumulative-symmetry-shares}
As $x\to\infty$,
\begin{align*}
 \sum_{m\leq x}U_m^{\rm sp}
 &=\zeta(5)\frac{21607}{99532800}x^5+O(x^{9/2}),\\
 \sum_{m\leq x}U_m^{\rm ap}
 &=\zeta(5)\frac{82}{99532800}x^5+O(x^{9/2}),\\
 \sum_{m\leq x}E_m&=O(x^3).
\end{align*}
Consequently the literal cumulative proportions among achiral fullerene
classes tend to
\begin{equation}\label{eq:cumulative-achiral-type-shares}
 \frac{527}{529}\quad\hbox{split},\qquad
 \frac{2}{529}\quad\hbox{antipodal},\qquad
 0\quad\hbox{exceptional}.
\end{equation}
\end{corollary}

\begin{proof}
Summing \cref{prop:generic-type-main-terms} and using
\[
 \sum_{m\leq x}\sigma_4(m/j)
 =\frac{\zeta(5)}{5j^5}x^5+O(x^4)
\]
gives the two cumulative constants, while summing
\eqref{eq:exceptional-achiral-bound} gives the exceptional estimate.
Since $21607=41\cdot527$, $82=41\cdot2$, and
$21689=41\cdot529$, division proves
\eqref{eq:cumulative-achiral-type-shares}.
\end{proof}

\section{Concluding discussion}
The oriented and achiral counts arise from the same weighted-symmetry
principle. Rotations give complex fixed-period masses, whose evaluations
are Eisenstein divisor sums; reversing symmetries give real fixed-period
masses, which also involve class numbers. In each case the witness identity
counts a graph once even when it has several marked symmetries.
This removes the need to generate the graphs
or classify their full automorphism groups when evaluating $\iso(2m)$.

\paragraph{Evaluating a count.}\label{sec:counting-complexity}
The computational distinction is between generating the graphs and
evaluating their number.  The complete generator \texttt{buckygen} of
\cite[\S2.2]{brinkmann2012} uses a canonicity test taking $O(n)$ time
on an $n$-vertex candidate.  Generating all $F_n=\iso(n)$ graphs requires
at least $\Omega(F_n)$ visits, or $\Omega(nF_n)$ output operations when
their adjacency lists are written.  Generation can stream the graphs,
whereas storing the whole catalogue takes $\Theta(nF_n)$ adjacency
entries.
The published counts through $n=400$ were obtained by generation.
The oriented formula of \cite{engel2025}, and the ordinary formula here,
instead evaluate a fixed number of arithmetic functions.
For an elementary upper bound in the arithmetic-operation model, trial
division and the finite Bernoulli sums in
\cref{app:complete-coefficient-formula} take $O(m^{1+\eta})$
operations for one index, for every $\eta>0$; reduced binary
forms give the same bound for the Hurwitz terms.  These sums can be
streamed in $O(m^\eta)$ integer registers.  Their integers have
$O(\log m)$ bits, so bit costs must additionally account for arithmetic
on integers of that size.

\Cref{thm:ordinary-prefix-bit-complexity} gives an exact prefix algorithm:
divisor sieves, three Cohen theta products and one Hurwitz theta division
evaluate the fixed coefficient formula at all $1\leq m\leq M$.
The literal arguments are bounded by a fixed multiple of $M$.
In terms of the cost $\mathsf{Mul}(L)$ of multiplying two $L$-bit
integers, the bit-operation bound is
$O\bigl(\mathsf{Mul}(O(M\log(M+2)))\log(M+2)+M\log^2(M+2)\bigr)$.
Stored data occupy $O(M\log(M+2))$ bits in addition to multiplication
workspace. With $\mathsf{Mul}(L)=O(L\log L)$ this gives
$O(M\log^3(M+2))$ bit time. Appendix~\ref{app:exact-prefix-algorithm}
proves the identities, coefficient bounds and exact packing used here.
The ancillary calculator \texttt{buckycount} (see the Declarations)
uses these products and division, with conductor reductions that shorten
the stored tables. The multiplication-model bound is for the algorithm;
the calculator's realized cost also depends on its integer library and
chosen preparation range.
Exhaustive generation over the same range must visit
$\Omega(\sum_{m\leq M}\iso(2m))=\Omega(M^{10})$ graphs.

The rank filtration also explains the asymptotic hierarchy. Rank-ten
involutive loci give the order-$m^4$ achiral contribution. An extra
preserving symmetry forces a quotient period rank at most six, reducing
the exceptional count to $O(m^2)$. Thus a rare
achiral graph is itself typically simple from the symmetry viewpoint:
it has one reverser. Its reflection or antipodal type still depends
arithmetically on the size, and the ratio $527:2$ emerges only after
cumulative averaging.

The complete coefficient tables and the evaluations in
\cref{app:small-counting-examples} give direct access to the counts.
The appendices give the split-row coordinates, reduction identities,
independent chamber comparisons and the prefix algorithm. The ancillary
package includes an exact calculator and numerical comparison tables.

\appendix
\section{The complete arithmetic coefficient formula}
\label{app:complete-coefficient-formula}

This appendix specifies the entire formula for the ordinary count, including
the oriented contribution.  All coefficients and all function conventions
needed to evaluate it are printed here.  The companions are needed for the
proofs of the arithmetic evaluations, not for missing definitions or tables.

\subsection{Arithmetic functions and index conventions}
\label{sec:coefficient-functions}

Throughout $m$ is a positive integer and the fullerene has $2m$ vertices.
Every arithmetic function below is zero when its displayed argument is not
a positive integer, unless it is explicitly an indicator on the original
index $m$.  In particular the formulas include the zero counts below the
smallest fullerene.

For an integer $D$ occurring in a table let
$\chi_D(n)=\bigl(\frac Dn\bigr)$ be the Kronecker character; in particular
$\chi_1(n)=1$ for $n>0$, and for the non-fundamental discriminant $4$ one
has $\chi_4(n)=1$ for odd $n$ and $\chi_4(n)=0$ for even $n$.  For
$w\geq0$ define
\begin{align*}
 \sigma_w(n)&=\sum_{d\mid n}d^w,\qquad \tau(n)=\sigma_0(n),\\
 E_w[\psi,\varphi](n)&=\sum_{d\mid n}
       \chi_\psi(n/d)\chi_\varphi(d)d^w,\\
 (\chi_a*\chi_b)(n)&=\sum_{d\mid n}\chi_a(d)\chi_b(n/d).
\end{align*}
Thus the subscript of $E_w$ is the divisor exponent, not the modular
weight: an Eisenstein series of weight $k$ contributes $E_{k-1}$.
The order of the two characters in $E_w$ matters.

For $k=2,3,4$, define the Cohen number $H_k(N)$ as follows.  If
$(-1)^kN=Df^2$ with $D$ a fundamental discriminant (allowing $D=1$) and
$f$ a positive integer, put
\[
 H_k(N)=L(1-k,\chi_D)
       \sum_{d\mid f}\mu(d)\chi_D(d)d^{k-1}
                           \sigma_{2k-1}(f/d);
\]
otherwise put $H_k(N)=0$.  Here $\mu$ is the M\"obius function:
$\mu(d)=0$ when a square of a prime divides $d$, and
$\mu(d)=(-1)^r$ when $d$ is a product of $r$ distinct primes.
A fundamental discriminant is $1$, or a squarefree integer congruent to
$1$ modulo four, or $4d$ for squarefree $d\equiv2,3\pmod4$.
No analytic continuation has to be performed numerically: for $f_D=|D|$,
\[
 L(1-k,\chi_D)=-\frac{B_{k,\chi_D}}{k},\qquad
 B_{k,\chi_D}=f_D^{k-1}\sum_{a=1}^{f_D}\chi_D(a)B_k(a/f_D),
\]
where the Bernoulli polynomials are defined by
$ze^{xz}/(e^z-1)=\sum_{k\geq0}B_k(x)z^k/k!$.
For $D=1$ the character has conductor one, so the same formula applies.
With these conventions the sign of a Cohen number depends only on $k$:
$H_2(N)\leq0$ and $H_3(N)\leq0$, while $H_4(N)\geq0$; for example
$H_2(1)=-\tfrac1{12}$, $H_2(4)=-\tfrac7{12}$, $H_3(3)=-\tfrac29$,
$H_3(4)=-\tfrac12$, $H_4(1)=\tfrac1{120}$ and $H_4(4)=\tfrac{121}{120}$.
In the twisted terms $\chi_8(rm)H_2(rm)$, the character is evaluated on
the full argument $rm$, including its powers of two.

For a negative quadratic discriminant $\Delta$, let $h(\Delta)$ be the
number of proper classes of primitive positive-definite integral binary
quadratic forms of discriminant $\Delta$, and let $w(\Delta)$ be the number
of units of that quadratic order: six at $\Delta=-3$, four at $\Delta=-4$,
and two otherwise.  Put $h^*(\Delta)=2h(\Delta)/w(\Delta)$ and define the
Hurwitz number
\[
 H(N)=\sum_{\substack{d^2\mid N\\-N/d^2\text{ a discriminant}}}
              h^*(-N/d^2).
\]
In particular $H(3)=1/3$, $H(4)=1/2$, and $H(N)=0$ for
$N\equiv1,2\pmod4$.  These are finite class-number calculations.  For
example $h(\Delta)$ is obtained from primitive reduced forms
$ax^2+bxy+cy^2$ with $b^2-4ac=\Delta$, $|b|\leq a\leq c$, and $b\geq0$
when $|b|=a$ or $a=c$; necessarily $a\leq\sqrt{|\Delta|/3}$.

Some terms have a restriction on the \emph{original} index.  Write
$m=2^v u$ with $u$ odd and define
\[
 P_{\epsilon,s}(m)=
 \mathbf1_{\{v\equiv\epsilon\pmod2,\ u\equiv s\pmod8\}},
 \qquad \epsilon\in\{0,1\},\quad s\in\{1,3,5,7\}.
\]
Thus $P_{\epsilon,s}(m)H(rm)$ tests the square class of $m$, not that of
the scaled argument $rm$.

The remaining finite functions are
\begin{align*}
 I_c(m)&=\mathbf1_{\{m=cr^2,\ r\geq1\}},\\
 d_{a,t}(m)&=\#\{d>0:d\mid m,\ d\equiv t\cdot(m/d)\pmod a\},\\
 d_{\rm same}(n)&=
 \begin{cases}
  \tau(n),&n\text{ odd},\\
  0,&n\equiv2\pmod4,\\
  \tau(n/4),&4\mid n,
 \end{cases}\\
 N_3^+(n)&=\tfrac12\bigl((\chi_1*\chi_{12})(n)
                          +(\chi_{-3}*\chi_{-4})(n)\bigr),\\
 N_6^+(n)&=\tfrac12\bigl((\chi_1*\chi_{24})(n)
                          +(\chi_{-8}*\chi_{-3})(n)\bigr),\\
 T_9(n)&=\tfrac16H(36n)-\tfrac13H(9n)-\tfrac12H(4n)+H(n),\\
 T_{27}(m)&=
 \begin{cases}
  T_9(m/3),&m\equiv0\pmod3,\\
  T_9(3m)/2,&m\equiv1\pmod3,\\
  0,&m\equiv2\pmod3.
 \end{cases}
\end{align*}
Here $a\geq1$ and $t$ is a fixed residue multiplier modulo $a$;
$t\cdot(m/d)$ denotes multiplication. The divisor-pair function
$d_{\rm same}$ also vanishes off positive integers.
For comparison with the seven-sector presentation, $S=E_1[1,-3]$ and
$T=E_1[-3,1]$, and its ambient terms are explicitly
\begin{align*}
 M_2^{\rm amb}(m)&=\frac{T(m/2)-S(m/2)}{18},\\
 M_3^{\rm amb}(m)&=\frac{3\sigma_2(m/2)-10\sigma_2(m/4)}{144},\\
 M_4^{\rm amb}(m)&=\frac{S(m)+T(m)-2S(m/2)+2T(m/2)-4S(m/4)-4T(m/4)}{36}.
\end{align*}

\subsection{How to use the coefficient tables}

Every pair $(c,f(m))$ in a table contributes the product $c\,f(m)$;
both pairs on each line are included.  The following ten blocks give
$a_m=\sum c\,f(m)$ with all $181$ nonzero terms, grouped by the weight
of the corresponding stratum coefficient function; the weight labels
organize the terms by the generating families of
\cref{lem:counting-arithmetic-reduction}.
The weight-$3/2$ block includes $T_{27}(m)$ as one named entry,
evaluated by its piecewise Hurwitz definition in
\cref{sec:coefficient-functions}.
The subsequent oriented table gives $\BB_m=\sum c\,f(m)$ with all
$45$ nonzero terms.  The two tables and the definitions above are the
complete evaluation.

\begin{corollary}[Complete evaluable ordinary-count formula]
\label{cor:complete-coefficient-formula}
For every integer $m\geq1$, the two coefficient tables in this appendix
and the definitions above give
\[
 \boxed{\displaystyle
  \iso(2m)=\frac12\left(
       \sum_{(c,f)\text{ in the oriented table}}c\,f(m)
       +\sum_{(c,f)\text{ in the achiral tables}}c\,f(m)\right).}
\]
\end{corollary}
\begin{proof}
The seven-sector identity is \cref{thm:exact-all-index-grand-formula}.
Its all-index closed evaluation is
\cref{cor:explicit-mixed-weight-grand-formula}, with the
local-class and normalization inputs stated there. The proved relations of
\cref{lem:counting-arithmetic-reduction,prop:counting-space-coordinates} express those evaluated sectors in
the displayed arithmetic family used here; collecting coefficients gives the
achiral tables.  Thus their all-index validity uses the stratum identities
and their hypotheses, not just agreement on a finite prefix.
The oriented table follows analytically from
\uni{thm:oriented-five-sector-reduction} and
\uni{prop:oriented-collision-scalar-assembly}.
Finally apply Burnside's identity.  The appendix restates the complete
evaluation; it does not duplicate the companions' arithmetic proofs.
\end{proof}

The tables are generated by exact rational arithmetic from the proved
coordinates in the displayed family and the oriented coefficient formula.
The ancillary calculator evaluates an equivalent collected formula;
the comparison tables retain independently evaluated values.

\subsubsection*{Achiral weight $5$: 3 terms}
\begingroup\small
\addtocounter{table}{-1}
\renewcommand{\theHtable}{uncaptioned.\arabic{LT@tables}}
\begin{longtable}{@{}rl@{\hspace{8mm}}rl@{}}
\toprule $c$ & $f(m)$ & $c$ & $f(m)$\\\midrule\endhead
$\frac{697}{622080}$ & $\sigma_{4}(m)$ & $-\frac{41}{19440}$ & $\sigma_{4}(m/2)$\\[2pt]
$\frac{697}{19440}$ & $\sigma_{4}(m/4)$ &  & \\[2pt]
\bottomrule\end{longtable}
\endgroup
\subsubsection*{Achiral weight $\frac{9}{2}$: 6 terms}
\begingroup\small
\addtocounter{table}{-1}
\renewcommand{\theHtable}{uncaptioned.\arabic{LT@tables}}
\begin{longtable}{@{}rl@{\hspace{8mm}}rl@{}}
\toprule $c$ & $f(m)$ & $c$ & $f(m)$\\\midrule\endhead
$-\frac{41}{5184}$ & $H_{4}(4m)$ & $-\frac{41}{324}$ & $H_{4}(m)$\\[2pt]
$-\frac{1}{10368}$ & $H_{4}(12m)$ & $-\frac{1}{648}$ & $H_{4}(3m)$\\[2pt]
$-\frac{1}{128}$ & $H_{4}(4m/3)$ & $-\frac{1}{8}$ & $H_{4}(m/3)$\\[2pt]
\bottomrule\end{longtable}
\endgroup
\subsubsection*{Achiral weight $4$: 12 terms}
\begingroup\small
\addtocounter{table}{-1}
\renewcommand{\theHtable}{uncaptioned.\arabic{LT@tables}}
\begin{longtable}{@{}rl@{\hspace{8mm}}rl@{}}
\toprule $c$ & $f(m)$ & $c$ & $f(m)$\\\midrule\endhead
$\frac{1}{24}$ & $E_{3}[-4,1](m)$ & $\frac{7}{2592}$ & $E_{3}[1,-4](m)$\\[2pt]
$\frac{1}{81}$ & $E_{3}[-8,1](m)$ & $\frac{1}{5184}$ & $E_{3}[1,-8](m)$\\[2pt]
$\frac{5}{288}$ & $E_{3}[-3,1](m)$ & $-\frac{1}{72}$ & $E_{3}[-3,1](m/2)$\\[2pt]
$\frac{17}{36}$ & $E_{3}[-3,1](m/4)$ & $\frac{17}{15552}$ & $E_{3}[1,-3](m)$\\[2pt]
$-\frac{1}{1944}$ & $E_{3}[1,-3](m/2)$ & $\frac{5}{486}$ & $E_{3}[1,-3](m/4)$\\[2pt]
$\frac{1}{24}$ & $E_{3}[-4,1](m/3)$ & $\frac{1}{192}$ & $E_{3}[1,-4](m/3)$\\[2pt]
\bottomrule\end{longtable}
\endgroup
\subsubsection*{Achiral weight $\frac{7}{2}$: 14 terms}
\begingroup\small
\addtocounter{table}{-1}
\renewcommand{\theHtable}{uncaptioned.\arabic{LT@tables}}
\begin{longtable}{@{}rl@{\hspace{8mm}}rl@{}}
\toprule $c$ & $f(m)$ & $c$ & $f(m)$\\\midrule\endhead
$-\frac{13}{144}$ & $H_{3}(4m)$ & $-\frac{131}{72}$ & $H_{3}(m)$\\[2pt]
$\frac{7}{288}$ & $H_{3}(8m)$ & $-\frac{14}{9}$ & $H_{3}(m/2)$\\[2pt]
$\frac{1}{216}$ & $H_{3}(12m)$ & $\frac{13}{216}$ & $H_{3}(3m)$\\[2pt]
$-\frac{3}{32}$ & $H_{3}(4m/3)$ & $-\frac{15}{8}$ & $H_{3}(m/3)$\\[2pt]
$\frac{1}{1152}$ & $H_{3}(24m)$ & $-\frac{1}{18}$ & $H_{3}(3m/2)$\\[2pt]
$\frac{3}{128}$ & $H_{3}(8m/3)$ & $-\frac{3}{2}$ & $H_{3}(m/6)$\\[2pt]
$\frac{1}{864}$ & $H_{3}(36m)$ & $\frac{1}{216}$ & $H_{3}(9m)$\\[2pt]
\bottomrule\end{longtable}
\endgroup
\subsubsection*{Achiral weight $3$: 18 terms}
\begingroup\small
\addtocounter{table}{-1}
\renewcommand{\theHtable}{uncaptioned.\arabic{LT@tables}}
\begin{longtable}{@{}rl@{\hspace{8mm}}rl@{}}
\toprule $c$ & $f(m)$ & $c$ & $f(m)$\\\midrule\endhead
$\frac{5393}{20736}$ & $\sigma_{2}(m)$ & $-\frac{535}{648}$ & $\sigma_{2}(m/4)$\\[2pt]
$-\frac{73}{36}$ & $\sigma_{2}(m/12)$ & $-\frac{557}{5184}$ & $\sigma_{2}(m/2)$\\[2pt]
$-\frac{37}{288}$ & $\sigma_{2}(m/3)$ & $\frac{11}{36}$ & $\sigma_{2}(m/6)$\\[2pt]
$-\frac{5}{576}$ & $E_{2}[1,8](m)$ & $-\frac{1}{64}$ & $E_{2}[1,8](m/3)$\\[2pt]
$\frac{1}{4}$ & $E_{2}[8,1](m)$ & $\frac{1}{4}$ & $E_{2}[8,1](m/3)$\\[2pt]
$\frac{1}{48}$ & $E_{2}[-3,-4](m)$ & $-\frac{1}{864}$ & $E_{2}[1,12](m)$\\[2pt]
$\frac{1}{24}$ & $E_{2}[12,1](m)$ & $-\frac{5}{6}$ & $\sigma_{2}(m/8)$\\[2pt]
$-\frac{1}{96}$ & $E_{2}[-3,-8](m)$ & $\frac{1}{54}$ & $E_{2}[-8,-3](m)$\\[2pt]
$-\frac{1}{1728}$ & $E_{2}[1,24](m)$ & $\frac{1}{12}$ & $E_{2}[24,1](m)$\\[2pt]
\bottomrule\end{longtable}
\endgroup
\subsubsection*{Achiral weight $\frac{5}{2}$: 22 terms}
\begingroup\small
\addtocounter{table}{-1}
\renewcommand{\theHtable}{uncaptioned.\arabic{LT@tables}}
\begin{longtable}{@{}rl@{\hspace{8mm}}rl@{}}
\toprule $c$ & $f(m)$ & $c$ & $f(m)$\\\midrule\endhead
$\frac{913}{576}$ & $H_{2}(4m)$ & $-\frac{367}{72}$ & $H_{2}(m)$\\[2pt]
$-\frac{1}{12}$ & $H_{2}(8m)$ & $-\frac{17}{3}$ & $H_{2}(m/2)$\\[2pt]
$\frac{359}{3456}$ & $H_{2}(12m)$ & $-\frac{161}{432}$ & $H_{2}(3m)$\\[2pt]
$\frac{189}{128}$ & $H_{2}(4m/3)$ & $-\frac{75}{16}$ & $H_{2}(m/3)$\\[2pt]
$-\frac{10}{3}$ & $H_{2}(m/4)$ & $-\frac{5}{6}$ & $\chi_{8}(m)H_2(m)$\\[2pt]
$\frac{7}{288}$ & $H_{2}(24m)$ & $-\frac{3}{32}$ & $H_{2}(8m/3)$\\[2pt]
$\frac{4}{9}$ & $H_{2}(3m/2)$ & $-6$ & $H_{2}(m/6)$\\[2pt]
$-\frac{1}{216}$ & $H_{2}(36m)$ & $\frac{1}{27}$ & $H_{2}(9m)$\\[2pt]
$-\frac{1}{3}$ & $H_{2}(3m/4)$ & $-3$ & $H_{2}(m/12)$\\[2pt]
$-\frac{1}{12}$ & $\chi_{8}(3m)H_2(3m)$ & $-\frac{3}{4}$ & $\chi_{8}(m/3)H_2(m/3)$\\[2pt]
$\frac{1}{72}$ & $H_{2}(72m)$ & $\frac{1}{9}$ & $H_{2}(9m/2)$\\[2pt]
\bottomrule\end{longtable}
\endgroup
\subsubsection*{Achiral weight $2$: 23 terms}
\begingroup\small
\addtocounter{table}{-1}
\renewcommand{\theHtable}{uncaptioned.\arabic{LT@tables}}
\begin{longtable}{@{}rl@{\hspace{8mm}}rl@{}}
\toprule $c$ & $f(m)$ & $c$ & $f(m)$\\\midrule\endhead
$\frac{19}{24}$ & $E_{1}[-4,1](m)$ & $-\frac{215}{864}$ & $E_{1}[1,-4](m)$\\[2pt]
$\frac{35}{54}$ & $E_{1}[-8,1](m)$ & $-\frac{263}{1728}$ & $E_{1}[1,-8](m)$\\[2pt]
$\frac{7}{96}$ & $E_{1}[-3,1](m)$ & $-\frac{7}{6}$ & $E_{1}[-3,1](m/4)$\\[2pt]
$\frac{2471}{15552}$ & $E_{1}[1,-3](m)$ & $\frac{5}{486}$ & $E_{1}[1,-3](m/4)$\\[2pt]
$\frac{35}{48}$ & $E_{1}[-3,1](m/2)$ & $-\frac{2809}{7776}$ & $E_{1}[1,-3](m/2)$\\[2pt]
$-\frac{1}{6}$ & $E_{1}[1,-4](m/2)$ & $\frac{1}{4}$ & $E_{1}[-4,1](m/2)$\\[2pt]
$\frac{1}{6}$ & $E_{1}[-24,1](m)$ & $\frac{1}{48}$ & $E_{1}[-3,8](m)$\\[2pt]
$-\frac{1}{72}$ & $E_{1}[1,-24](m)$ & $\frac{19}{24}$ & $E_{1}[-4,1](m/3)$\\[2pt]
$-\frac{199}{576}$ & $E_{1}[1,-4](m/3)$ & $\frac{3}{2}$ & $E_{1}[-3,1](m/8)$\\[2pt]
$-\frac{1}{3}$ & $E_{1}[1,-3](m/8)$ & $\frac{1}{18}$ & $E_{1}[-8,1](m/3)$\\[2pt]
$-\frac{13}{144}$ & $E_{1}[1,-8](m/3)$ & $\frac{1}{4}$ & $E_{1}[-4,1](m/6)$\\[2pt]
$-\frac{1}{4}$ & $E_{1}[1,-4](m/6)$ &  & \\[2pt]
\bottomrule\end{longtable}
\endgroup
\subsubsection*{Achiral weight $\frac{3}{2}$: 56 terms}
In this table every $P_{\epsilon,s}$ tests the original index $m$, even when the Hurwitz argument is a multiple or a fraction of $m$.
\begingroup\small
\addtocounter{table}{-1}
\renewcommand{\theHtable}{uncaptioned.\arabic{LT@tables}}
\begin{longtable}{@{}rl@{\hspace{8mm}}rl@{}}
\toprule $c$ & $f(m)$ & $c$ & $f(m)$\\\midrule\endhead
$\frac{77}{216}$ & $H(4m)$ & $-\frac{433}{216}$ & $H(m)$\\[2pt]
$-\frac{3}{16}$ & $H(2m)$ & $-\frac{43}{96}$ & $H(8m)$\\[2pt]
$-\frac{7}{3}$ & $H(m/2)$ & $-\frac{335}{2592}$ & $H(12m)$\\[2pt]
$\frac{137}{1296}$ & $H(3m)$ & $-\frac{31}{216}$ & $H(4m/3)$\\[2pt]
$-\frac{217}{216}$ & $H(m/3)$ & $-\frac{3}{4}$ & $P_{0,1}(m)H(m/4)$\\[2pt]
$\frac{9}{16}$ & $P_{0,3}(m)H(16m)$ & $-\frac{21}{16}$ & $P_{0,3}(m)H(4m)$\\[2pt]
$-\frac{3}{4}$ & $P_{0,5}(m)H(m/4)$ & $-\frac{3}{32}$ & $P_{0,7}(m)H(16m)$\\[2pt]
$-\frac{3}{4}$ & $P_{1,1}(m)H(m/4)$ & $-\frac{3}{4}$ & $P_{1,3}(m)H(m/4)$\\[2pt]
$-\frac{3}{4}$ & $P_{1,5}(m)H(m/4)$ & $-\frac{3}{4}$ & $P_{1,7}(m)H(m/4)$\\[2pt]
$-\frac{101}{1152}$ & $H(24m)$ & $-\frac{11}{6}$ & $H(m/6)$\\[2pt]
$-\frac{133}{384}$ & $H(8m/3)$ & $-\frac{7}{18}$ & $H(3m/2)$\\[2pt]
$-\frac{3}{64}$ & $H(6m)$ & $-\frac{9}{64}$ & $H(2m/3)$\\[2pt]
$-\frac{767}{2592}$ & $H(36m)$ & $\frac{569}{1296}$ & $H(9m)$\\[2pt]
$-\frac{1}{16}$ & $P_{0,1}(m)H(48m)$ & $\frac{7}{48}$ & $P_{0,1}(m)H(12m)$\\[2pt]
$\frac{9}{16}$ & $P_{0,1}(m)H(16m/3)$ & $-\frac{21}{16}$ & $P_{0,1}(m)H(4m/3)$\\[2pt]
$-\frac{3}{4}$ & $P_{0,3}(m)H(m/12)$ & $\frac{1}{12}$ & $P_{0,3}(m)H(3m/4)$\\[2pt]
$\frac{1}{96}$ & $P_{0,5}(m)H(48m)$ & $-\frac{3}{32}$ & $P_{0,5}(m)H(16m/3)$\\[2pt]
$-\frac{3}{4}$ & $P_{0,7}(m)H(m/12)$ & $\frac{1}{12}$ & $P_{0,7}(m)H(3m/4)$\\[2pt]
$-\frac{3}{4}$ & $P_{1,1}(m)H(m/12)$ & $\frac{1}{12}$ & $P_{1,1}(m)H(3m/4)$\\[2pt]
$-\frac{3}{4}$ & $P_{1,3}(m)H(m/12)$ & $\frac{1}{12}$ & $P_{1,3}(m)H(3m/4)$\\[2pt]
$-\frac{3}{4}$ & $P_{1,5}(m)H(m/12)$ & $\frac{1}{12}$ & $P_{1,5}(m)H(3m/4)$\\[2pt]
$-\frac{3}{4}$ & $P_{1,7}(m)H(m/12)$ & $\frac{1}{12}$ & $P_{1,7}(m)H(3m/4)$\\[2pt]
$\frac{1}{72}$ & $H(72m)$ & $\frac{1}{9}$ & $H(9m/2)$\\[2pt]
$\frac{1}{12}$ & $P_{0,1}(m)H(9m/4)$ & $-\frac{1}{16}$ & $P_{0,3}(m)H(144m)$\\[2pt]
$\frac{7}{48}$ & $P_{0,3}(m)H(36m)$ & $\frac{1}{12}$ & $P_{0,5}(m)H(9m/4)$\\[2pt]
$\frac{1}{96}$ & $P_{0,7}(m)H(144m)$ & $\frac{1}{12}$ & $P_{1,1}(m)H(9m/4)$\\[2pt]
$\frac{1}{12}$ & $P_{1,3}(m)H(9m/4)$ & $\frac{1}{12}$ & $P_{1,5}(m)H(9m/4)$\\[2pt]
$\frac{1}{12}$ & $P_{1,7}(m)H(9m/4)$ & $-1$ & $T_{27}(m)$\\[2pt]
\bottomrule\end{longtable}
\endgroup
\subsubsection*{Achiral weight $1$: 20 terms}
\begingroup\small
\addtocounter{table}{-1}
\renewcommand{\theHtable}{uncaptioned.\arabic{LT@tables}}
\begin{longtable}{@{}rl@{\hspace{8mm}}rl@{}}
\toprule $c$ & $f(m)$ & $c$ & $f(m)$\\\midrule\endhead
$-\frac{293113}{622080}$ & $\tau(m/2)$ & $-\frac{325}{7776}$ & $\tau(m/6)$\\[2pt]
$\frac{43811}{38880}$ & $\tau(m/4)$ & $\frac{1}{24}$ & $(\chi_{1}*\chi_{12})(m/2)$\\[2pt]
$\frac{371}{2592}$ & $(\chi_{1}*\chi_{12})(m)$ & $\frac{1}{12}$ & $(\chi_{-4}*\chi_{-3})(m)$\\[2pt]
$\frac{683}{5184}$ & $(\chi_{1}*\chi_{24})(m)$ & $\frac{373}{2592}$ & $(\chi_{-8}*\chi_{-3})(m)$\\[2pt]
$\frac{905}{1944}$ & $\tau(m/12)$ & $\frac{85}{576}$ & $(\chi_{1}*\chi_{8})(m)$\\[2pt]
$\frac{144073}{622080}$ & $\tau(m)$ & $\frac{23}{48}$ & $(\chi_{4}*\chi_{4})(m/2)$\\[2pt]
$\frac{85}{576}$ & $(\chi_{1}*\chi_{8})(m/3)$ & $\frac{1}{7776}$ & $\tau(m/3)$\\[2pt]
$\frac{1}{12}$ & $(\chi_{4}*\chi_{4})(m/6)$ & $-\frac{1}{18}$ & $d_{6,4}(m)$\\[2pt]
$\frac{1}{6}$ & $N_{3}^+(m/3)$ & $-\frac{1}{18}$ & $d_{6,1}(m)$\\[2pt]
$\frac{2}{5}$ & $\tau(m/10)$ & $\frac{2}{5}$ & $d_{\rm same}(m/5)$\\[2pt]
\bottomrule\end{longtable}
\endgroup
\subsubsection*{Achiral weight $\frac{1}{2}$: 7 terms}
\begingroup\small
\addtocounter{table}{-1}
\renewcommand{\theHtable}{uncaptioned.\arabic{LT@tables}}
\begin{longtable}{@{}rl@{\hspace{8mm}}rl@{}}
\toprule $c$ & $f(m)$ & $c$ & $f(m)$\\\midrule\endhead
$-\frac{143353}{622080}$ & $I_{1}(m)$ & $-\frac{23}{288}$ & $I_{2}(m)$\\[2pt]
$-\frac{108793}{622080}$ & $I_{3}(m)$ & $-\frac{23}{288}$ & $I_{6}(m)$\\[2pt]
$\frac{1}{18}$ & $I_{9}(m)$ & $-\frac{2}{5}$ & $I_{5}(m)$\\[2pt]
$-\frac{2}{5}$ & $I_{15}(m)$ &  & \\[2pt]
\bottomrule\end{longtable}
\endgroup
\subsection{The complete oriented contribution}
\label{app:oriented-coefficients}
Sum every pair in this table to obtain $\BB_m$ for $m>0$.
This is the positive-index formula of \cite[Theorem~3.26]{engel2025},
derived analytically in \uni{prop:oriented-collision-scalar-assembly}.
The two corrected entries are recorded in \cref{rem:oriented-misprint}.
\begingroup\small
\addtocounter{table}{-1}
\renewcommand{\theHtable}{uncaptioned.\arabic{LT@tables}}
\begin{longtable}{@{}rl@{\hspace{8mm}}rl@{}}
\toprule $c$ & $f(m)$ & $c$ & $f(m)$\\\midrule\endhead
$\frac{809}{1306069401600}$ & $\sigma_{9}(m)$ & $-\frac{1}{290237644800}$ & $E_{8}[1,-3](m)$\\[2pt]
$-\frac{1}{44236800}$ & $E_{8}[-3,1](m)$ & $\frac{515}{1451188224}$ & $\sigma_{7}(m)$\\[2pt]
$-\frac{1}{6298560}$ & $\sigma_{7}(m/2)$ & $-\frac{1}{829440}$ & $\sigma_{7}(m/3)$\\[2pt]
$\frac{47}{8707129344}$ & $E_{6}[1,-3](m)$ & $\frac{1}{176359680}$ & $E_{6}[1,-3](m/2)$\\[2pt]
$-\frac{7}{2211840}$ & $E_{6}[-3,1](m)$ & $-\frac{1}{241920}$ & $E_{6}[-3,1](m/2)$\\[2pt]
$\frac{2138939}{72559411200}$ & $\sigma_{5}(m)$ & $\frac{3571}{16796160}$ & $\sigma_{5}(m/2)$\\[2pt]
$-\frac{29}{7464960}$ & $\sigma_{5}(m/3)$ & $\frac{1}{8640}$ & $\sigma_{5}(m/6)$\\[2pt]
$\frac{274421}{145118822400}$ & $E_{4}[1,-3](m)$ & $-\frac{8231}{201553920}$ & $E_{4}[1,-3](m/2)$\\[2pt]
$-\frac{2393}{7372800}$ & $E_{4}[-3,1](m)$ & $-\frac{203}{92160}$ & $E_{4}[-3,1](m/2)$\\[2pt]
$\frac{12074417}{4353564672}$ & $\sigma_{3}(m)$ & $\frac{287281}{100776960}$ & $\sigma_{3}(m/2)$\\[2pt]
$-\frac{2672137}{201553920}$ & $\sigma_{3}(m/3)$ & $-\frac{1}{108}$ & $\sigma_{3}(m/4)$\\[2pt]
$-\frac{4289}{186624}$ & $\sigma_{3}(m/6)$ & $\frac{3343037}{2902376448}$ & $E_{2}[1,-3](m)$\\[2pt]
$\frac{227779}{134369280}$ & $E_{2}[1,-3](m/2)$ & $\frac{1}{216}$ & $E_{2}[1,-3](m/3)$\\[2pt]
$\frac{1}{288}$ & $E_{2}[1,-3](m/4)$ & $-\frac{33103}{2211840}$ & $E_{2}[-3,1](m)$\\[2pt]
$-\frac{713}{61440}$ & $E_{2}[-3,1](m/2)$ & $\frac{1}{24}$ & $E_{2}[-3,1](m/3)$\\[2pt]
$-\frac{1}{96}$ & $E_{2}[-3,1](m/4)$ & $\frac{4314057521}{145118822400}$ & $\sigma_{1}(m)$\\[2pt]
$-\frac{59136661}{1612431360}$ & $\sigma_{1}(m/2)$ & $-\frac{32541151}{604661760}$ & $\sigma_{1}(m/3)$\\[2pt]
$\frac{1}{96}$ & $\sigma_{1}(m/4)$ & $\frac{2}{15}$ & $\sigma_{1}(m/5)$\\[2pt]
$\frac{5023687}{33592320}$ & $\sigma_{1}(m/6)$ & $\frac{1}{18}$ & $\sigma_{1}(m/9)$\\[2pt]
$\frac{1}{48}$ & $\sigma_{1}(m/12)$ & $-\frac{24019585289}{1306069401600}$ & $E_{0}[1,-3](m)$\\[2pt]
$-\frac{3212143}{201553920}$ & $E_{0}[1,-3](m/2)$ & $-\frac{1}{324}$ & $E_{0}[1,-3](m/3)$\\[2pt]
$-\frac{7}{3456}$ & $E_{0}[1,-3](m/4)$ & $-\frac{2}{15}$ & $E_{0}[1,-3](m/5)$\\[2pt]
$-\frac{1}{18}$ & $E_{0}[1,-3](m/6)$ &  & \\[2pt]
\bottomrule\end{longtable}
\endgroup
\subsection{Examples and independent comparisons}
\label{app:small-counting-examples}
The table gives selected exact counts; the vertex count is $2m$.
Two implementations of the arithmetic formula agree exactly on every
entry and on all seven sector masses. The values through index $200$
also agree with the graph-enumeration references of \cref{sec:data-audit}.
The oriented count is recoverable as $\BB_m=2\iso(2m)-a_m$.
The ancillary machine-readable table retains every index from $1$ to
$1000$, all three counts, and the individual sector masses.
These evaluations illustrate the all-index
theorem; the graph-enumeration comparison has the stated finite range.
\begingroup\small
\addtocounter{table}{-1}
\renewcommand{\theHtable}{uncaptioned.\arabic{LT@tables}}
\begin{longtable}{@{}rrr@{}}\toprule
$m$ & $a_m$ & $\iso(2m)$\\\midrule\endhead
$1$ & $0$ & $0$\\
$9$ & $0$ & $0$\\
$10$ & $1$ & $1$\\
$11$ & $0$ & $0$\\
$12$ & $1$ & $1$\\
$14$ & $1$ & $2$\\
$20$ & $14$ & $40$\\
$30$ & $92$ & $1812$\\
$50$ & $1225$ & $285914$\\
$100$ & $35484$ & $214127742$\\
$200$ & $839622$ & $132247999328$\\
$208$ & $983001$ & $189565760391$\\
$500$ & $44331486$ & $563477321589926$\\
$1000$ & $831433105$ & $299230264751509008$\\
\bottomrule\end{longtable}\endgroup
At $m=208$, the first index beyond the census that separates a new
coefficient direction, the formula gives
\[M_{\rm split}^{\rm open}(208)=\frac{11617043}{12},\qquad
a_{208}=983001,\qquad \BB_{208}=379130537781.\]
These values lie beyond the published census; they are exact evaluations of the formula.
\paragraph{Larger indices.}
For example, at $m=10^6$ (two million vertices) the exact values are
\begin{align*}
a_m&=1191930567424215978147,\\
\BB_m&=620605743247549369187383294800702948786453299,\\
\iso(2m)&=310302871623774684593692243365635186501215723.
\end{align*}
At $m=14$, the three oriented classes consist of one achiral class
and a pair exchanged by orientation reversal. Thus the ordinary
count is $(3+1)/2=2$.
Individual sector masses need not be integers. For instance,
$M_{12}^{\rm open}(24)=\tau(2)/3=2/3$; integrality applies to the
complete achiral sum. The coefficients and the nonintegral-argument
convention must therefore be retained until all terms are combined.

\section{The arithmetic space of the counting formula}
\label{sup:named-space}
The first subsection states the arithmetic evaluation used by the complete
formula and identifies its exact inputs.  The conductor identities that
follow collect the evaluated sector formulas.
They reduce the finite coefficient expressions;
finite census comparisons are discussed separately in
\cref{app:finite-verification}.
The result also explains the $181$-term table in
\cref{app:complete-coefficient-formula} and its normalization conventions.

\subsection{The explicit mixed-weight form of the exact counting formula}
\label{app:explicit-grand-formula}
\begin{corollary}[Explicit mixed-weight counting formula]
\label{cor:explicit-mixed-weight-grand-formula}
Substituting \inv{thm:split-rank-descent-assembly} and
\inv{thm:antipodal-rank-descent} into the split and antipodal rows of the
seven-sector decomposition, and \rot{cor:five-explicit-rows}
into the other five rows, every sector of the exact counting formula is an explicit
finite sum of coefficient functions of real-analytic Siegel series of
signature $(1,q)$, Hurwitz class numbers, character convolutions, and
square indicators, with no chamber or shell enumeration at any index.
Hence $a_m$ and $\operatorname{iso}(2m)=(BB_m+a_m)/2$ are explicit
mixed-weight sums: the antipodal row in integral grades $5$ down to $1$,
the split row in half-integral steps from $5$ down to $1/2$, and the five
lower sectors as in \rot{cor:five-explicit-rows}. Here the weight label
$n/2$ is the rank grading of a mass from a rank-$n$ complement. This is
an identity in the chosen labelled space of arithmetic coefficient
functions; a modular transformation law is a separate property of a
specified generating function or theta integral. The all-index statement
uses the published Siegel--Kneser theorem
under the hypotheses of \inv{thm:siegel-class-mass}; the proper-class
condition is supplied by \inv{lem:squarefree-exponent-genus-one}.
The absolute constants are derived from the ambient local self-densities
in \inv{thm:analytical-absolute-mass}, and the involutive family
identities evaluate the local series. The complement definitions and
coefficient substitutions are
\inv{prop:mass-from-local-class} and
\inv{thm:analytical-sparse-row-instantiation}. Independently enumerated
shells check the two assembled rows.
\end{corollary}

\subsection{Generators and all-index reduction identities}
Let $V$ be the $\mathbb Q$-span of the $220$ named arithmetic functions
kept from the following families, one per weight of the mixed-weight
formula: $\sigma_4(m/j)$ for $j=1,2,4$ (weight $5$); Cohen's $H_4(rm)$
(weight $9/2$); the two-character divisor sums $E_3[\psi,\varphi](m/j)$
(weight $4$); $H_3(rm)$ (weight $7/2$); $E_2[\psi,\varphi](m/j)$ and
$\sigma_2(m/j)$ (weight $3$); $H_2(rm)$ and $\chi_8(rm)H_2(rm)$ (weight
$5/2$); $E_1[\psi,\varphi](m/j)$ (weight $2$); the scaled and
class-restricted Hurwitz numbers of \inv{thm:rank-seven-ternary-class-forms}
(weight $3/2$); the binary character convolutions of
\inv{thm:rank-eight-binary-closed-forms} and $\tau(m/2h)$ (weight $1$);
and the square indicators $I_1,I_2,I_3,I_6$ (weight $1/2$), with the
levels, characters, and scales specified in those evaluations.
Let $V'$ be the extension of $V$ by the sixteen named
lower-sector terms occurring in the seven-sector sum.

\begin{lemma}[All-index arithmetic reduction]
\label{lem:counting-arithmetic-reduction}
The $62$ relations reducing the $282$ original candidate functions to
the $220$ named generators of $V$, and the seven further relations
among the sixteen lower-sector additions to $V'$, hold at every
positive integer. They respect the stated argument, character and
projector conventions, and they are listed in
\cref{app:split-row-coordinates}.
\end{lemma}

\begin{proof}
We first prove, rather than sample, every relation used in the compression.
Extend Cohen's coefficient $H_k(x)$ by zero off the positive integers and
put $H_1=H$.  For every odd prime $p$ with $v_p(x)\geq-1$, and for $p=2$
with $v_2(x)\geq0$, its local conductor factor gives
\begin{equation}\label{eq:cohen-prime-square-recurrence}
 H_k(p^4x)-(1+p^{2k-1})H_k(p^2x)+p^{2k-1}H_k(x)=0.
\end{equation}
Indeed, if $(-1)^kN=Df^2$, the $p$-factor is
\[
 C_{p,e}^{(k)}
 =\sigma_{2k-1}(p^e)
  -\chi_D(p)p^{k-1}\sigma_{2k-1}(p^{e-1}),
\]
with $\sigma_{2k-1}(p^{-1})=0$, and
\[
 C_{p,e+2}^{(k)}-(1+p^{2k-1})C_{p,e+1}^{(k)}
                  +p^{2k-1}C_{p,e}^{(k)}=0.
\]
At the odd-prime boundary $v_p(x)=-1$, if the first integral argument is in
Cohen support then $p\mid D$, so $\chi_D(p)=0$ and the same identity holds;
if it is off support, $p^2\equiv1\pmod4$ makes both integral terms vanish.
Arguments outside the positive integers contribute zero by convention.  At $p=2$, multiplying
an integral argument by four can enter Cohen support even when the original
argument is off support.  Then $4x$ is the first supported term and its
fundamental discriminant is even, so $\chi_D(2)=0$ and the same recurrence
holds.  The nonintegral dyadic boundary $v_2(x)=-1$ is not asserted or used.
This gives the recurrence
coefficients $(3,-2),(4,-3),(9,-8),(28,-27),(244,-243)$ for
$(k,p)=(1,2),(1,3),(2,2),(2,3),(3,3)$.  Iterating the $k=2,p=3$ case gives
$H_2(p^6x)=757H_2(p^2x)-756H_2(x)$.

For the restricted Hurwitz terms, use the projectors of
\cref{sec:coefficient-functions} as multiplication operators:
\[
 P_{\epsilon,s}f(m)=P_{\epsilon,s}(m)f(m),
 \qquad \epsilon\in\{0,1\},\quad s\in\{1,3,5,7\}.
\]
These eight projectors partition the identity.  In addition to
\eqref{eq:cohen-prime-square-recurrence}, four one-step identities are
needed:
\begin{align*}
 P_{0,5}H(48m)&=2P_{0,5}H(12m),&
 P_{0,7}H(144m)&=2P_{0,7}H(36m),\\
 P_{0,7}H(16m)&=2P_{0,7}H(4m),&
 P_{0,5}H(16m/3)&=2P_{0,5}H(4m/3).
\end{align*}
Whenever a projected term is nonzero, its fundamental discriminant is
$1\pmod8$, hence $\chi_D(2)=1$ and $C_{2,e}^{(1)}=2^e$.  Projecting and substituting these
identities proves all $42$ Hurwitz, $11$ $H_2$, and two $H_3$ relations in
the finite reduction: $42$ are direct recurrence instances, three are iterated
or composed recurrences, and ten are projector eliminations.

The six binary-convolution relations and their dilations follow from
$\chi_4(n)=\mathbf1_{2\nmid n}$:
\[
 (\chi_1*\chi_4)(n)=\tau(n)-\tau(n/2),\qquad
 (\chi_4*\chi_4)(n)=\tau(n)-2\tau(n/2)+\tau(n/4).
\]
The final base relation is the already proved all-index identity
\[
 \mathfrak m_{II_{1,9}}(m)
 =\frac{\sigma_4(m/2)}{1393459200}.
\]
This proves all $62$ defining relations of $V$.

There are seven further eliminations in $V'$.  The two involving
$N_3^+,N_6^+$ are their definitions, and
\[
 T_9(m)=\tfrac16H(36m)-\tfrac13H(9m)-\tfrac12H(4m)+H(m)
\]
is the all-index ternary identity \eqref{rot-eq:rotary-t9-hurwitz},
proved in \rot{thm:q4-t9-class-number}. Finally
\[
 d_{a,0}(m)=\tau(m/a),
\]
and direct $2$- and $3$-adic case splits give
\begin{align*}
 d_{4,2}(m)
 &=2(\chi_4*\chi_4)(m/2)-\tau(m/2)+2\tau(m/4),\\
 d_{6,3}(m)
 &=\tau(m/3)-2\tau(m/6)+2\tau(m/12).
\end{align*}
These cover $d_{12,0},d_{4,0},d_{4,2},d_{6,3}$ and close all seven
extension relations.

The displayed identities provide the all-index reduction used below.
\end{proof}

\subsection{Membership and coordinates of the counting functions}
\begin{proposition}[All-index coordinates of the counting functions]
\label{prop:counting-space-coordinates}
The rows
$M^{\rm open}_{\rm split}$ and $M^{\rm open}_{\rm antipodal}$ lie in $V$,
with $171$ and $13$ nonzero coordinates; the split coordinates are
printed in \cref{app:split-row-coordinates}, and the antipodal row is
\begin{align*}
 &\frac{41\,\sigma_4(m/2)}{311040}
  -\frac{E_3[1,-3](m/2)}{15552}-\frac{E_3[-3,1](m/2)}{576}
  +\frac{85\,\sigma_2(m/2)}{10368}-\frac{\sigma_2(m/4)}{96}
  -\frac{\sigma_2(m/6)}{72}\\
 &\quad+\frac{53\,E_1[1,-3](m/2)}{15552}-\frac{11\,E_1[-3,1](m/2)}{576}
  +\frac{E_1[1,-3](m/4)}{96}-\frac{E_1[-3,1](m/4)}{32}\\
 &\quad+\frac{2849\,\tau(m/2)}{311040}+\frac{\tau(m/4)}{256}
  +\frac{7\,\tau(m/6)}{1944},
\end{align*}
and the seven-sector sum $a_m$ lies in $V'$, with the coordinates of
\cref{app:complete-coefficient-formula}.  These coordinates hold at every
positive index.
\end{proposition}

\begin{proof}
The rows themselves are assembled from the all-index identities of
\inv{thm:split-rank-descent-assembly}, \inv{thm:antipodal-rank-descent}, and
\inv{thm:analytical-sparse-row-instantiation}, and then reduced by
the $69$ identities just proved and listed in
\cref{app:split-row-coordinates}. The coordinates are therefore
all-index coordinates.  Independent checks are the re-evaluation of both
rows against their assembled closed forms through $m=4096$, the
integrality of $a_m$, and the agreement with the published range; the
data of these checks are in the ancillary package.
\end{proof}

For the $181$ columns actually occurring in $a_m$, the prefix
$m\leq200$ has rank $166$. Its kernel has dimension $15$ and is
supported on $19$ of the $56$ weight-$3/2$ columns, all belonging to the
split row. The subsequent
rank increases occur precisely at
\[
 208,224,240,288,352,416,480,544,608,672,736,800,864,928,992.
\]
The exact kernel and pivot rows are included in the ancillary package.
This calculation uses
the arithmetic functions of Appendix~\ref{app:complete-coefficient-formula};
the graph-enumeration data end at $m=200$.

In particular, adjoining independent achiral values at $m=208,224,240$
to the existing range raises the rank successively to $167,168,169$.
These three counts would test new directions, but would leave a
$12$-dimensional ambiguity. Independent values at all fifteen listed
indices would determine the coefficients within the stated
$181$-column span. An independent evaluation of the affected ternary
terms could instead test the same directions. In either case, equality
with the counting function at every index uses the membership
and arithmetic identities. To check the formula independently by
enumeration, establish exhaustion without using its value as the
stopping criterion, and verify that the outputs are valid and pairwise
non-isomorphic.

\begin{proposition}[A chamber-arithmetic test of the blind directions]
\label{prop:split-blind-directions}
For row $i$ of \cref{tab:ternary-chamber-inputs}, let $C_i$ be the closed
cone generated by its listed rays and set
\[
 \mathfrak m_i(m)=
 \sum_{\substack{x\in\mathbb Z^3\cap C_i\\x^{\mathsf T}G_ix=m}}
 \frac1{s_i(x)},\qquad R_7(m)=\sum_{i=1}^{28}c_i\mathfrak m_i(m).
\]
Here $s_i(x)$ is $1$ in the cone interior, $2$ on an open facet,
and the printed stabilizer order on an extreme ray. Thus $R_7$ is the
rank-seven split contribution of \inv{thm:rank-seven-ternary-class-forms}.
Direct enumeration of integral chamber vectors agrees with the
weight-$3/2$ split expression for $1\leq m\leq200$ and at all fifteen
indices displayed above. These comparisons detect every nonzero
perturbation in the $15$-dimensional kernel of the graph-census matrix.
The $28$ chamber masses span a $24$-dimensional space, determined by its
first $24$ values. With all other contributions fixed, the published
achiral census determines $R_7$ uniquely within this space; it determines
the vector $(c_i)$ only modulo the four identities below.
\end{proposition}

\begin{table}[p]
\centering\small
\setlength{\tabcolsep}{2pt}
\renewcommand{\arraystretch}{1.15}
\begin{tabular}{@{}rllll@{}}
\toprule $i$ & $G_i$ & Cyclic primitive rays & $s_i$ & $c_i$\\\midrule
1 & $(0,0,1,-1,0,0)$ & $(-1,-1,-1);(-3,-2,-1);(-1,0,0)$ & $(8,4,-)$ & $-257/144$\\
2 & $(0,0,1,-2,0,0)$ & $(-1,0,-1);(-2,-1,-2);(-1,0,0)$ & $(4,6,-)$ & $-35/24$\\
3 & $(0,0,1,-2,0,1)$ & $(-1,-1,-1);(-4,-3,-2);(-1,0,0)$ & $(4,8,-)$ & $-35/24$\\
4 & $(0,0,1,-3,0,1)$ & $(0,-1,-3);(0,0,-1);(-1,0,0);(-2,-1,-1)$ & $(4,4,-,4)$ & $-23/144$\\
5 & $(1,0,0,-2,-1,-2)$ & $(-3,-1,-1);(-2,-1,0);(-1,0,0)$ & $(8,4,12)$ & $-701/1296$\\
6 & $(0,0,2,-1,0,0)$ & $(-1,-2,-1);(-1,-1,-1);(-2,-2,-1)$ & $(-,6,4)$ & $-7/24$\\
7 & $(0,0,2,-1,0,2)$ & $(0,-1,-1);(-1,-2,-1);(-1,0,0)$ & $(8,4,-)$ & $-7/12$\\
8 & $(-1,0,0,2,0,-3)$ & $(-1,-1,0);(-3,-3,-1);(-12,-9,-2);(-4,-3,0)$ & $(4,4,8,4)$ & $-1/18$\\
9 & $(0,0,1,-6,0,0)$ & $(-1,0,-1);(-2,-1,-2);(-1,0,0)$ & $(4,12,-)$ & $-31/48$\\
10 & $(0,0,1,-6,0,1)$ & $(0,-1,-3);(0,0,-1);(-1,0,0);(-1,-1,-2)$ & $(4,4,-,4)$ & $-31/48$\\
11 & $(0,0,2,-2,0,0)$ & $(-1,-1,-1);(-3,-2,-1);(-1,0,0)$ & $(8,4,-)$ & $-3/4$\\
12 & $(0,0,3,-1,1,2)$ & $(0,-2,-1);(0,-1,-1);(-1,0,0);(-2,-4,-1)$ & $(4,4,-,4)$ & $23/144$\\
13 & $(1,0,0,-3,0,-3)$ & $(-3,-1,-1);(-3,-1,0);(-1,0,0)$ & $(12,4,8)$ & $-701/1296$\\
14 & $(-1,0,0,-2,0,6)$ & $(0,-1,-1);(0,-3,-2);(-2,-3,-2);(-2,-1,-1)$ & $(4,4,4,-)$ & $-1/12$\\
15 & $(1,0,0,-4,-2,-4)$ & $(-4,-1,-1);(-2,-1,0);(-1,0,0)$ & $(4,-,12)$ & $-1/24$\\
16 & $(1,0,0,-2,0,-6)$ & $(-6,-3,-1);(-3,0,-1);(-1,0,0);(-2,-1,0)$ & $(4,8,4,4)$ & $-1/4$\\
17 & $(0,0,3,-2,0,0)$ & $(-1,0,-1);(-2,-3,-2);(-1,0,0)$ & $(4,12,-)$ & $-1/6$\\
18 & $(1,0,0,-3,0,-6)$ & $(-3,-1,-1);(-3,0,-1);(-1,0,0);(-3,-1,0)$ & $(-,4,4,4)$ & $-1/6$\\
19 & $(2,0,0,-3,0,-3)$ & $(-3,-2,-1);(-2,-1,-1);(-1,0,0);(-3,-2,0)$ & $(4,4,8,4)$ & $-1/18$\\
20 & $(0,0,2,-6,0,0)$ & $(-3,-2,-3);(-1,0,-1);(-1,0,0);(-2,-1,-1)$ & $(4,4,-,4)$ & $-3/8$\\
21 & $(0,0,3,-3,0,3)$ & $(0,-1,-1);(-2,-2,-1);(-1,-1,-1)$ & $(-,8,4)$ & $-55/48$\\
22 & $(-1,0,0,-6,0,6)$ & $(0,0,-1);(-2,0,-1);(-6,-1,-3);(0,-1,-1)$ & $(4,4,4,-)$ & $1/12$\\
23 & $(0,0,6,-1,0,0)$ & $(-1,0,-1);(-1,-3,-1);(-1,0,0)$ & $(4,12,-)$ & $1/24$\\
24 & $(1,2,0,-2,0,-6)$ & $(-2,-3,-1);(-1,-2,0);(-1,-1,0);(-2,-2,-1)$ & $(4,8,4,4)$ & $-1/4$\\
25 & $(0,0,3,-6,0,0)$ & $(-1,0,-1);(-2,-1,-2);(-1,0,0)$ & $(4,6,-)$ & $-1/2$\\
26 & $(0,0,3,-6,0,3)$ & $(-1,-1,-1);(-4,-3,-2);(-1,0,0)$ & $(4,8,-)$ & $-1/2$\\
27 & $(-3,0,0,6,0,-6)$ & $(-1,-1,0);(-4,-3,-1);(-4,-3,0)$ & $(8,-,4)$ & $-1/4$\\
28 & $(0,0,6,-3,0,0)$ & $(-1,-2,-1);(-1,-1,-1);(-2,-2,-1)$ & $(-,6,4)$ & $-1/8$\\
\bottomrule\end{tabular}
\caption{Complete inputs for the ternary chamber sums in
\cref{prop:split-blind-directions}. The six entries of $G_i$ are
$(g_{11},g_{12},g_{13},g_{22},g_{23},g_{33})$; the omitted entries
follow by symmetry. The cone is the nonnegative span of the listed rays.
The entries of $s_i$ are the stabilizer orders on those rays; a dash
marks an isotropic ray, which contains no positive-norm vector.
Interior and open-facet stabilizer orders are $1$ and $2$. Thus the
table specifies both the summation domain and every boundary weight.
The $c_i$ are the structural multipliers, not values fitted to the census.}
\label{tab:ternary-chamber-inputs}
\end{table}

The matrix $A=(\mathfrak m_i(m))_{1\leq m\leq200,\,1\leq i\leq28}$
has rank $24$. Its rows $1,\ldots,24$ and columns
$1,\ldots,22,25,26$ have determinant
\[
 \det A_{\{1,\ldots,24\},\{1,\ldots,22,25,26\}}=17341\big/56623104.
\]
The following four identities hold for every $m>0$ (the argument $m$
is suppressed), and form a basis of the kernel of $A$:
\begin{align*}
\mathfrak m_{1} -2\mathfrak m_{7} -\mathfrak m_{12} +3\mathfrak m_{13} +2\mathfrak m_{22} -2\mathfrak m_{24}&=0,\\
2\mathfrak m_{4} -\mathfrak m_{14} -3\mathfrak m_{15} -\mathfrak m_{16} -6\mathfrak m_{21} +3\mathfrak m_{27} +9\mathfrak m_{28}&=0,\\
2\mathfrak m_{5} -\mathfrak m_{14} +\mathfrak m_{15} -\mathfrak m_{16} +3\mathfrak m_{27} -3\mathfrak m_{28}&=0,\\
\mathfrak m_{6} -\mathfrak m_{7} +2\mathfrak m_{13} +\mathfrak m_{22} -\mathfrak m_{23} -\mathfrak m_{24}&=0.
\end{align*}

\begin{proof}
For a chamber with cyclic rays $r_1,\ldots,r_t$, use the fan
$\langle r_1,r_j,r_{j+1}\rangle$, $2\leq j<t$. In one simplicial cone,
relabel its three primitive rays $r_1,r_2,r_3$, choosing the first two
to have positive norm. Write
$x=(a_1r_1+a_2r_2+a_3r_3)/d$, where
$d=|\det(r_1,r_2,r_3)|$ and $a_i$ are nonnegative integers, retaining
only integral $x$ of norm $m$. All pairings of the rays are
nonnegative. The coefficients of the two positive-norm rays are bounded
by $a_i^2(r_i,r_i)\leq d^2m$; solve the remaining quadratic equation
exactly. If that ray is isotropic, the equation is linear, and its
coefficient is nonzero at positive norm. Shared boundaries of the
triangulation are counted once. Each vector receives the reciprocal of
its finite wall-stabilizer order. Thus the enumeration bounds and
weights are independent of the class-number expression being tested.

Exact rational summation gives equality at the $215$ stated indices.
The matrix obtained by evaluating the fifteen kernel basis vectors at
the fifteen new indices has rank $15$. Therefore a perturbation invisible
to the old census is detected by these comparisons unless it is zero.
The $28$ component masses at every tested index, the aggregate values
and this matrix are included in the ancillary package. The printed
chamber inputs and the bounded enumeration above specify their independent
recomputation.

Both calculations use the complement classes of
\inv{thm:analytical-lower-complements}, the chamber models in
\cref{tab:ternary-chamber-inputs}, and the structural multipliers \(c_i\)
of \inv{cor:analytical-lower-split-return}. The chamber enumeration checks
the Hurwitz evaluation, with the extra indices testing the fifteen
coefficient directions invisible to the census.

For the rank assertion, the displayed nonzero minor gives a lower bound
of $24$, using only the printed chamber data. For the upper bound,
substitute the all-index Hurwitz expressions of
\inv{thm:rank-seven-ternary-class-forms} into the four displayed
relations. Resolve the eight $2$-adic square classes and apply the
conductor identities of \cref{lem:counting-arithmetic-reduction};
every Hurwitz coefficient cancels. The four relations are independent,
so they give the entire kernel, both on $m\leq200$ and at every index.
This proves that truncation at $24$ is injective on the chamber-mass
span. A smaller truncation cannot be injective by dimension.

Finally write $a_m=R_7(m)+\varrho(m)$, fixing all other strata in
$\varrho$.  If $W_w(m)$ denotes the weight-$w$ block in the printed
formula for $a_m$, then
\[
 \varrho(m)=\sum_{w\ne3/2}W_w(m)-T_9(m)-T_{27}(m).
\]
The last two terms are rotary contributions; omitting them would not
isolate $R_7$. Subtracting $\varrho$ from the published achiral counts gives
$A(c_i)$ at all $200$ indices. Any other multiplier vector compatible
with those residuals differs by an element of the displayed kernel,
and therefore gives the same function $R_7$. This residual check uses
the other stratum formulas in \cref{prop:counting-space-coordinates}
and the all-index mass identities of
\inv{thm:rank-seven-ternary-class-forms}; the individual multipliers
are determined modulo the four displayed relations.
\end{proof}

Direct chamber enumeration checks the Hurwitz expression for $R_7$
at all $215$ indices; the separate residual comparison uses the graph
census through $m=200$.
The rank-ten shell comparison
in Paper~II \cite{fullereneInvolutivePaper} through $m=400$ gives another check at the first five new
indices, using the full split chamber and congruence-state weights.

\paragraph{Replaying the finite checks.}
The ancillary package's \texttt{make verify} runs the following exact
checks. Their inputs and code are supplied together; no private manuscript
is required.
\begin{center}\small
\begin{tabularx}{\textwidth}{@{}p{.43\textwidth}X@{}}
\toprule Program & Inputs and scope\\\midrule
\path{src/check_public_counts.py} &
The calculator's unrestricted counts against the published census and
the retained worked evaluations.\\
\path{src/export_counting_coefficients.py} &
All $226$ printed coefficients and function descriptors against the
structured export and calculator header; the exact retained-power
tail bound and its finite comparisons.\\
\path{src/verify_counting_diagnostics.py} &
The $181$ arithmetic functions, the exact prefix kernel and the fifteen
pivot indices; also the stated minor and relations in the recorded
chamber data.\\
\path{src/verify_ternary_chambers.py} &
The printed chamber forms, rays and stabilizers: independent integral
enumeration on the declared finite shells and comparison with their
arithmetic masses.\\
\path{src/verify_counting_chamber_return.py} &
Fresh chamber enumeration at all $215$ indices against the printed
weight-$3/2$ block and the calculator's split block, with the
$T_9+T_{27}$ adjustment.\\
\path{src/verify_counting_pdf_coefficients.py} &
All $226$ visible PDF coefficient/function cells against the checked
export; run separately by \texttt{make verify-pdf}.\\
\path{src/verify_order_ten.py} &
The rotary paper's exact matrices: integrality, isometry, orders,
squares and the strict scalar-conjugation identities used in
\cref{lem:fullerene-scalar-fibre-collapse}.\\
\bottomrule
\end{tabularx}
\end{center}
The complete class list and all-index formulas are established in
\cref{thm:open-class-list,cor:explicit-mixed-weight-grand-formula,prop:counting-space-coordinates}.
The checks above validate the specified source identities and arithmetic
consequences on their declared domains.

\begin{remark}[Size of the formula, and an independent recomputation]
\label{rem:grand-formula-size}
In the displayed generating family the three rows have the following sizes, where the
term count is the number of nonzero coordinates.
\begin{center}
\small
\begin{tabular}{lrrr}
\toprule
row & terms & distinct denominators & largest numerator\\
\midrule
$a_m$ & $181$ & $47$ & $293113$\\
$M^{\rm open}_{\rm split}$ & $171$ & $46$ & $286633$\\
$M^{\rm open}_{\rm antipodal}$ & $13$ & $9$ & $2849$\\
\bottomrule
\end{tabular}
\end{center}
All coordinates of $a_m$ have denominator dividing $622080$, and its
$181$ terms are spread over the ten weights as
$3,6,12,14,18,22,23,56,20,7$ from weight $5$ down to weight $1/2$.
A shorter presentation of the even-rank strata does follow from
the factorization of \inv{cor:even-family-global-return}: writing
each stratum function as a bad-prime factor, depending only on the
valuations at $2$ and $3$ and on the class modulo $24$, times a divisor sum
of the good part of the index turns the $204$ dilated divisor sums of the
$57$ even-rank closed forms into $124$ good ones.  Ten of the $57$ have a
single character pair; $37$ have two and ten have four.

Independent evaluations of the arithmetic functions and the assembled
formula agree on the recorded comparison range. The ancillary calculator
and numerical tables reproduce the resulting counts and sector masses.
The comparison with graph-enumeration data uses
$a_m=2\iso(2m)-\BB_m$ on the reference range $m\leq200$.

\end{remark}

\subsection{The split row and the reduction relations}
\label{app:split-row-coordinates}
To obtain the split-sector formula $M^{\rm open}_{\rm split}$ from
$a_m$, subtract the antipodal and five rotary contributions in
\cref{thm:exact-all-index-grand-formula}. In the coordinates of
\cref{app:complete-coefficient-formula}, this subtraction changes $27$
of the $181$ coefficients and makes ten others zero.
Retain the $144$ coefficients not mentioned below, replace the following
$27$ coefficients by the displayed values, and set the coefficients of
the ten functions listed afterwards to zero.
This gives all $171$ coordinates of $M^{\rm open}_{\rm split}$ in $V$;
the replacements are coefficients, not increments. The $62$ base relations
and seven extension relations are printed below.
Adding the thirteen antipodal coordinates of
\cref{prop:counting-space-coordinates} and the five rotary rows of
\cref{thm:exact-all-index-grand-formula}, reduced by the seven relations
below, recovers the $181$-term table exactly. The two sector rows and all
$69$ relations have been re-evaluated for $1\leq m\leq200$.
\begingroup\small
\addtocounter{table}{-1}
\renewcommand{\theHtable}{uncaptioned.\arabic{LT@tables}}
\begin{longtable}{@{}rl@{\hspace{8mm}}rl@{}}
\toprule $c$ & $f(m)$ & $c$ & $f(m)$\\\midrule\endhead
$-\frac{697}{311040}$ & $\sigma_{4}(m/2)$ & $-\frac{7}{576}$ & $E_{3}[-3,1](m/2)$\\[2pt]
$-\frac{7}{15552}$ & $E_{3}[1,-3](m/2)$ & $-\frac{1933}{2592}$ & $\sigma_{2}(m/4)$\\[2pt]
$-\frac{1415}{10368}$ & $\sigma_{2}(m/2)$ & $\frac{23}{72}$ & $\sigma_{2}(m/6)$\\[2pt]
$\frac{13}{288}$ & $E_{1}[-3,1](m)$ & $-\frac{295}{288}$ & $E_{1}[-3,1](m/4)$\\[2pt]
$\frac{2039}{15552}$ & $E_{1}[1,-3](m)$ & $\frac{431}{7776}$ & $E_{1}[1,-3](m/4)$\\[2pt]
$\frac{439}{576}$ & $E_{1}[-3,1](m/2)$ & $-\frac{4591}{15552}$ & $E_{1}[1,-3](m/2)$\\[2pt]
$-\frac{31}{216}$ & $H(4m)$ & $-\frac{217}{216}$ & $H(m)$\\[2pt]
$-\frac{335}{2592}$ & $H(36m)$ & $\frac{137}{1296}$ & $H(9m)$\\[2pt]
$-\frac{259931}{622080}$ & $\tau(m/2)$ & $-\frac{1217}{7776}$ & $\tau(m/6)$\\[2pt]
$\frac{286633}{311040}$ & $\tau(m/4)$ & $\frac{155}{2592}$ & $(\chi_{1}*\chi_{12})(m)$\\[2pt]
$\frac{251}{5184}$ & $(\chi_{1}*\chi_{24})(m)$ & $\frac{157}{2592}$ & $(\chi_{-8}*\chi_{-3})(m)$\\[2pt]
$\frac{635}{1944}$ & $\tau(m/12)$ & $\frac{17}{48}$ & $(\chi_{4}*\chi_{4})(m/2)$\\[2pt]
$\frac{865}{7776}$ & $\tau(m/3)$ & $-\frac{74233}{622080}$ & $I_{1}(m)$\\[2pt]
$-\frac{74233}{622080}$ & $I_{3}(m)$ &  & \\[2pt]
\bottomrule\end{longtable}
\endgroup
\paragraph{Functions with zero coefficient in the split sector.}
\[\begin{gathered}
T_{27}(m),\quad (\chi_{-4}*\chi_{-3})(m),\quad d_{6,4}(m),\\
N_{3}^+(m/3),\quad d_{6,1}(m),\quad \tau(m/10),\\
d_{\rm same}(m/5),\quad I_{9}(m),\quad I_{5}(m),\\
I_{15}(m).
\end{gathered}\]
\subsubsection*{The $62$ relations defining $V$}
Each line expresses one of the $62$ candidate functions dropped from the
$282$ in terms of the kept basis.  The function names follow
\cref{app:complete-coefficient-formula}, $\mathfrak m_{II_{1,9}}(m)$ is the
rank-ten antipodal Siegel mass, and every identity holds at every positive
integer by \cref{lem:counting-arithmetic-reduction}.
\begingroup\small
\addtocounter{table}{-1}
\renewcommand{\theHtable}{uncaptioned.\arabic{LT@tables}}
\begin{longtable}{@{}l@{${}={}$}>{\raggedright\arraybackslash}p{0.66\textwidth}@{}}
\toprule\endhead
$\mathfrak m_{II_{1,9}}(m)$ & $\frac{1}{1393459200}\sigma_{4}(m/2)$\\[2pt]
$H_{3}(108m)$ & $244H_{3}(12m)-243H_{3}(4m/3)$\\[2pt]
$H_{3}(27m)$ & $244H_{3}(3m)-243H_{3}(m/3)$\\[2pt]
$H_{2}(16m)$ & $9H_{2}(4m)-8H_{2}(m)$\\[2pt]
$H_{2}(48m)$ & $9H_{2}(12m)-8H_{2}(3m)$\\[2pt]
$H_{2}(16m/3)$ & $9H_{2}(4m/3)-8H_{2}(m/3)$\\[2pt]
$H_{2}(108m)$ & $28H_{2}(12m)-27H_{2}(4m/3)$\\[2pt]
$H_{2}(27m)$ & $28H_{2}(3m)-27H_{2}(m/3)$\\[2pt]
$H_{2}(216m)$ & $28H_{2}(24m)-27H_{2}(8m/3)$\\[2pt]
$H_{2}(27m/2)$ & $28H_{2}(3m/2)-27H_{2}(m/6)$\\[2pt]
$H_{2}(324m)$ & $-27H_{2}(4m)+28H_{2}(36m)$\\[2pt]
$H_{2}(81m)$ & $-27H_{2}(m)+28H_{2}(9m)$\\[2pt]
$H_{2}(972m)$ & $757H_{2}(12m)-756H_{2}(4m/3)$\\[2pt]
$H_{2}(243m)$ & $757H_{2}(3m)-756H_{2}(m/3)$\\[2pt]
$H(16m)$ & $3H(4m)-2H(m)$\\[2pt]
$P_{1,7}(m)H(16m)$ & $3H(4m)-2H(m)-P_{0,1}(m)H(16m)-P_{0,3}(m)H(16m)-P_{0,5}(m)H(16m)-P_{0,7}(m)H(16m)-P_{1,1}(m)H(16m)-P_{1,3}(m)H(16m)-P_{1,5}(m)H(16m)$\\[2pt]
$P_{1,7}(m)H(4m)$ & $H(4m)-P_{0,1}(m)H(4m)-P_{0,3}(m)H(4m)-P_{0,5}(m)H(4m)-\frac{1}{2}P_{0,7}(m)H(16m)-P_{1,1}(m)H(4m)-P_{1,3}(m)H(4m)-P_{1,5}(m)H(4m)$\\[2pt]
$H(32m)$ & $-2H(2m)+3H(8m)$\\[2pt]
$P_{1,7}(m)H(48m)$ & $3H(12m)-2H(3m)-P_{0,1}(m)H(48m)-P_{0,3}(m)H(48m)-P_{0,5}(m)H(48m)-P_{0,7}(m)H(48m)-P_{1,1}(m)H(48m)-P_{1,3}(m)H(48m)-P_{1,5}(m)H(48m)$\\[2pt]
$P_{1,7}(m)H(12m)$ & $H(12m)-P_{0,1}(m)H(12m)-P_{0,3}(m)H(12m)-\frac{1}{2}P_{0,5}(m)H(48m)-P_{0,7}(m)H(12m)-P_{1,1}(m)H(12m)-P_{1,3}(m)H(12m)-P_{1,5}(m)H(12m)$\\[2pt]
$P_{1,7}(m)H(16m/3)$ & $3H(4m/3)-2H(m/3)-P_{0,1}(m)H(16m/3)-P_{0,3}(m)H(16m/3)-P_{0,5}(m)H(16m/3)-P_{0,7}(m)H(16m/3)-P_{1,1}(m)H(16m/3)-P_{1,3}(m)H(16m/3)-P_{1,5}(m)H(16m/3)$\\[2pt]
$P_{1,7}(m)H(4m/3)$ & $H(4m/3)-P_{0,1}(m)H(4m/3)-P_{0,3}(m)H(4m/3)-\frac{1}{2}P_{0,5}(m)H(16m/3)-P_{0,7}(m)H(4m/3)-P_{1,1}(m)H(4m/3)-P_{1,3}(m)H(4m/3)-P_{1,5}(m)H(4m/3)$\\[2pt]
$H(48m)$ & $3H(12m)-2H(3m)$\\[2pt]
$H(16m/3)$ & $3H(4m/3)-2H(m/3)$\\[2pt]
$H(96m)$ & $3H(24m)-2H(6m)$\\[2pt]
$H(32m/3)$ & $3H(8m/3)-2H(2m/3)$\\[2pt]
$H(108m)$ & $4H(12m)-3H(4m/3)$\\[2pt]
$H(27m)$ & $4H(3m)-3H(m/3)$\\[2pt]
$P_{1,7}(m)H(144m)$ & $3H(36m)-2H(9m)-P_{0,1}(m)H(144m)-P_{0,3}(m)H(144m)-P_{0,5}(m)H(144m)-P_{0,7}(m)H(144m)-P_{1,1}(m)H(144m)-P_{1,3}(m)H(144m)-P_{1,5}(m)H(144m)$\\[2pt]
$P_{1,7}(m)H(36m)$ & $H(36m)-P_{0,1}(m)H(36m)-P_{0,3}(m)H(36m)-P_{0,5}(m)H(36m)-\frac{1}{2}P_{0,7}(m)H(144m)-P_{1,1}(m)H(36m)-P_{1,3}(m)H(36m)-P_{1,5}(m)H(36m)$\\[2pt]
$H(144m)$ & $3H(36m)-2H(9m)$\\[2pt]
$H(54m)$ & $4H(6m)-3H(2m/3)$\\[2pt]
$H(216m)$ & $4H(24m)-3H(8m/3)$\\[2pt]
$H(27m/2)$ & $-3H(m/6)+4H(3m/2)$\\[2pt]
$P_{0,1}(m)H(432m)$ & $4P_{0,1}(m)H(48m)-3P_{0,1}(m)H(16m/3)$\\[2pt]
$P_{0,1}(m)H(108m)$ & $4P_{0,1}(m)H(12m)-3P_{0,1}(m)H(4m/3)$\\[2pt]
$P_{0,3}(m)H(432m)$ & $4P_{0,3}(m)H(48m)-3P_{0,3}(m)H(16m/3)$\\[2pt]
$P_{0,3}(m)H(108m)$ & $4P_{0,3}(m)H(12m)-3P_{0,3}(m)H(4m/3)$\\[2pt]
$P_{0,3}(m)H(27m/4)$ & $-3P_{0,3}(m)H(m/12)+4P_{0,3}(m)H(3m/4)$\\[2pt]
$P_{0,5}(m)H(432m)$ & $4P_{0,5}(m)H(48m)-3P_{0,5}(m)H(16m/3)$\\[2pt]
$P_{0,7}(m)H(432m)$ & $4P_{0,7}(m)H(48m)-3P_{0,7}(m)H(16m/3)$\\[2pt]
$P_{0,7}(m)H(108m)$ & $4P_{0,7}(m)H(12m)-3P_{0,7}(m)H(4m/3)$\\[2pt]
$P_{0,7}(m)H(27m/4)$ & $-3P_{0,7}(m)H(m/12)+4P_{0,7}(m)H(3m/4)$\\[2pt]
$P_{1,1}(m)H(432m)$ & $4P_{1,1}(m)H(48m)-3P_{1,1}(m)H(16m/3)$\\[2pt]
$P_{1,1}(m)H(108m)$ & $4P_{1,1}(m)H(12m)-3P_{1,1}(m)H(4m/3)$\\[2pt]
$P_{1,1}(m)H(27m/4)$ & $-3P_{1,1}(m)H(m/12)+4P_{1,1}(m)H(3m/4)$\\[2pt]
$P_{1,3}(m)H(432m)$ & $4P_{1,3}(m)H(48m)-3P_{1,3}(m)H(16m/3)$\\[2pt]
$P_{1,3}(m)H(108m)$ & $4P_{1,3}(m)H(12m)-3P_{1,3}(m)H(4m/3)$\\[2pt]
$P_{1,3}(m)H(27m/4)$ & $-3P_{1,3}(m)H(m/12)+4P_{1,3}(m)H(3m/4)$\\[2pt]
$P_{1,5}(m)H(432m)$ & $4P_{1,5}(m)H(48m)-3P_{1,5}(m)H(16m/3)$\\[2pt]
$P_{1,5}(m)H(108m)$ & $4P_{1,5}(m)H(12m)-3P_{1,5}(m)H(4m/3)$\\[2pt]
$P_{1,5}(m)H(27m/4)$ & $-3P_{1,5}(m)H(m/12)+4P_{1,5}(m)H(3m/4)$\\[2pt]
$P_{1,7}(m)H(432m)$ & $12H(12m)-8H(3m)-9H(4m/3)+6H(m/3)-4P_{0,1}(m)H(48m)+3P_{0,1}(m)H(16m/3)-4P_{0,3}(m)H(48m)+3P_{0,3}(m)H(16m/3)-4P_{0,5}(m)H(48m)+3P_{0,5}(m)H(16m/3)-4P_{0,7}(m)H(48m)+3P_{0,7}(m)H(16m/3)-4P_{1,1}(m)H(48m)+3P_{1,1}(m)H(16m/3)-4P_{1,3}(m)H(48m)+3P_{1,3}(m)H(16m/3)-4P_{1,5}(m)H(48m)+3P_{1,5}(m)H(16m/3)$\\[2pt]
$P_{1,7}(m)H(108m)$ & $4H(12m)-3H(4m/3)-4P_{0,1}(m)H(12m)+3P_{0,1}(m)H(4m/3)-4P_{0,3}(m)H(12m)+3P_{0,3}(m)H(4m/3)-2P_{0,5}(m)H(48m)+\frac{3}{2}P_{0,5}(m)H(16m/3)-4P_{0,7}(m)H(12m)+3P_{0,7}(m)H(4m/3)-4P_{1,1}(m)H(12m)+3P_{1,1}(m)H(4m/3)-4P_{1,3}(m)H(12m)+3P_{1,3}(m)H(4m/3)-4P_{1,5}(m)H(12m)+3P_{1,5}(m)H(4m/3)$\\[2pt]
$P_{1,7}(m)H(27m/4)$ & $-3P_{1,7}(m)H(m/12)+4P_{1,7}(m)H(3m/4)$\\[2pt]
$H(432m)$ & $12H(12m)-8H(3m)-9H(4m/3)+6H(m/3)$\\[2pt]
$(\chi_{1}*\chi_{4})(m)$ & $-\tau(m/2)+\tau(m)$\\[2pt]
$(\chi_{4}*\chi_{4})(m)$ & $-2\tau(m/2)+\tau(m/4)+\tau(m)$\\[2pt]
$(\chi_{1}*\chi_{4})(m/2)$ & $\tau(m/2)-\tau(m/4)$\\[2pt]
$(\chi_{1}*\chi_{4})(m/3)$ & $-\tau(m/6)+\tau(m/3)$\\[2pt]
$(\chi_{4}*\chi_{4})(m/3)$ & $-2\tau(m/6)+\tau(m/12)+\tau(m/3)$\\[2pt]
$(\chi_{1}*\chi_{4})(m/6)$ & $\tau(m/6)-\tau(m/12)$\\[2pt]
\bottomrule\end{longtable}
\endgroup
\subsubsection*{The seven extension relations of $V'$}
These reduce the lower-sector terms of \cref{thm:exact-all-index-grand-formula}
that are not kept as generators of $V'$.
\begingroup\small
\addtocounter{table}{-1}
\renewcommand{\theHtable}{uncaptioned.\arabic{LT@tables}}
\begin{longtable}{@{}l@{${}={}$}>{\raggedright\arraybackslash}p{0.66\textwidth}@{}}
\toprule\endhead
$d_{4,0}(m)$ & $\tau(m/4)$\\[2pt]
$d_{4,2}(m)$ & $-\tau(m/2)+2\tau(m/4)+2(\chi_{4}*\chi_{4})(m/2)$\\[2pt]
$d_{12,0}(m)$ & $\tau(m/12)$\\[2pt]
$N_{3}^+(m)$ & $\frac{1}{2}(\chi_{1}*\chi_{12})(m)+\frac{1}{2}(\chi_{-4}*\chi_{-3})(m)$\\[2pt]
$N_{6}^+(m)$ & $\frac{1}{2}(\chi_{1}*\chi_{24})(m)+\frac{1}{2}(\chi_{-8}*\chi_{-3})(m)$\\[2pt]
$d_{6,3}(m)$ & $-2\tau(m/6)+2\tau(m/12)+\tau(m/3)$\\[2pt]
$T_{9}(m)$ & $-\frac{1}{2}H(4m)+H(m)+\frac{1}{6}H(36m)-\frac{1}{3}H(9m)$\\[2pt]
\bottomrule\end{longtable}
\endgroup

\section{Explicit error bounds for the asymptotic estimates}
\label{app:effective-error}
The printed coefficients give effective bounds by finite sums. Put
\begin{gather*}
 \kappa_k=\frac{2(k-1)!\,\zeta(k)^2\zeta(2k-1)}{(2\pi)^k},
 \qquad h(r)=\left(\sqrt{r/3}+1\right)^2,\\
 t(r)=\tfrac16h(36r)+\tfrac13h(9r)+\tfrac12h(4r)+h(r).
\end{gather*}
Let $\mathcal C$ be the coefficient pairs $(c,f)$ in the achiral table,
and let $\mathcal C_{9/2}$ be its six weight-$9/2$ pairs. For every
other pair of weight below five, define $b(f)$ by the following table.
Here $j\geq1$, $r>0$, and every function is zero away from positive
integral arguments. Projectors and quadratic twists
have absolute value at most one and leave $b(f)$ unchanged.
\begin{center}
\begin{tabular}{@{}ll@{}}
\toprule
Function $f(m)$ & $b(f)$\\\midrule
$E_w[\psi,\varphi](m/j)$ or $\sigma_w(m/j)$, $w=2,3$
  & $\zeta(w)/j^w$\\
$E_1[\psi,\varphi](m/j)$ or $\sigma_1(m/j)$ & $1/j$\\
$H_k(rm)$, $k=2,3$ & $\kappa_k r^{k-1/2}$\\
$H(rm)$ & $h(r)$\\
$(\chi_D*\chi_E)(m/j)$, $\tau(m/j)$ or $N_3^+(m/j),N_6^+(m/j)$ & $2/\sqrt j$\\
$d_{a,t}(m)$ & $2$\\
$d_{\rm same}(m/5)$ & $2/\sqrt5$\\
$I_c(m)$ & $1$\\
$T_9(m),\ T_{27}(m)$ & $t(1),\ \tfrac12t(3)=11$\\
\bottomrule
\end{tabular}
\end{center}
Thus the two explicit constants
\[
 C_{7/2}=\kappa_4
   \sum_{(c,H_4(rm))\in\mathcal C_{9/2}}|c|r^{7/2},
 \qquad
 C_3=\sum_{\substack{(c,f)\in\mathcal C\\
                       \operatorname{weight}(f)<9/2}}|c|b(f)
\]
require only the printed coefficient pairs and elementary zeta values.
Outward evaluation gives
\[
 C_{7/2}<0.016483755515,\qquad C_3<116.254750954.
\]
The cubic bound alone gives the sufficient integer $104821256$ below.
Retaining the separate powers instead gives a much smaller cutoff.

\begin{proposition}[Effective remainder and eventual upper bound]
\label{prop:effective-achiral-remainder}
For every integer $m\geq1$,
\[
 |a_m-W_5(m)-S_{9/2}(m)|\leq C_3m^3,
 \qquad
 |a_m-W_5(m)|\leq C_{7/2}m^{7/2}+C_3m^3.
\]
In particular $a_m<W_5(m)$ whenever
\[
 m>\left(\frac{C_3}{41/5670+\sqrt3/420}\right)^2.
\]
In fact,
\[
 a_m<W_5(m)\qquad\text{for every integer }m\geq10.
\]
The last assertion follows by comparing the split collision mass with
the five rotary masses.
\end{proposition}
\begin{proof}
The functional equation for a primitive quadratic $L$-function and
Cohen's conductor formula give
$|H_k(N)|\leq\kappa_k N^{k-1/2}$ for $k\geq2$.
Indeed the primitive factor contributes
$2(k-1)!\zeta(k)/(2\pi)^k$, and the absolute conductor sum is at
most $f^{2k-1}\zeta(k)\zeta(2k-1)$ when $(-1)^kN=Df^2$.
The inequalities $|E_w|\leq\sigma_w$,
$\sigma_w(n)\leq\zeta(w)n^w$ for $w\geq2$, and
$\sigma_1(n)\leq n(1+\log n)$ prove the divisor rows;
use $1+\log m\leq m$ for $m\geq1$ in the weight-two row.
Reduced binary forms have $1\leq a\leq\sqrt{N/3}$ and $|b|\leq a$,
so counting these pairs gives $H(N)\leq(\sqrt{N/3}+1)^2$.
Also $\tau(n)\leq2\sqrt n$, which bounds the convolution and
restricted-divisor rows. Square indicators are bounded by one.
Finally the displayed Hurwitz expression for $T_9$ gives
$|T_9(rm)|\leq t(r)m$ for $m\geq1$.  Since
$t(r)=4r+8\sqrt{r/3}+2$, the two nonzero branches of $T_{27}$ give
\[
 |T_{27}(m)|\leq
 \begin{cases}
  t(1/3)m=6m,&m\equiv0\pmod3,\\
  \tfrac12t(3)m=11m,&m\equiv1\pmod3.
 \end{cases}
\]
The remaining branch vanishes, so $b(T_{27})=11$ suffices uniformly.
Each of these bounds has the form $b(f)m^{w}$ with an exponent
$w\leq3$, and $m^{w}\leq m^3$ for integers $m\geq1$; every row therefore
satisfies $|f(m)|\leq b(f)m^3$, which is how the single remainder
$C_3m^3$ absorbs all the sub-leading terms.
Apply the triangle inequality to the finite coefficient sums.
The first sufficient threshold follows from the negative leading-error
bound in \cref{cor:achiral-leading-error-sign}.

For the sharper assertion put $R(m)=a_m-W_5(m)-S_{9/2}(m)$.
Retain the exponents in the preceding estimates: the divisor terms of
weights four and three contribute powers $m^3,m^2$, the Cohen terms
contribute $m^{5/2},m^{3/2}$, and the reduced-form bound for $H(rm)$
contributes $m,\sqrt m,1$. The weight-two bound is
$m(1+\log m)/j$, and $1+\log m\leq\sqrt m$ for $m\geq16$;
indeed the difference is positive at $16$ and increasing thereafter.
The $T_9,T_{27}$ bounds are linear, and the divisor-pair and square
terms contribute $\sqrt m$ and $1$. Summing their absolute coefficients
in Appendix~A and rounding each grouped constant upward gives
\begin{equation}\label{eq:retained-power-majorant}
\begin{aligned}
 |R(m)|\leq P(m):={}&\frac{11}{100}m^3+\frac{63}{100}m^{5/2}
       +\frac32m^2+13m^{3/2}\\
       &+37m+40\sqrt m+27\qquad(m\geq16).
\end{aligned}
\end{equation}
The finite coefficient comparisons, including the square-root and
zeta bounds rounded upward, are recorded in
\path{data/public_numerical/counting_coefficients.json}.
Every monomial of $P(m)/m^{7/2}$ decreases for $m>0$, and
\[
 \frac{P(225)}{225^{7/2}}
 =\frac{1240093}{113906250}<\frac{11}{1000}
 <\frac{1279}{113400}
 =\frac{41}{5670}+\frac{17}{4200}
 <\frac{41}{5670}+\frac{\sqrt3}{420}.
\]
Thus the negative $S_{9/2}$ term dominates $R$ for every $m\geq225$.

To prove the sign down to $m=10$, mark an unoriented norm-minus-one
line perpendicular to a split vector of norm $m$ in $I_{1,9}$.
Such a line blocks every congruence state by the split obstruction
criterion of Paper~II; taking an integral multiple preserves the
perpendicular. Distinct such lines are orthogonal: in the negative
definite perpendicular their unit generators have integral inner
product of absolute value strictly less than one. Thus a blocked
vector receives at most nine marks. Orbit--stabilizer, applied in each
fiber of the forgetful map, bounds the split collision mass below by
one ninth of the marked mass, including stabilizer weights.

A unit line splits off integrally with complement $I_{1,8}$.
Its normal group has order two and all $183680$ states are blocked,
so the marked mass is $91840\,\mathfrak m_{I_{1,8}}(m)$.
For clarity, the odd-unimodular Cohen return
\inv{lem:uniform-odd-dyadic-return} at determinant one, oddity
$1-8\equiv1\bmod8$ and $k=4$ gives the shape
$H_4(4m)+16H_4(m)$; there is no bad odd prime. At norm one the
unimodular definite complements are $I_8,E_8$, each with its full
orthogonal stabilizer. Since $H_4(1)=1/120$ and
$H_4(4)=(1+2^7-2^3)/120$, the scalar multiplying this shape is
\[
 \frac{1/(2^8\,8!)+1/696729600}{137/120}=\frac1{11612160}.
\]
Using $H_4(4m)\ge(4m)^{7/2}/140$ from
\cref{cor:achiral-leading-error-sign} and $H_4(m)\ge0$, we obtain
\[
 M_{\rm split}^{\rm amb}(m)-M_{\rm split}^{\rm open}(m)
 \ \ge\ \frac{41}{51030}m^{7/2}>\frac{m^{7/2}}{1250}.
\]

Let $U$ be the sum of the three ambient masses of types $2,3,4$ and
the two open masses of types $10,12$. The seven-sector formula and
$M^{\rm open}\le M^{\rm amb}$ show that it suffices to bound
$U<m^{7/2}/1250$. Write $m=2^a b$ with $b$ odd. Use $S,T$ from
\rot{thm:siegel-q234-ambient}, the divisor sums with $\chi_{-3}$ on
the divisor and the complementary divisor, respectively. When $a=0$,
$36(M_2^{\rm amb}+M_4^{\rm amb})=S(b)+T(b)$. For $a\ge1$ it equals
\[
 \frac{4(-2)^a-7}{3}S(b)+\frac{4\cdot2^a-7(-1)^a}{3}T(b).
\]
Indeed, substitute $S(2^j)=(1-(-2)^{j+1})/3$,
$T(2^j)=(-1)^j S(2^j)$ in the two ambient formulas, with absent
dilations zero. Both coefficients have absolute value at most
$5\cdot2^a/2$. Since the divisors of $b$ are odd,
$\sigma_1(b)/b\le1+\frac12\log b$ by integral comparison.
For $m\ge16$, $1+\frac12\log m\le\frac23\sqrt m$
by comparison at $16$ and differentiation. Consequently
\[
 M_2^{\rm amb}+M_4^{\rm amb}\le\tfrac{5}{54}m^{3/2}.
\]
The type-three formula gives
$M_3^{\rm amb}\le\zeta(2)m^2/192<m^2/112$.
The two remaining open masses sum to at most
\[
 \tfrac25\{\tau(m/10)+\tau(m/5)\}+\tfrac13\tau(m/12)
 \le \left(\tfrac4{5\sqrt{10}}+\tfrac4{5\sqrt5}
                 +\tfrac2{3\sqrt{12}}\right)\sqrt m<\sqrt m .
\]
Dividing by $m^{7/2}$ gives a decreasing majorant, already at $16$
bounded by $5/(54\cdot256)+1/(112\cdot64)+1/4096<1/1250$.

For odd $m\ge11$ only types four and ten can contribute.
Here $1+\frac12\log m\le\sqrt m$, and hence
$U\le m^{3/2}/18+2\sqrt m/5$, since $2/\sqrt5<1$.
After division by $m^{3/2}$ it suffices that
$1/18+2/(5m)<m^2/1250$, which holds at $11$ and persists.
For the three small even arguments the same short divisor formulas
give the following upper bounds, retaining $d_{\rm same}(2)=0$:
\[
\begin{array}{c|ccc}
m&M_2^{\rm amb}+M_4^{\rm amb}&M_3^{\rm amb}&M_{10}^{\rm open}+M_{12}^{\rm open}\\ \hline
10&10/9&13/24&\le2/5\\
12&1/3&25/72&\le1/3\\
14&0&25/24&0
\end{array}
\]
Each row sums to less than $9/4<10^{7/2}/1250$. This proves the
strict inequality throughout the stated range.
\end{proof}

\section{Finite verification and reference audits}
\label{app:finite-verification}

The ancillary file
\path{data/public_numerical/counting_coefficients.json} exports the
$45+181$ coefficient entries as exact structured data. Appendix~A,
including its function and index conventions, is the mathematical
specification; the export and calculator must agree with it exactly.
The package's \texttt{SHA256SUMS} binds all three representations,
and its verifier rejects a disagreement.

The companion citations refer to coordinated release~1.4.1 at
\url{https://doi.org/10.5281/zenodo.24289727}. Each source
package includes \texttt{COMPANIONS.json}, which identifies the companion
source snapshots by SHA-256 and records the resolved statement labels.
The theorem numbers refer to those snapshots; later versions can be
matched by their stable source labels and content identifiers.
The comparisons in this appendix check the formulas independently on
finite ranges; the all-index argument is
\cref{thm:exact-all-index-grand-formula}.

\subsection{Reference-data audit}\label{sec:data-audit}

The reference tables of the ancillary package list the $191$ indices
$m=10,\dots,200$; all counts vanish below $m=10$.  For the \emph{oriented} counts, OEIS A057210 revision 23, accessed on
28 August 2026, agrees exactly on all 191 rows with Table~3 of the arXiv
version~3 \cite{engel2024preprint,oeisA057210}.  The executable Sage
reconstruction embedded in that source gives the same sequence.  Our audit
also identifies 28 smaller values in $m=68,\ldots,100$ as a superseded 2006
pre-bug-fix \texttt{fullgen} table; those values are retained only as a
provenance warning and are not used.

The \emph{ordinary} comparison values $\iso(2m)$ are taken from
A007894 \cite{oeisA007894,brinkmann2012}. They check the formula for
$m\leq200$. Completeness of the finite arithmetic classifications
is established from their stated domains, independently of these values.

\subsection{Cross-sector comparisons and implementation checks}\label{app:cross-sector-verification}

The rotary formulas and their normalizations are proved in
\rot{thm:rotary-rows}. In particular, the $Q_3$ return uses the genuine
ambient atoms rather than all abstract root planes
(\rot{thm:rotary-atom-radical}, \rot{cor:q3-analytical-open-return});
the $Q_4$ ternary mass is the class-number expression of
\rot{thm:q4-t9-class-number}. These distinctions are already incorporated
in the coefficient formula of \cref{cor:complete-coefficient-formula}.

The split coordinate inventories independently check the normal and
complement classification of \inv{thm:analytical-lower-complements}.
The comparison in \cref{app:finite-assembly} tests the assembly across
sectors.

\subsection{Exact finite assembly of the seven sectors}
\label{app:finite-assembly}
\begin{proposition}[Exact finite open-mass assembly]
\label{prop:finite-assembly}
For the seven masses in \cref{thm:exact-all-index-grand-formula},
the specialization of the counting identity at $m=30$ is
\begin{equation}\label{eq:finite-assembly}
\begin{split}
 a_{30}={}&M_{\rm split}^{\rm open}(30)+M_{\rm antipodal}^{\rm open}(30)
       +M_{2}^{\rm open}(30)+M_{3}^{\rm open}(30)+M_{4}^{\rm open}(30)\\
     &+M_{10}^{\rm open}(30)+M_{12}^{\rm open}(30).
\end{split}
\end{equation}
The seven summands at this index are respectively
\[
 \frac{251}{3},\quad\frac{38}{15},\quad\frac23,\quad4,
 \quad\frac13,\quad\frac45,\quad0,
\]
whose sum is $92=a_{30}$.

\end{proposition}

\begin{proof}
Substitute $m=30$ in the seven arithmetic expressions of
\cref{thm:exact-all-index-grand-formula}. Their numerators over the
common denominator $15$ are $1255,38,10,60,5,12,0$.
Their sum is $1380$, giving $a_{30}=92$.
\end{proof}
Independent lattice and graph enumerations agree with the seven-sector
sum throughout $0\leq m\leq200$. The comparison uses reciprocal
centralizer-stabilizer weights and the Burnside target
$2\iso(2m)-\BB_m$.

\subsection{Symmetry comparisons and finite-data cautions}
\label{app:symmetry-comparisons}

The independent point-group census checks each sector separately through
$C_{150}$. A graph of point group $G$ contributes
$r_j(G)/|G^+|$ to sector $j$, where $r_j(G)$ counts its reversing
elements of that type. This is orbit--stabilizer applied to
\cref{lem:witness-formula}; the complete weights and the $462$ comparisons
are given in \rot{tab:point-group-sector-weights} and
\rot{prop:rotary-point-group-census}.

Finite sizes can also mislead about the main term. At $m=200$,
$W_5(m)/a_m\simeq2.30$: the leading term is still more than twice the
count. The negative contribution in
\cref{cor:achiral-leading-error-sign} explains this slow approach.

\section{An exact prefix algorithm}
\label{app:exact-prefix-algorithm}

The arithmetic identities in this appendix evaluate the coefficient
formula at all indices through a prescribed bound. The class-number
relations are classical; the purpose here is to give the complete
division and bit-complexity argument used in
\cref{sec:counting-complexity}.

For these auxiliary power series, extend the positive-index Hurwitz
numbers of \cref{app:complete-coefficient-formula} by
$H(0)=-1/12$ and $H(n)=0$ for $n<0$, and put
\[
 \mathcal H(q)=\sum_{n\geq0}H(n)q^n,
 \qquad \theta(q)=\sum_{r\in\mathbb Z}q^{r^2},
 \qquad d_{\min}(n)=\sum_{d\mid n}\min(d,n/d).
\]
The auxiliary constant term does not change any positive-index counting
convention.

\begin{proposition}[The complete Hurwitz quotient]
\label{prop:complete-hurwitz-quotient}
Define $C(q)=\sum_{n\geq0}c_nq^n$ by
\[
\begin{aligned}
 c_0&=-1,\\
 c_{4n}&=24\sigma_1(n)-12d_{\min}(n) &&(n\geq1),\\
 c_n&=4\sigma_1(n)-6d_{\min}(n) &&(n\geq1\text{ odd}),\\
 c_n&=0 &&(n\equiv2\pmod4).
\end{aligned}
\]
Then $12\mathcal H(q)=C(q)/\theta(q)$ as formal power series.
\end{proposition}
\begin{proof}
The Kronecker--Hurwitz and Eichler relations, with the present
normalization, are
\[
 \sum_{r\in\mathbb Z}H(4n-r^2)
   =2\sigma_1(n)-d_{\min}(n) \quad(n\geq1),
\]
\[
 \sum_{r\in\mathbb Z}H(n-r^2)
   =\frac{\sigma_1(n)}3-\frac{d_{\min}(n)}2
   \quad(n\geq1\text{ odd});
\]
see \cite[\S5]{cohen1975} and
\cite[equations~(1.1)--(1.3)]{mertens2014classRelations}.
They evaluate the coefficients of $\theta\mathcal H$ at indices
divisible by four and at odd indices. At an index congruent to two
modulo four, subtracting a square leaves residue one or two, where
$H$ vanishes. The constant coefficient is $-1/12$. This exhausts all
indices. Since $\theta(0)=1$, division in the formal power-series
ring is unique and gives the result.
\end{proof}

\begin{samepage}
\begin{proposition}[Exact packed division]
\label{prop:exact-hurwitz-packed-division}
Let $N\geq0$, write $[F]_N$ for truncation through degree $N$, and
put $U=C+\theta$. Choose a power of two $R\geq2$ with $R>12N$.
For $Q=R^{N+1}$, the base-$R$ digits of the unique residue $Z\in[0,Q)$
defined by
\[
 Z\equiv[U]_N(R)[\theta]_N(R)^{-1}\pmod Q
\]
are $(0,12H(1),\ldots,12H(N))$.
\end{proposition}
\end{samepage}
\begin{proof}
The reduced forms defining $H(n)$, including the imprimitive forms,
satisfy $1\leq a\leq A=\lfloor\sqrt{n/3}\rfloor$ and $|b|\leq a$,
with $c$ fixed by $4ac-b^2=n$. Their weights are at most one.
For $A\geq1$ this gives
\[
 0\leq H(n)\leq\sum_{a=1}^A(2a+1)
   =A^2+2A\leq3A^2\leq n.
\]
The values at $n=1,2$ vanish, and $12H(n)$ is integral by the
class-number weights. Set $P=12\mathcal H+1$, so $P(0)=0$ and
$U=\theta P$. Evaluation at $q=R$ defines a homomorphism
\[
 \mathbb Z[q]/(q^{N+1})\longrightarrow\mathbb Z/(R^{N+1}),
\]
and $[\theta]_N(R)$ is odd, hence invertible modulo $Q$.
Every coefficient of $[P]_N$ lies in $[0,R)$, so $[P]_N(R)$ lies
in $[0,Q)$ and has precisely those coefficients as its digits.
The congruence and uniqueness of this representative prove the result.
At $N=0$ the residue is zero. Intermediate carries do not affect the
argument; the bound on the final digits is what identifies the return.
\end{proof}

\begin{theorem}[Bit complexity of the ordinary prefix]
\label{thm:ordinary-prefix-bit-complexity}
For every $M\geq1$, the list $(\iso(2m))_{1\leq m\leq M}$ can be
computed deterministically and exactly in
\[
 O\!\left(\mathsf{Mul}(O(M\log(M+2)))\log(M+2)
          +M\log^2(M+2)\right)
\]
bit operations, where $\mathsf{Mul}(L)$ is the cost of multiplying
two $L$-bit integers. Stored data occupy $O(M\log(M+2))$ bits in
addition to the multiplication workspace. The multiplication bound
$\mathsf{Mul}(L)=O(L\log L)$ of \cite[Thm.~1.1]{harveyHoeven2021}
gives $O(M\log^3(M+2))$ bit time.
All constants may depend on the fixed coefficient formula.
\end{theorem}
\begin{proof}
First construct $C$ through degree $N$ by walking the ordered pairs
$de\leq N$ in a divisor sieve. Add $d$ to the $\sigma_1(de)$ entry and
$\min(d,e)$ to the $d_{\min}(de)$ entry. There are
$O(N\log(N+2))$ visits. Both entries are bounded by
$N(1+\log N)$, so this takes $O(N\log^2(N+2))$ bit operations.
For fixed-width packing choose $R$ as the least power of two with
\[
 R>\max\left(12N,2,\max_{0\leq n\leq N}|[q^n]U|\right),
\]
with its exponent rounded up to a whole number of bytes.
The extra bounds ensure that the input coefficients fit into single
blocks; the preceding proposition only needs the bound on the output
digits.
Then $\log R=O(\log(N+2))$ and the packed integers have
$L=O(N\log(N+2))$ bits.

Use fixed-width buffers to pack the positive and negative parts of
$[U]_N$ separately, then subtract their integer values. Pack
$[\theta]_N$ similarly. For an odd integer $t$, its inverse modulo
$2^b$ is obtained by starting at $x=1\pmod2$ and doubling the known
precision with
\[
 x\longmapsto x(2-tx),\qquad 1-tx(2-tx)=(1-tx)^2.
\]
Reduce operands and results at each precision. There are $O(\log L)$
stages, all using integers of $O(L)$ bits. One final multiplication
and reduction give $Z$, and a single buffer export gives its digits.
Thus Hurwitz numbers through $N$ cost
$O(\mathsf{Mul}(O(L))\log L+N\log^2(N+2))$ bit operations.
Packing and unpacking require linear buffer work, rather than
successive shifts of a growing full-size integer.

The Cohen coefficients have equally explicit products. Set
$\mathcal H_k(q)=\zeta(1-2k)+\sum_{n\geq1}H_k(n)q^n$ for
$k=2,3,4$, and define
\[
 \mathsf E_2(q)=1-24\sum_{n\geq1}\sigma_1(n)q^n,
 \qquad \mathsf E_4(q)=1+240\sum_{n\geq1}\sigma_3(n)q^n,
\]
\[
 A(q)=\sum_{n\geq1}\left(\sum_{d\mid n}\chi_{-4}(n/d)d^2\right)q^n,
 \qquad
 B(q)=\sum_{n\geq1}\left(\sum_{d\mid n}\chi_{-4}(d)d^2\right)q^n.
\]
Cohen's identities \cite[Corollary~3.2, Propositions~4.1--4.2]{cohen1975}
give
\[
 \begin{aligned}
 120\mathcal H_2&=\theta\mathsf E_2(q^4)-6q\theta',\\
 252\mathcal H_3&=-\theta(1+2A-4B),\\
 240\mathcal H_4&=\theta\mathsf E_4(q^4).
 \end{aligned}
\]
Each divisor kernel has polynomially bounded integer coefficients and
can be sieved in $O(N\log^2(N+2))$ bit operations. Fixed powers of
the divisors can be precomputed by finite differences in $O(N)$
additions.

For completeness, exact multiplication also handles the signs in these
kernels. If $K(q)=\sum_{j=0}^Nk_jq^j$ is an integer kernel, let
$D=\max_j|k_j|$ and write $[\theta]_N=\sum_{j=0}^Nt_jq^j$,
$S=\sum_jt_j$. If $D=0$ the product is zero. Otherwise choose a
power of two $R>2DS$ and multiply
\[
 \sum_{j=0}^N(k_j+D)R^j
 \quad\hbox{by}\quad
 \sum_{j=0}^Nt_jR^j.
\]
Every coefficient of the corresponding polynomial product lies in
$[0,2DS]$, so the integer product has no carries between coefficient
digits. For $j\leq N$, subtract $D\sum_{a=0}^jt_a$ from digit $j$
to recover $[q^j]K\theta$. Since $D$ is polynomial in $N$ and
$S=O(\sqrt{N+1})$, each digit has $O(\log(N+2))$ bits. The three
Cohen tables therefore need a fixed number of $O(N\log(N+2))$-bit
products and linear buffer operations. The derivative term is supported
on squares and is added directly.

Finally use the coefficient expression in
\cref{app:complete-coefficient-formula}. Its literal dilations require
Hurwitz arguments at most $144M$, and Cohen arguments at most
$72M$, $36M$, $12M$ for $k=2,3,4$, respectively. The apparent rotary
exceptions are also covered: \cref{lem:counting-arithmetic-reduction}
writes $T_9$ as a fixed Hurwitz combination,
and $T_{27}$ uses that combination at $m/3$ or $3m$.
The remaining functions are fixed-character divisor sums, character
convolutions, congruence-restricted divisor sums, fixed dilations and
square indicators. Their sieves have the same
$O(M\log^2(M+2))$ bit bound; their powers and periodic characters
are precomputed. The residue and valuation projectors are evaluated
from a least-prime table. All values are polynomially bounded, and
the fixed rational coefficients have a common denominator independent
of $M$. The fixed linear combination at each index therefore costs
$O(\log(M+2))$ bit operations. Taking $N=144M$ in the constructions
above gives the stated time and storage bounds.
\end{proof}

\section*{Declarations}
\subsection*{Use of artificial intelligence}
During preparation of this work, the author used Fable 5.1 (Claude Code) and
Astra 6.0 (OpenAI Codex) and Deep Review (Google) to assist with exploration,
manuscript revision and verification code.  After using these tools, the
author reviewed and edited the content as needed and takes full
responsibility for the content of the publication.
\subsection*{Competing interests}
The author declares no conflict of interest.
\subsection*{Data and code availability}
The archive DOI for the fullerene counting paper series and its code and data packages is
\url{https://doi.org/10.5281/zenodo.24289727}.
The ancillary package of this article accompanies it as ancillary files.
It contains the source of the exact calculators and comparison checkers
named in the text, the machine-readable tables, and the finite inventories,
kernels and witnesses that the text calls ancillary, together with the
build instructions that rerun them.  Every reference in the text to the
ancillary package, to an ancillary calculator or table, or to an
independent finite calculation refers to this package.
The calculator \texttt{buckycount} is also available at
\url{https://github.com/krk/buckycount}.

\bibliographystyle{plain}
\bibliography{references}
\end{document}